\documentclass[11pt]{amsart}
\usepackage[T1]{fontenc}
\usepackage{amsmath,amssymb,amsthm,mathrsfs}
\usepackage{tikz-cd}                
\usepackage{float}                 
\usepackage[margin=1in]{geometry}
\usepackage[hidelinks]{hyperref}
\newtheorem{theorem}{Theorem}[section]
\newtheorem{proposition}[theorem]{Proposition}
\newtheorem{lemma}[theorem]{Lemma}
\newtheorem{corollary}[theorem]{Corollary}

\theoremstyle{definition}
\newtheorem{definition}[theorem]{Definition}

\theoremstyle{remark}

\newcommand{\Mq}{\mathsf q}                      
\newcommand{\Mt}{\mathsf t}                      
\newcommand{\K}{\mathbb K}                       
\newcommand{\bv}{{\mathbf v}}                    
\newcommand{\betaz}{\beta_0}                     

\newcommand{\La}{\Lambda}                        
\newcommand{\bl}{{\boldsymbol{\lambda}}}         
\newcommand{\bmu}{{\boldsymbol{\mu}}}
\newcommand{\Add}[1]{A(#1)}                      
\newcommand{\Rem}[1]{R(#1)}                      
\newcommand{\bwt}[1]{\zeta_{#1}}                 
\newcommand{\CY}{\mathcal Y}
\newcommand{\Yfn}[1]{\CY_{#1}}                   
\newcommand{\res}[2]{r_{#1}(#2)}                 
\newcommand{\cores}[2]{r^*_{#1}(#2)}             
\newcommand{\Cartan}[1]{H_{#1}}                  
\newcommand{\Newton}[1]{E_{#1}}                  

\newcommand{\cornermeas}[1]{\rho_{#1}}           
\newcommand{\mom}[2]{\mu_{#1}(#2)}               
\newcommand{\srcnorm}[1]{n_{#1}}                 
\newcommand{\Gram}[1]{G(#1)}                     
\newcommand{\ip}[2]{\left\langle #1,#2\right\rangle}
\newcommand{\ipinf}[2]{\left\langle #1,#2\right\rangle_\infty}   

\newcommand{\ceig}[1]{c_{#1}}                    
\newcommand{\bplus}[1]{{#1}^{+}}                 
\newcommand{\Remp}[1]{R^{+}(#1)}                 
\newcommand{\Wfib}[1]{W_{#1}}                    
\newcommand{\Wdeg}[2]{W_{#1,#2}}                 
\newcommand{\Fdeg}[2]{\mathcal F_{#1,#2}}         
\newcommand{\bx}[3]{(#1,#2)_{#3}}                

\newcommand{\Ptwo}{\mathbb P^2}
\newcommand{\skyO}{k_0}                          
\newcommand{\Kloc}[1]{K_T(#1)_{\mathrm{loc}}}    
\newcommand{\framing}{\mathbb C^m_{\bv}}         
\newcommand{\Gies}[2]{\mathcal M(#1,#2)}         
\newcommand{\Zc}{\mathfrak Z}
\newcommand{\Hecke}[1]{\Zc_{#1,#1+1}}            
\newcommand{\HeckeTwo}[1]{\Zc_{#1,#1+2}}         
\newcommand{\Parab}[3]{\mathcal M^{\mathrm{par}}(#1,#2,#3)}  
\newcommand{\Fock}[1]{\mathcal F_{#1}}            
\newcommand{\Sh}{F}                              
\newcommand{\Fsh}[1]{\Sh_{#1}}                   
\newcommand{\Wsp}[1]{W_{#1}}                     
\newcommand{\HH}{\mathcal H}
\newcommand{\Hsp}[1]{\HH_{#1}}                   
\newcommand{\Csec}[1]{\La^{\otimes #1}\otimes\K^{#1}}   
\newcommand{\Wc}{\mathcal W}                     
\newcommand{\Tc}{\mathcal T}                     
\newcommand{\Vone}{V_1}                          
\newcommand{\auxe}[1]{e_{#1}}                    

\newcommand{\kket}[1]{\lvert #1\rangle\!\rangle} 
\newcommand{\fix}[1]{[#1]}                       
\newcommand{\pt}[2]{[#1,#2]}                     

\newcommand{\pminus}{p^-}                        
\newcommand{\pplus}{p^+}                         
\newcommand{\taut}{\ell}                         
\newcommand{\dtaut}{\taut^{\vee}}                
\newcommand{\tautone}{\ell_1}
\newcommand{\tauttwo}{\ell_2}
\newcommand{\pulm}{p^{-*}}                       
\newcommand{\pushp}{p^+_*}                       
\newcommand{\pim}{\pi^-}                         
\newcommand{\pip}{\pi^+}                         
\newcommand{\pimpull}{\pi^{-*}}                  
\newcommand{\pippush}{\pi^+_*}                   
\newcommand{\pippull}{\pi^{+*}}                  
\newcommand{\incl}{j}                            
\newcommand{\Shift}{S}                           
\newcommand{\Omap}{\Omega}                       
\newcommand{\OmapG}{\Omega_{\mathcal G}}         
\newcommand{\Llet}{L}                            
\newcommand{\Gmod}[1]{\mathcal G_{#1}}           
\newcommand{\Trunc}[1]{\mathfrak t_{#1}}         

\newcommand{\zz}[1]{z_{#1}}
\newcommand{\yy}[1]{y_{#1}}
\newcommand{\zone}{\zz{1}}
\newcommand{\yone}{\yy{1}}
\newcommand{\dplus}{d_+}
\newcommand{\dminus}{d_-}
\newcommand{\ZNtwo}[1]{\mathfrak Z_{(#1)}}       

\newcommand{\xplus}[1]{x^+_{#1}}                 
\newcommand{\sigH}{\sigma_H}                     
\newcommand{\bu}{{\mathbf u}}                    
\newcommand{\pk}[2]{p_{#1}^{(#2)}}               
\newcommand{\GMP}[1]{P_{#1}}                     
\newcommand{\GMPn}[1]{\widetilde P_{#1}}         
\newcommand{\pvac}{\mathfrak p_\varnothing}      
\newcommand{\Nekbare}[2]{N_{#1,#2}}              
\newcommand{\Nek}[2]{\widetilde N_{#1,#2}}       
\newcommand{\Hiso}{\iota}                        
\newcommand{\kappam}[1]{\kappa_{#1}}             

\newcommand{\LL}{\mathbf L}
\newcommand{\Lax}[1]{\LL_{#1}}                   
\newcommand{\GLax}[1]{\mathcal L_{#1}}           
\newcommand{\QQ}{\mathsf Q}
\newcommand{\RaiseFull}[1]{\QQ_{#1}(z)}          
\newcommand{\Raise}[1]{\widehat\QQ_{#1}}         
\newcommand{\Akernel}[1]{\mathsf A_{#1}}         
\newcommand{\CC}{\mathsf C}
\newcommand{\Comp}[1]{\CC_{#1}}                  
\newcommand{\DD}{\mathcal D}
\newcommand{\Twist}[1]{\DD_{#1}}                 
\newcommand{\piz}{\pi_0}                         
\newcommand{\pile}{\pi_{\leq0}}                  
\newcommand{\TT}{\mathbb T}
\newcommand{\Transfer}[2]{\TT_{#1}(#2)}          
\newcommand{\TransferMat}[2]{T_{#1}(#2)}         
\newcommand{\Pser}[1]{\mathcal P_{#1}}           
\newcommand{\Dq}{\Delta_{q_3}}                   

\newcommand{\rhov}{\rho^{(0,m)}_{\bv}}           
\newcommand{\rhou}{\rho^{(m,0)}_{\bu}}           
\newcommand{\Miki}{\mathcal S}                   
\newcommand{\sigV}{\sigma_V}                     
\newcommand{\Xcur}{\mathcal X}                   
\newcommand{\Xk}[1]{\mathcal X^{(#1)}}           
\newcommand{\Xhat}[1]{\widehat{\mathcal X}{}^{(#1)}}   
\newcommand{\MM}{\mathsf M}
\newcommand{\Miura}[1]{\MM_{#1}}                 
\newcommand{\VV}{\mathsf V}
\newcommand{\Vconj}[1]{\VV_{#1}}                 
\newcommand{\Rlast}[1]{\mathcal R_{#1}}          
\newcommand{\Miuraop}{\mathfrak M}               
\newcommand{\Miuraproj}{\Miuraop_{\leq0}}         
\newcommand{\Meig}[2]{\varphi_{#1}^{#2}}         
\newcommand{\Shuf}{\mathcal A^{+}}               
\DeclareMathOperator{\Sym}{Sym}

\newcommand{\Qcr}{Q}                             
\newcommand{\Qan}{Q^\dagger}                     
\newcommand{\PsiNS}[2]{\Psi_{#1}^{#2}}           
\newcommand{\Plethmap}{\Upsilon}                 
\newcommand{\Supmap}{\Gamma}                     
\newcommand{\Exp}{\operatorname{Exp}}            
\newcommand{\Negut}{Neguț}

\DeclareMathOperator{\Res}{Res}

\title{Nazarov--Sklyanin Lax operators at higher rank}
\author{Ryan Mickler}
\email{ry.mickler@gmail.com}
\date{\today}

\begin{document}
\maketitle

\begin{abstract}
We construct a rank-\(m\) Nazarov--Sklyanin Lax operator for generalized Macdonald functions.
The operator acts on the horizontal Fock module of quantum toroidal \(\mathfrak{gl}_1\), extended by an
auxiliary polynomial space with \(m\) components. A geometric realization is given by the tensor product
with the dual tautological line on the equivariant \(K\)-theory of a Hecke correspondence studied by \Negut.
The eigenbasis is indexed by an \(m\)-tuple of partitions and an addable
corner, whose content is the eigenvalue.
We give an anti-normal factorization and prove that the Lax operator is self-adjoint
for the transported Euler pairing. With unit framing weight, the rank-one
specialization recovers the Nazarov--Sklyanin factorization.
We then show that the Lax operator is conjugate to the
companion matrix of a projected q-Miura operator of the deformed \(W\)-algebra of type
\(\mathfrak{gl}_m\).
This also gives an explicit presentation
of the \(K\)-theory of the Hecke correspondence. \end{abstract}

\clearpage
\tableofcontents
\clearpage

\section{Introduction}\label{sec:intro}

Nazarov and Sklyanin \cite{NS} constructed a Lax operator for the Macdonald symmetric functions,
acting on \(\La[z^{-1}]\), where \(\La\) is the ring of symmetric functions.
The variable \(z^{-1}\) has degree
one, and the operator preserves the sum of the symmetric-function degree and the auxiliary
degree. Its eigenvectors are indexed by a partition and an addable corner, with eigenvalue
the content of that corner. Compressing its resolvent to the constant polynomials in \(z^{-1}\)
gives commuting operators on \(\La\), diagonal in the Macdonald basis.

In the normalization of \cite{M26}, the Lax operator has the factorization
\begin{equation}\label{eq:rankoneNS}
 \LL=\Qan\,\Twist{\Mt/\Mq}\,\pile\Qcr.
\end{equation}
Here \(\Qcr\) is a multiplication series and \(\Qan\) a differential series,
\(\pile\) keeps the non-positive powers of \(z\), and \(\Twist{\Mt/\Mq}\) multiplies
\(z^{-k}\) by \((\Mt/\Mq)^k\). The formulas for \(\Qcr\) and \(\Qan\) are given in
Section~\ref{sec:rankone}. The operator \(\LL\) is self-adjoint for the pairing defined there.
The transfer operator, the resolvent \((1-u^{-1}\LL)^{-1}\) projected to \(\La\), acts on the integral form \(J_\lambda\) of the Macdonald function by
\(1-\Yfn{\lambda}(1/u)^{-1}\), where \(\Yfn{\lambda}\) is the Nekrasov \(Y\)-function attached
to the fixed partition \(\lambda\), defined in Section~\ref{sec:mainresults}.

The generalized Macdonald functions of Bourgine, Cassia and Stoyan \cite{BCS} are labelled by
\(m\)-tuples of partitions \(\bl\). Their commuting Hamiltonians suggest seeking a Lax
construction with eigenvectors indexed by \((\bl,a)\), where \(a\) is an addable corner, and
with a compressed resolvent governed by the rank-\(m\) Nekrasov \(Y\)-function.
Replacing \(\La\) by \(\La^{\otimes m}\) in the rank-one auxiliary space gives too few states.
In each total degree the number of these corner pairs is \(m\) times the dimension of
\(\La^{\otimes m}[z^{-1}]\). We therefore need \(m\) auxiliary components, whose geometric
realization is established by the freeness theorem below. The compressed resolvent becomes
an \(m\times m\) matrix, and the \(Y\)-function appears in a determinant identity.
The problem is to realize this spectrum by an explicit horizontal operator,
obtain its anti-normal factorization, and find a conjugation that expresses
the matrix entries directly in terms of horizontal vertex currents.

The corner pairs are the torus-fixed points of the Hecke correspondence of
\Negut{} \cite{N1,N2}, which parametrizes framed sheaves of rank \(m\) on
\(\Ptwo\) with a subsheaf of colength one at the origin.
At the fixed point \((\bl,a)\), the dual tautological line has weight \(\ceig{a}\).
Tensoring with this line gives the geometric Lax operator of
Definition~\ref{def:operator}. The fixed-point classes form an eigenbasis with
the required corner-content spectrum and are orthogonal for a twisted Euler pairing.

The Fock model of Theorem~\ref{thm:fock} expresses the geometric operator on
the horizontal module. In these coordinates, the rank-\(m\) operator \(\Lax{m}\)
has the anti-normal factorization of Theorem~\ref{thm:factorization}.
An explicit graded change of coordinates conjugates \(\Lax{m}\) to \(\Miura{m}\), the companion
matrix of the projected q-Miura operator of the deformed \(W\)-algebra of type
\(\mathfrak{gl}_m\) (Theorem~\ref{thm:companion}).

Theorem~\ref{thm:kernel} constructs the eigenvectors of \(\Miura{m}\) by a recursion
over removable corners. Their last coordinates form a basis of the projected
q-Miura kernel at the inverse contents of addable corners. The conjugacy also
determines the Hecke \(K\)-theory \(\Wsp{m}\) as a module over the free algebra
on the dual tautological line and a shift
(Corollaries~\ref{cor:structure} and~\ref{cor:presentation}).
Theorem~\ref{thm:transfer} gives the transfer determinant.

\subsection{Main results}\label{sec:mainresults}

Throughout, \(q_1,q_2,v_1,\dots,v_m\) are the equivariant parameters of
\cite{BCS}, with \(q_1=\Mt^{-1}\), \(q_2=\Mq\), \(q_3=(q_1q_2)^{-1}=\Mt/\Mq\)
and \(\gamma=q_3^{1/2}\). We work over
\(\K:=\mathbb Q(q_1,q_2,v_1,\dots,v_m)(\gamma^{1/2})\), which includes the
normalizing powers of \(\gamma\) in the representations of \cite{BCS}.
We assume that \(q_1\) and \(q_2\) are generic and that
\(v_\alpha/v_\beta\notin q_1^{\mathbb Z}q_2^{\mathbb Z}\) for \(\alpha\neq\beta\).
We extend scalars to \(\K\) throughout.
An \(m\)-tuple of partitions \(\bl\) has addable and removable corners
\(\Add{\bl}\) and \(\Rem{\bl}\). We call the component index \(\alpha\) a colour.
Rows and columns are indexed from \(0\), so a
box \(b=\bx{i}{j}{\alpha}\) in row \(i\) and column \(j\) of colour \(\alpha\) has
weight \(\bwt{b}=v_\alpha q_1^{i}q_2^{j}\) and content
\(\ceig{b}=\bwt{b}^{-1}=v_\alpha^{-1}\Mt^i\Mq^{-j}\). In rank one with \(v_1=1\)
this is the content \([s]=\Mt^i\Mq^{-j}\) of a box \(s=(i,j)\) in \cite{M26}.
(In \cite{BCS} rows and columns are indexed from \(1\), so their box \((i,j)\) in the \(\alpha\)-th partition
is our \(\bx{i-1}{j-1}{\alpha}\).) As in \cite{M26}, \(\bplus{b}\) is the box
\(b+(1,1)\) of the same colour, and the outer corners of \(\bl\) are
\(\Remp{\bl}=\{\bplus{r}:r\in\Rem{\bl}\}\). Thus \(\ceig{\bplus{r}}=q_3\ceig{r}\).

\paragraph{The Hecke correspondence.}
We use the framed moduli spaces and Hecke correspondences described by
\Negut{} \cite[Sections~2 and~3.1]{N1}. Let \(\Gies{m}{n}\) be the moduli space of framed
torsion-free sheaves on \(\Ptwo\) of rank \(m\) and second Chern class \(n\).

The torus-fixed points of \(\Gies{m}{n}\) are indexed by the \(m\)-tuples \(\bl\) with \(|\bl|=n\). Put
\[
 \Fock{m}:=\bigoplus_n\Kloc{\Gies{m}{n}}.
\]
If \(i_\bl\) is the inclusion of the fixed
point represented by the sheaf \(\Fsh{\bl}\), write
\(\fix{\bl}:=(i_\bl)_*1=[\mathcal O_{\{\Fsh{\bl}\}}]\) for its structure-sheaf class.
By \cite{N1,BCS} and Proposition~\ref{prop:vertical},
\(\Fock{m}\) is the vertical module of \cite{BCS} of level \((0,m)\), and
\(\fix{\bl}\) is its vertical basis vector. We work in the horizontal module
instead, through Proposition~\ref{prop:horizontal}.

\begin{definition}[The Hecke correspondence]\label{def:space}
Let \(\Hecke{n}\) be the Hecke correspondence of \cite{N1} with support at the
origin, consisting of pairs \(\Sh'\subset\Sh\) with \(\Sh\in\Gies{m}{n}\),
\(\Sh'\in\Gies{m}{n+1}\) and \(\Sh/\Sh'\cong\skyO\), where \(\skyO\) is the
length-one skyscraper sheaf at the origin.
The correspondence has maps
\(\pminus(\Sh'\subset\Sh):=\Sh\) and \(\pplus(\Sh'\subset\Sh):=\Sh'\).
The tautological line bundle \(\taut\) has fibre \(H^0(\Ptwo,\Sh/\Sh')\) at
\(\Sh'\subset\Sh\), and its dual \(\dtaut\) has class \(\taut^{-1}\) in \(K\)-theory. Put
\[
 \Wsp{m}:=\bigoplus_n\Kloc{\Hecke{n}},
\]
and write \(\Wdeg{m}{N}:=\Kloc{\Hecke{N}}\) for its summand of degree \(N\).
\end{definition}

Let \(\HeckeTwo{N}\) be the space of two-step flags \(\Sh''\subset\Sh'\subset\Sh\)
with \(\Sh\in\Gies{m}{N}\) and both quotients isomorphic to \(\skyO\) \cite{N1,N2}, with maps
\(\pim(\Sh''\subset\Sh'\subset\Sh):=(\Sh'\subset\Sh)\) and
\(\pip(\Sh''\subset\Sh'\subset\Sh):=(\Sh''\subset\Sh')\).
The quotient lines \(\tautone\) and \(\tauttwo\) have fibres
\(H^0(\Ptwo,\Sh/\Sh')\) and \(H^0(\Ptwo,\Sh'/\Sh'')\), respectively.
The associated push-pull shift
\(\Shift\colon\Wdeg{m}{N}\to\Wdeg{m}{N+1}\) is defined in
Definition~\ref{def:shift}. Tensor product with \(\dtaut\) preserves the degree.
We write \(\taut^{-j}\) for the class
of \((\dtaut)^{\otimes j}\).

\begin{figure}[H]
\[
\begin{tikzcd}[column sep=2.4em, row sep=3.6em]
 & \Kloc{\Zc_{N-1,N+1}} \arrow[dr, "\pippush"] & & \Kloc{\HeckeTwo{N}} \arrow[dr, "\pippush"] & \\
 \Wdeg{m}{N-1} \arrow[ur, "\pimpull"] \arrow[rr, "\Shift"] \arrow[dr, shift left=0.7ex, "\pushp"]
 & & \Wdeg{m}{N} \arrow[ur, "\pimpull"] \arrow[rr, "\Shift"] \arrow[dr, shift left=0.7ex, "\pushp"]
     \arrow[dl, shift left=0.7ex, "\pminus_*"] \arrow["\dtaut", loop, out=122, in=58, distance=3.2em]
 & & \Wdeg{m}{N+1} \arrow[dl, shift left=0.7ex, "\pminus_*"] \\
 & \Fdeg{m}{N} \arrow[ur, shift left=0.7ex, "\pulm"] \arrow[ul, shift left=0.7ex, "\pplus{}^*"]
 & & \Fdeg{m}{N+1} \arrow[ur, shift left=0.7ex, "\pulm"] \arrow[ul, shift left=0.7ex, "\pplus{}^*"] &
\end{tikzcd}
\]
\caption{The Hecke correspondences in consecutive degrees, with
\(\Fdeg{m}{N}:=\Kloc{\Gies{m}{N}}\). The loop is tensor product with \(\dtaut\).
The paths through \(\Fdeg{m}{N}\) and through two-step flags differ by the shift,
\(\pulm\pushp=\pippush\pimpull+\Shift\) (Lemma~\ref{lem:excess}).}\label{fig:Wm}
\end{figure}

\begin{lemma}[{see Jiang \cite{Jiang} and Jiang--Li--Zhao \cite{JLZ}}]\label{lem:generation}
\begin{enumerate}
\item[(i)] As a graded \(\K[\Shift]\)-module, \(\Wsp{m}\) is free, with decomposition
\[
 \Wsp{m}=\bigoplus_{k\geq0}\Shift^k\bigoplus_{j=0}^{m-1}\taut^{-j}\,\pulm\Fock{m}.
\]
Its generator space is \(\bigoplus_{j=0}^{m-1}\taut^{-j}\pulm\Fock{m}\).
\item[(ii)] The operators \(\dtaut\) and \(\Shift\) do not commute:
\[
 \dtaut\,\Shift-q_3\,\Shift\,\dtaut=\dtaut\,\pulm\pushp.
\]
\end{enumerate}
\end{lemma}

Section~\ref{sec:fixedpoints} describes \(\pminus\) as a projectivization, and
Section~\ref{sec:freeness} proves part~(i) directly by fixed-point localization.

Iterating part~(ii) expresses \(\dtaut\,\Shift^k\) as \(q_3^k\Shift^k\dtaut\) plus terms
with fewer factors of \(\Shift\), through \(\pulm\pushp\) \eqref{eq:iterated}.
This relation determines the action of the geometric Lax operator in polynomial coordinates
and is used to construct the companion recurrence in Theorem~\ref{thm:companion}.

\Negut{} \cite[(3.35), (3.39)]{N3} chooses a twisted Euler pairing on the framed moduli
spaces for which the geometric raising and lowering correspondences satisfy
the adjoint identities of Lemma~\ref{lem:negut}.
The determinant twist adjusts Euler adjunction to the normalization of these two
correspondences. We pull the twist back along \(\pminus\) to the Hecke correspondence
and rescale each degree to recover the Nazarov--Sklyanin pairing when \(m=1\) and \(v_1=1\).
The compatibility of the Hecke pairing with the shift,
\(\ip{\Shift A}{\Shift B}=\Mq\,\ip{A}{B}\), is then a consequence
(Lemma~\ref{lem:shiftip}). The one-box adjunction is explained in Section~\ref{sec:form}.
On \(\Wdeg{m}{n}\) define the twisted Euler pairing
\begin{equation}\label{eq:eulerpairing}
 \ip{A}{B}:=(-1)^{mn}\Mq^{n}\,\chi_{\Hecke{n}}\bigl(A\otimes B\otimes
 \pulm\det(q_3\,\framing\otimes\Tc^\vee)\bigr),
\end{equation}
where \(\Tc\) is the tautological bundle of \(\Gies{m}{n}\) and \(\framing\) is the framing. Summands
of different degrees are orthogonal. All pairings below are symmetric and \(\K\)-bilinear,
and are nondegenerate at generic parameters.

\begin{definition}[The geometric Lax operator]\label{def:operator}
The geometric Lax operator is tensor product with the dual tautological line,
\[
 \GLax{m}:=\dtaut\otimes(-)\colon\Wsp{m}\to\Wsp{m}.
\]
\end{definition}

The operator \(\GLax{m}\) is self-adjoint, since
\(\ip{\dtaut A}{B}\) and \(\ip{A}{\dtaut B}\) are the same Euler characteristic in
\eqref{eq:eulerpairing}. Figure~\ref{fig:Wm}
shows \(\GLax{m}\) with the maps into and out of \(\Wsp{m}\).

\begin{proposition}
\label{prop:fixed}
The torus-fixed points of \(\Hecke{n}\) are the flags
\(\Fsh{\bl+a}\subset\Fsh{\bl}\), with \(|\bl|=n\) and \(a\in\Add{\bl}\).
Their structure-sheaf classes \(\pt{\bl}{a}\) form an eigenbasis of \(\GLax{m}\), with
\[
 \GLax{m}\,\pt{\bl}{a}=\ceig{a}\,\pt{\bl}{a}.
\]
\end{proposition}

\begin{proof}
The fixed framed sheaves are sums of monomial ideals indexed by the coloured partitions
\cite[Sections~5.2--5.3]{N1}. A fixed subsheaf of colength one is obtained by adding one
box \(a\), and the quotient line has weight \(\bwt{a}\).
Thus its dual has weight \(\ceig{a}=\bwt{a}^{-1}\). The fixed-point classes form a basis
after localization, and tensor product with the dual line acts by this weight.
\end{proof}

The rank-\(m\) Nekrasov \(Y\)-function of \cite{BCS} is
\[
 \Yfn{\bl}(w):=\frac{\prod_{a\in\Add{\bl}}(1-\bwt{a}/w)}
                  {\prod_{b\in\Remp{\bl}}(1-\bwt{b}/w)}.
\]
We use the residue coefficients \(\res{\bl}{a}\) and \(\cores{\bl}{r}\) of
\cite{BCS}, associated with \(\Yfn{\bl}\) and its reciprocal, with the
conventions specified in Lemma~\ref{lem:residues}.

For each multipartition \(\bl\), let \(\srcnorm{\bl}\) be
\((\Mq\gamma^m)^{|\bl|}\) times its vertical norm in \cite[(2.25)]{BCS}.
Proposition~\ref{prop:euler} identifies this with the geometric source norm
and expresses the Hecke norm as \(\srcnorm{\bl}/\res{\bl}{a}\).

\paragraph{The Fock model.}

The horizontal Fock space consists of symmetric functions in \(m\) independent
alphabets. In power-sum coordinates,
\[
 \La^{\otimes m}=\K[\pk{k}{\alpha}\mid k\geq1,\ 1\leq\alpha\leq m].
\]
The Heisenberg action uses multiplication by power sums and differentiation
with respect to them, and its vertex operators are formal series of such operators.

Let \(\mathcal E\) be the quantum toroidal algebra of \(\mathfrak{gl}_1\) in
\cite[Section~2.1]{BCS}, with currents \(x^\pm(w)=\sum_{k\in\mathbb Z}x^\pm_k w^{-k}\).
Its vertical representation \(\rhov\), of level \((0,m)\), acts on \(\Fock{m}\).
For horizontal weights \(\bu=(u_1,\dots,u_m)\), its representation \(\rhou\),
of level \((m,0)\), acts on \(\La^{\otimes m}\). The relation between these two
actions uses the Miki automorphism (Section~\ref{sec:bcs}).

For \(u_\alpha:=-\gamma v_\alpha\), let
\(\Hiso\colon\Fock{m}\to\La^{\otimes m}\) be the graded linear isomorphism of
Bourgine--Cassia--Stoyan \cite[Section~3.3]{BCS}, with the normalization of
Proposition~\ref{prop:horizontal}. The images
\(\kket{\bl}:=\Hiso\fix{\bl}\) form a basis of normalized generalized
Macdonald functions.

\begin{definition}[The auxiliary Fock space]\label{def:auxiliary}
Let \(V:=\La^{\otimes m}[z^{-1}]\) be the auxiliary polynomial space.
Set \(\Hsp{m}:=V\otimes\K^m\), where
\(\auxe{0},\dots,\auxe{m-1}\) is the standard basis of \(\K^m\).
Give \(\pk{k}{\alpha}\) degree \(k\), \(z^{-1}\) degree one, and each \(\auxe{j}\) degree zero.
The vectors \(z^{-k}\kket{\bl}\otimes\auxe{j}\) have total degree \(|\bl|+k\).
Let \(\piz\) be extraction of the constant term in \(z^{-1}\), so that
\(\piz\Hsp{m}=\Csec{m}\). We call this subspace the constant sector.
\end{definition}

In the isomorphism below, \(\auxe{j}\) represents the generator \(\taut^{-j}\) in
Lemma~\ref{lem:generation}(i). The framing scalar is
\(\betaz:=(-1)^{m+1}\prod_{\alpha=1}^m v_\alpha^{-1}\).

\begin{theorem}[The Fock model of \(\Wsp{m}\)]\label{thm:fock}
\begin{enumerate}
\item[(i)] The map
\[
 \Omap\colon\Hsp{m}\to\Wsp{m},\qquad
 z^{-k}\kket{\bl}\otimes\auxe{j}\longmapsto\Shift^k\,\taut^{-j}\,\pulm\fix{\bl},
\]
is an isomorphism of graded \(\K[\Shift]\)-modules, with \(\Shift\) acting on \(\Hsp{m}\)
by multiplication by \(z^{-1}\). The image of the constant sector under \(\Omap\) is
\[
 \Omap(\Csec{m})=\bigoplus_{j<m}\taut^{-j}\pulm\Fock{m}.
\]
\item[(ii)] The pairing
\[
 \begin{aligned}
 &\ip{z^{-k}\kket{\bl}\otimes\auxe{i}}{z^{-l}\kket{\bmu}\otimes\auxe{j}}\\
 &\qquad=\delta_{kl}\,\delta_{\bl\bmu}\,\Mq^k\,\srcnorm{\bl}
       \sum_{a\in\Add{\bl}}\res{\bl}{a}\ceig{a}^{i+j}
 \end{aligned}
\]
is the pullback under \(\Omap\) of the Euler pairing \eqref{eq:eulerpairing}.
The pairing is nondegenerate at generic parameters.
\end{enumerate}
\end{theorem}

\begin{definition}[The Nazarov--Sklyanin Lax operator]\label{def:NSLax}
The Nazarov--Sklyanin Lax operator is the geometric Lax operator in the Fock model of \(\Wsp{m}\),
\[
 \Lax{m}:=\Omap^{-1}\circ\GLax{m}\circ\Omap\colon\Hsp{m}\to\Hsp{m}.
\]
\end{definition}

The first main result constructs \(\Lax{m}\) from a raising matrix and its adjoint,
extending the rank-one factorization \eqref{eq:rankoneNS}.
The finite auxiliary matrix in this factorization comes from the Nekrasov
\(Y\)-function. The coefficients of \(\Comp{m}\) are expressed through the elementary Newton polynomials
\begin{equation}\label{eq:newtoncorners}
 \sum_{r\geq0}\Newton{r}(\bl)\,u^r
 :=\frac{\prod_{a\in\Add{\bl}}(1+u\,\ceig{a})}
          {\prod_{b\in\Remp{\bl}}(1+u\,\ceig{b})}
 =\frac{u^m}{\prod_\alpha v_\alpha}\,\Yfn{\bl}(-u).
\end{equation}
Thus the same \(Y\)-function gives the Newton coefficients through its normalized
expansion at the origin. The factor \(u^m\) removes its pole there.

\begin{theorem}[The Nazarov--Sklyanin anti-normal form]\label{thm:factorization}
On \(\Hsp{m}\),
\[
 \boxed{\;\Lax{m}=\Raise{m}^\dagger\,\Comp{m}\,\Twist{q_3}\,\Raise{m}\;}
\]
Here \(\Twist{q_3}\) multiplies \(z^{-k}\) by \(q_3^k\), and
\(\Raise{m}=\pile\circ\RaiseFull{m}\), with \(\RaiseFull{m}\) an \(m\times m\)
matrix of formal power series in \(z\) with operator coefficients on \(\La^{\otimes m}\),
built from a Heisenberg exponential, a projected current and a finite difference recursion
(Theorem~\ref{thm:raising}). \(\Raise{m}^\dagger\) is its adjoint for the pairing
of Theorem~\ref{thm:fock}, explicit through the horizontal anti-involution of
\cite{BCS} (Theorem~\ref{thm:adjointexplicit}).
On \(\kket{\bl}\otimes\K^m\), \(\Comp{m}\) is the companion matrix of the degree-\(m\)
polynomial extracted from the Nekrasov \(Y\)-function by
\begin{equation}\label{eq:companionpolynomial}
 \begin{aligned}
 p_\bl(x)&:=-\betaz\bigl(1+[\Yfn{\bl}(1/x)]_{>0}\bigr)\\
 &=x^m+\sum_{k=1}^{m-1}(-1)^k\Newton{k}(\bl)\,x^{m-k}-\betaz.
 \end{aligned}
\end{equation}
Here \([\,\cdot\,]_{>0}\) retains the strictly positive powers of \(x\) in the Laurent
expansion at infinity, and the constant term of \(p_\bl\) is fixed to \(-\betaz\).
Thus \(\Comp{m}\) represents multiplication by \(x\) on \(\K[x]/(p_\bl(x))\)
in the basis \(1,x,\dots,x^{m-1}\). It acts coefficientwise in \(z^{-1}\),
with entries \eqref{eq:companionentries}.
\(\Comp{m}\) and \(\Twist{q_3}\) are commuting and self-adjoint. For \(m=1\) and \(v_1=1\) this
is the Nazarov--Sklyanin operator \(\Qan\Twist{\Mt/\Mq}\pile\Qcr\) of \cite{M26}.
\end{theorem}

We call this the Nazarov--Sklyanin \emph{anti-normal} form of \(\Lax{m}\): the factor \(\Raise{m}^\dagger\),
which removes boxes, stands to the left of \(\Raise{m}\), which adds them, the opposite of normal
ordering.
The rank-two factors and their power-sum formulas are displayed in
Section~\ref{sec:example}, equations~\eqref{eq:ranktwohalves}--\eqref{eq:ranktwolax}.

\paragraph{The q-Miura form.}

In the rank-one NS form \eqref{eq:rankoneNS}, the dilation and projection
separate the raising half-vertex operator \(\Qcr\) from the lowering
half-vertex operator \(\Qan\). Conjugation by \(\Qan\) gives
\begin{equation}\label{eq:rankoneconjugacy}
 (\Qan)^{-1}\LL\Qan
 =\Twist{q_3}\pile\Qcr(z)\Qan
 =\Twist{q_3}\pile\eta^-(q_2\gamma z).
\end{equation}
The two half-vertex operators now occur as the normal-ordered product
defining a single horizontal current \cite[(2.34)]{BCS}.
This observation motivates our higher-rank construction: we seek a
conjugation of the NS operator whose matrix entries have explicit
expressions in the colour currents \(\eta^-_\alpha\).

On the tensor product of \(m\) Fock spaces, the fused horizontal currents
realize the deformed \(W\)-algebra of type \(\mathfrak{gl}_m\)
\cite{FR,FF,AKOS}. These currents provide the coefficients of the q-Miura
operator. Following \Negut{} \cite[Section~5.1]{N3},
we write the q-Miura operator as a formal difference operator, using the
horizontal normalization of \cite{BCS}. Write \(\Xk{k}(w)\) for the fused horizontal currents,
whose formulas are given in Section~\ref{sec:qmiura}. For a formal parameter \(y\), set
\[
 \Miuraop(y):=1+\sum_{k=1}^m y^k\gamma^k\Twist{q_3}^k\Xk{k}(q_2\gamma^mz).
\]
The dilation acts by \(\Twist{q_3}f(z)=f(z/q_3)\), and the factors in each summand
are composed in the displayed order. The action of \(\Miuraop(y)\) on auxiliary polynomials is
\[
 \Miuraop(y)\colon V
 \longrightarrow\La^{\otimes m}((z))[y].
\]
Here \(((z))\) denotes formal Laurent series with finitely many negative powers.
The annihilation series truncate on every polynomial input, while the creation
exponentials can produce infinitely many positive powers of \(z\).
The projection \(\pile\) discards those powers, giving
\begin{gather*}
 \Miuraproj(y):=\pile\Miuraop(y)=1-\sum_{k=1}^m y^k\Xhat{k}
 \quad\text{in }\operatorname{End}_{\K}(V)[y],\\
 \Xhat{k}:=-\gamma^k\Twist{q_3}^k\pile\Xk{k}(q_2\gamma^mz).
\end{gather*}
Let \(\Miura{m}\) be the companion matrix of \(\Miuraproj(y)\), acting on \(V^m\) by
\[
 \Miura{m}(\xi_1,\dots,\xi_m):=\Bigl(\sum_{k=1}^m\Xhat{k}\xi_k,\,\xi_1,\,\dots,\,\xi_{m-1}\Bigr).
\]

For each generalized Macdonald basis state \(\kket{\bl}\), the operator
\(\Comp{m}\) preserves \(\kket{\bl}\otimes\K^m\).
On this subspace, \(\Comp{m}\) is represented by the scalar companion matrix of \(p_\bl\).
The entries of \(\Miura{m}\) act on
the full auxiliary polynomial space and can change the Fock and auxiliary
degrees while preserving their sum. This second companion matrix represents
the full Lax operator in the coordinates of Theorem~\ref{thm:companion}.

\begin{theorem}[The q-Miura form]\label{thm:companion}
There is an explicit graded \(\K\)-linear isomorphism
\[
 \Vconj{m}\colon V^m\longrightarrow\Hsp{m}.
\]
It is triangular in the target order \(\auxe{1},\dots,\auxe{m-1},\auxe{0}\) and satisfies
\[
 \boxed{\;\Lax{m}\,\Vconj{m}=\Vconj{m}\,\Miura{m}.\;}
\]
For \(m=1\) and \(v_1=1\), \(\Vconj{1}=\Qan\) and
\(\Miura{1}\) is the conjugated rank-one operator in \eqref{eq:rankoneconjugacy}.
\end{theorem}

We call \(\Lax{m}=\Vconj{m}\Miura{m}\Vconj{m}^{-1}\) the q-Miura form of \(\Lax{m}\). The two forms
are related by \Negut's \cite{N3} factorization (Theorem~\ref{thm:negutfactorization}),
which gives Proposition~\ref{prop:WH} in our conventions.
The rank-two currents, companion matrix, and conjugating map are explicit in
Section~\ref{sec:example}, equations~\eqref{eq:ranktwovertices},
\eqref{eq:ranktwomiuracoefficients}, and~\eqref{eq:ranktwoconjugator}.

\paragraph{Eigenvectors and the kernel.}
On total degree \(N\), the companion matrix \(\Miura{m}\) has an eigenbasis
indexed by the pairs \((\bl,a)\) with \(|\bl|=N\) and \(a\in\Add{\bl}\).
Theorem~\ref{thm:kernel} constructs these vectors recursively. For each
nonzero \(c\), the last coordinates of the eigenvectors with eigenvalue \(c\) form a basis of
\(\ker\Miuraproj(c^{-1})\). The kernel basis vector indexed by \((\bl,a)\)
is uniquely determined by its constant term \(\betaz^{-1}\kket{\bl}\).

Corollary~\ref{cor:corneroperator} constructs this kernel basis using the inverse
of the lowering factor in the factorization of \Negut{} \cite[Theorem~1.2]{N3}.

\paragraph{The structure of \(\Wsp{m}\).}
The operators \(\dtaut\) and \(\Shift\) make \(\Wsp{m}\) a left module over the free algebra
\(\K\langle\Llet,\Shift\rangle\), with \(\Llet\) acting by \(\dtaut\) and \(\Shift\) by the shift.
Theorem~\ref{thm:companion} determines its underlying \(\K[\Llet]\)-module by a monic operator
polynomial of degree \(m\). Set
\begin{equation}\label{eq:modulequotient}
 \Gmod{m}:=V[\Llet]\big/\Llet^m\Miuraproj(\Llet^{-1})\,V[\Llet],
\end{equation}
the polynomials in a variable \(\Llet\) with coefficients in \(V\), modulo the image of
\(\Llet^m\Miuraproj(\Llet^{-1})=\Llet^m-\sum_{k=1}^m\Llet^{m-k}\Xhat{k}\), whose coefficients act on
\(V\) first. Corollary~\ref{cor:structure} identifies \(\Wsp{m}\) with this
quotient and determines the shift action in companion coordinates.
The full \(\K\langle\Llet,\Shift\rangle\)-module presentation is
Corollary~\ref{cor:presentation}.

\begin{theorem}[The transfer determinant]\label{thm:transfer}
The transfer operator
\[
 \Transfer{m}{u}:=\piz\,(1-u^{-1}\Lax{m})^{-1}\big|_{\Csec{m}}
\]
preserves each \(\kket{\bl}\otimes\K^m\). In the basis
\(\{\kket{\bl}\otimes\auxe{j}\}_{j=0}^{m-1}\) of this subspace, its matrix
\(\TransferMat{\bl}{u}\) satisfies
\[
 \boxed{\;\det\bigl(I_m-\TransferMat{\bl}{u}\bigr)=\Yfn{\bl}(1/u)^{-1},\;}
\]
and its \(00\) entry is \(1-\Yfn{\bl}(1/u)^{-1}\).
\end{theorem}

\subsection{Relation to the literature}
The operators \(\zone\) of \cite{CGM,GGS}, with their eigenbases and spectra, are known.
Proposition~\ref{prop:GGS} identifies the geometric Lax operator with the GGS operator
\(\zone^{-1}\). In rank one, Theorem~\ref{thm:rankone} identifies the Lax operator
with the CGM operator \(\zone\), using a plethystic change of variables and the
different parameter convention specified in Section~\ref{sec:cgm}.
Under the isometry of Section~\ref{sec:superspace}, the rank-one eigenvectors
correspond to Macdonald polynomials in superspace of fermionic degree one,
labelled by superpartitions. We use the
Hecke correspondences and two-step flags of \cite{N1,N2} and the representations, generalized Macdonald
functions and Miki automorphism of \cite{BCS}. The q-Miura operator and its fused currents are those of
the deformed \(W\)-algebras \cite{FR,FF,AKOS,AFHKSY,FHSSY}, and its factorization into factors that add
boxes, preserve them and remove them is \cite[Theorem~1.2]{N3}.

Jiang \cite[Corollary~8.15]{Jiang} proves a projectivization formula.
The iterated decomposition of Lemma~\ref{lem:generation}(i) is in
\cite[Proposition~7.7]{JLZ} for moduli of sheaves of any rank on a smooth projective surface, and part~(ii)
is the excess formula of \cite[Propositions~2.27 and~3.4]{N2}. The identity \cite[Lemma~6.2]{Nsurf}, in
the \(K\)-theory of the fibre product of the Hecke correspondence with itself over the moduli space, is
used there to prove that the \(W\)-currents vanish above degree \(m\). Read as an identity of operators on
\(\Wsp{m}\), the projectivization formula gives a relation of degree \(m\) in \(\taut^{-1}\) with Cartan coefficients and a remainder
through two-step flags (compare Lemma~\ref{lem:expansion}).

For \(m\geq2\), the main new results are Theorems~\ref{thm:factorization},
\ref{thm:companion}, and~\ref{thm:kernel}, and the module presentations of
Corollaries~\ref{cor:structure} and~\ref{cor:presentation}.
Also new in these ranks are the isomorphism of Theorem~\ref{thm:fock}, with \(\Shift\)
acting by \(z^{-1}\), the description of the pairing by corner moments, and the
identification of the currents of \cite{N3} with the fused currents in Section~\ref{sec:negut}.

\paragraph{Organization.}
Section~\ref{sec:operator} develops the Hecke geometry and its polynomial coordinates.
Section~\ref{sec:horizontal} gives the Fock model, the pairing, and the raising and lowering series.
Section~\ref{sec:laxforms} proves the q-Miura conjugacy and then the anti-normal factorization.
Section~\ref{sec:consequences} gives the module presentations, eigenvectors, and transfer determinant,
followed by an example in rank two. Section~\ref{sec:comparison} compares the construction with
\cite{M26,GGS,CGM} and Macdonald polynomials in superspace.

\section{The Hecke correspondence}\label{sec:operator}

In rank one, \(\Gies{1}{n}\) is the Hilbert scheme of \(n\) points in \(\mathbb C^2\).
Its Hecke correspondence \(\Hecke{n}\), from Definition~\ref{def:space},
parametrizes nested ideals with a length-one quotient at the origin.
For higher rank, framed-sheaf moduli replace the Hilbert schemes. We use the flag
correspondences developed by \Negut{} \cite{N1,N2} to describe successive modifications of
these sheaves. At a fixed point, a one-step modification adds a box, and the tautological
line records its weight. Its dual therefore gives the required corner-content eigenvalue.
Two-step flags define the shift that realizes the auxiliary polynomial variable.

The parameters are generic, so distinct coloured boxes have distinct weights.

\paragraph{Smoothness and localization.}
The spaces \(\Gies{m}{n}\), \(\Hecke{n}\) and \(\HeckeTwo{n}\) are smooth, of
dimensions \(2mn\), \(m(2n+1)-1\) and \(m(2n+2)-1\), with isolated torus-fixed
points. For \(\Gies{m}{n}\) and \(\Hecke{n}\) see \cite[Sections~2
and~3.1]{N1}. We use the equivariant framed version of the constructions in
\cite{N2}, as discussed there after Proposition~4.8. The framed moduli spaces
have a universal sheaf and unobstructed deformations \cite[Section~2]{N1},
providing the representability and smoothness inputs used in those constructions.
Allowing the quotients to sit at any point of \(\mathbb C^2\) gives the framed
spaces \(\ZNtwo{x}\) and \(\ZNtwo{x,x}\). The smoothness statements of
\cite[Propositions~2.26 and~2.27]{N2} apply in this setting, and translations of \(\mathbb C^2\)
identify them with \(\Hecke{n}\times\mathbb C^2\) and
\(\HeckeTwo{n}\times\mathbb C^2\).

For such a space \(X\), write \([x]\) for the class of the skyscraper sheaf at
a fixed point \(x\), and \(T_xX\) for the tangent representation. These classes
form a basis of \(\Kloc{X}\), and \([x]|_x=\Lambda_{-1}(T_x^\vee X)\), where
\(\Lambda_{-1}(V)=\sum_i(-1)^i[\Lambda^iV]\). For a character
\(V=\sum_iw_i\) this is \(\Lambda_{-1}(V^\vee)=\prod_i(1-w_i^{-1})\). A line
bundle acts on \([x]\) by its weight at \(x\). For an equivariant map \(f\colon X\to Y\)
of such spaces, the pullback and localized equivariant pushforward have the formulas
\begin{equation}\label{eq:localisation}
 f^*[y]=\sum_{f(x)=y}\frac{\Lambda_{-1}(T_y^\vee Y)}{\Lambda_{-1}(T_x^\vee X)}\,[x],
 \qquad f_*[x]=[f(x)],
\end{equation}
the sum running over the fixed points \(x\) over \(y\).

\subsection{Fixed points}\label{sec:fixedpoints}

We compute the fixed-point coefficients of the Hecke pullback. They are the starting
point for the raising action and for the polynomial coordinates on \(\Wsp{m}\).
The universal complex \(\Wc\) of \cite[(2.4)]{N1} is the derived restriction of the
universal framed sheaf to the origin. Its projectivization parametrizes the
colength-one quotients supported there, giving the Hecke correspondence
\cite[Proposition~2.19]{N2}.

The fixed sheaf at \(\bl\) is \(\Fsh{\bl}=\bigoplus_\alpha I_{\lambda^{(\alpha)}}\),
where \(I_\lambda\subset\mathbb C[x,y]\) is the monomial ideal of \(\lambda\) and
the \(\alpha\)-th summand has framing weight \(v_\alpha\).
The minimal monomial generators of each ideal lie at its addable corners.
The relations between consecutive generators occur at their least common multiples,
which lie at the outer corners. These generators and relations give a two-term free
resolution of \(\Fsh{\bl}\). Restricting this resolution to the origin gives a complex
representing the derived fibre \(\Wc|_\bl\). The alternating sum of the generator
and relation weights is
\begin{equation}\label{eq:Wcharacter}
 [\Wc|_\bl]=\sum_{a\in\Add{\bl}}\bwt{a}-\sum_{b\in\Remp{\bl}}\bwt{b}
 =\sum_\alpha v_\alpha-(1-q_1)(1-q_2)\sum_{b\in\bl}\bwt{b},
\end{equation}
as in \cite[(2.6) and (5.4)]{N1}. A fixed subsheaf of \(\Fsh{\bl}\) of colength one at the
origin is \(\Fsh{\bl+a}\) with \(a\in\Add{\bl}\), and the line
\(\taut=\Fsh{\bl}/\Fsh{\bl+a}\) has weight \(\bwt{a}\) there.

By \eqref{eq:localisation}, the coefficients of \(\pulm\fix{\bl}\) are ratios of
cotangent Euler classes. We calculate them from the virtual relative tangent class
of this projectivization.

\begin{proposition}\label{prop:pullback}
At the fixed point \(\pt{\bl}{a}\) of \(\Hecke{n}\),
\[
 T_{\pt{\bl}{a}}\Hecke{n}-T_\bl\Gies{m}{n}
 =\taut\otimes[\Wc|_\bl]^\vee-1
 =\sum_{b\in\Add{\bl},\,b\neq a}\frac{\bwt{a}}{\bwt{b}}
  -\sum_{b\in\Remp{\bl}}\frac{\bwt{a}}{\bwt{b}},
\]
and
\begin{equation}\label{eq:pullfixed}
 \pulm\fix{\bl}
 =\sum_{a\in\Add{\bl}}
  \frac{\prod_{b\in\Remp{\bl}}(1-\bwt{b}/\bwt{a})}
       {\prod_{b\in\Add{\bl},\,b\neq a}(1-\bwt{b}/\bwt{a})}\,\pt{\bl}{a}
 =\sum_{a\in\Add{\bl}}\res{\bl}{a}\,\pt{\bl}{a}.
\end{equation}
\end{proposition}

\begin{proof}
The map \(\Psi\) in the three-term complex \cite[(2.4)--(2.5)]{N1} is surjective, so its kernel
is a \(T\)-equivariant vector bundle. Replacing the last two terms by this kernel gives the
global \(T\)-equivariant presentation \(\Wc\simeq[E^{-1}\to E^0]\) on \(\Gies{m}{n}\)
of \cite[Section~2.4]{N1}, where \(E^{-1}\) is the degree \(-1\) term and \(E^0=\ker\Psi\).
By \cite[Proposition~2.19]{N2}, the total space \(\Hecke{n}\) is cut out in
\(\mathbb P(E^0)\) by the induced section of \(\taut\otimes(E^{-1})^\vee\).
The dimension \(2mn+m-1\) of \(\Hecke{n}\), with \(\operatorname{rank}E^0-\operatorname{rank}E^{-1}=m\),
shows that this section on the total projective bundle is regular.
The relative Euler sequence and the normal sequence therefore give the virtual relative tangent class
\[
 [T\Hecke{n}]-\pulm[T\Gies{m}{n}]
 =\taut\otimes\pulm\bigl([E^0]^\vee-[E^{-1}]^\vee\bigr)-1
 =\taut\otimes\pulm[\Wc]^\vee-1.
\]
Restricting this class to \(\pt{\bl}{a}\) and using \eqref{eq:Wcharacter} gives the displayed
character. The coefficients in \eqref{eq:pullfixed} then follow from \eqref{eq:localisation}, and the last
step holds because the coefficient is the residue of \(1/(w\Yfn{\bl}(w))\) at
\(w=\bwt{a}\), which is \(\res{\bl}{a}\) in \cite{BCS}.
\end{proof}

\begin{proposition}[{following \Negut{} \cite{N1} and Bourgine--Cassia--Stoyan \cite{BCS}}]\label{prop:vertical}
For the level-\((0,m)\) vertical action \(\rhov\) of the quantum toroidal algebra
\(\mathcal E\) on \(\Fock{m}\), the modes of \(x^+(w)\) are the one-box adding
correspondences. In the normalization of \cite{BCS}, for \(k\in\mathbb Z\),
\[
 \rhov(\xplus{k})=\pushp\bigl(\taut^k\cdot\pulm(-)\bigr).
\]
In the fixed-point basis this reads
\[
 \pushp(\taut^k\pulm\fix{\bl})
 =\sum_{a\in\Add{\bl}}\bwt{a}^k\res{\bl}{a}\fix{\bl+a}.
\]
\end{proposition}

\begin{proof}
Since \(\pushp\pt{\bl}{a}=\fix{\bl+a}\)
by \eqref{eq:localisation} and \(\taut^k\) acts on \(\pt{\bl}{a}\) by
\(\bwt{a}^k\), Proposition~\ref{prop:pullback} gives
\[
 \pushp(\taut^k\pulm\fix{\bl})=\sum_a\bwt{a}^k\res{\bl}{a}\fix{\bl+a}.
\]
The vertical current \(x^+(z)=\sum_k\xplus{k}z^{-k}\) of \cite[Section~2.2.1]{BCS} acts by
\[
 x^+(z)\fix{\bl}=\sum_{a\in\Add{\bl}}\delta(\bwt{a}/z)\res{\bl}{a}\fix{\bl+a}.
\]
Taking the coefficient of \(z^{-k}\) gives the assertion.
\end{proof}

\subsection{The shift}\label{sec:shift}

The shift relates consecutive Hecke spaces. We pull a class on
\(\Hecke{n}\) back to the space of flags \(\Sh''\subset\Sh'\subset\Sh\),
multiply by the class in \eqref{eq:shiftdef}, and push forward by forgetting
\(\Sh\). On the fixed-point class \(\pt{\bmu}{r}\), the result is a linear combination
of \(\pt{\bmu+r}{b}\), where \(b\) runs through the addable corners of \(\bmu+r\).
Lemma~\ref{lem:shiftfixed} computes the coefficients of this linear combination.

\begin{definition}[The shift]\label{def:shift}
The shift associated with the two-step flag space
\(\HeckeTwo{n}\) of Section~\ref{sec:mainresults} is
\begin{equation}\label{eq:shiftdef}
 \Shift:=\pippush\bigl(-q_1q_2\,\tautone\tauttwo^{-1}\cdot\pimpull(-)\bigr)
 \colon\Wdeg{m}{n}\longrightarrow\Wdeg{m}{n+1}.
\end{equation}
\end{definition}

There are two ways to pass between these Hecke spaces: the composition
\(\pulm\pushp\) goes through the moduli space of the intermediate sheaf \(\Sh'\),
whereas \(\pippush\pimpull\) goes through the flag space.
The excess-intersection formula inserts the Euler factor
\(1-q_1q_2\tautone\tauttwo^{-1}\) into the flag push-pull.
Consequently, the difference of the two operators is precisely the shift defined above.

\begin{lemma}\label{lem:excess}
On \(\Wsp{m}\),
\[
 \pulm\pushp
 =\pippush\bigl((1-q_1q_2\,\tautone\tauttwo^{-1})\cdot\pimpull(-)\bigr)
 =\pippush\pimpull+\Shift.
\]
\end{lemma}

\begin{proof}
The square
\[
 \begin{array}{ccc}
 \HeckeTwo{n}&\xrightarrow{\;\pip\;}&\Zc_{n+1,n+2}\\[2pt]
 {\scriptstyle\pim}\big\downarrow&&\big\downarrow{\scriptstyle\pminus}\\[2pt]
 \Hecke{n}&\xrightarrow{\;\pplus\;}&\Gies{m}{n+1}
 \end{array}
\]
commutes, since both composites send \(\Sh''\subset\Sh'\subset\Sh\) to \(\Sh'\).
The displayed square is a fibre square of sets. Translations of \(\mathbb C^2\) identify the square
\cite[(2.37)]{N2} with the product of the displayed square and \(\mathbb C^2\). By
\cite[Proposition~2.27]{N2} that square is a derived fibre square with excess
bundle \(\tautone\otimes\tauttwo^{-1}\otimes p_S^*\omega_S\).
Here \(p_S\) is the common support map to the surface \(S\), and \(\omega_S\)
is the canonical bundle in the notation of \cite{N2}.
Restricting the support point to the origin gives the displayed square with the restricted
excess bundle. In \cite{N2}, the lines are numbered from the smallest sheaf, so our
\(\tautone\) and \(\tauttwo\) are the second and first lines of \cite{N2}. The canonical bundle of \(\mathbb C^2\) has weight
\(q_1q_2\) at the origin. So the excess line is \(q_1q_2\tautone\tauttwo^{-1}\),
and the excess intersection formula \cite[Proposition~3.4]{N2} is the first
equality. The second follows from \eqref{eq:shiftdef}.
\end{proof}

\begin{lemma}\label{lem:shiftfixed}
The shift satisfies
\begin{equation}\label{eq:displacement}
 \Shift\,\taut-q_3\,\taut\,\Shift=\pulm\pushp\,\taut.
\end{equation}
For an \(m\)-tuple \(\bmu\) and \(r\in\Add{\bmu}\), the fixed-point action of \(\Shift\) is
\[
 \Shift\,\pt{\bmu}{r}
 =\sum_{b\in\Add{\bmu+r}}
   \res{\bmu+r}{b}\,\frac{\ceig{b}}{\ceig{b}-\ceig{\bplus{r}}}\,
   \pt{\bmu+r}{b}.
\]
\end{lemma}

\begin{proof}
Since \(\pimpull\taut=\tautone\) and \(\pippull\taut=\tauttwo\), the projection formula gives
\[
 \taut\,\Shift=\pippush\bigl(-q_1q_2\,\tautone\cdot\pimpull(-)\bigr)
 =-q_3^{-1}\,\pippush\pimpull\,\taut .
\]
So \(\pippush\pimpull\taut=-q_3\taut\Shift\), and Lemma~\ref{lem:excess}
composed with \(\taut\) gives \eqref{eq:displacement}. The fixed points of \(\HeckeTwo{n}\) over
\(\pt{\bmu}{r}\) are the flags \(\Fsh{\bmu+r+b}\subset\Fsh{\bmu+r}\subset\Fsh{\bmu}\),
\(b\in\Add{\bmu+r}\), and \(\pip\) sends them to \(\pt{\bmu+r}{b}\). With \eqref{eq:localisation},
\(\Shift\) sends \(\pt{\bmu}{r}\) into the span of the \(\pt{\bmu+r}{b}\).
Write \(\Shift\pt{\bmu}{r}=\sum_bs_b\pt{\bmu+r}{b}\). By \eqref{eq:pullfixed},
\(\pulm\pushp\pt{\bmu}{r}=\pulm\fix{\bmu+r}=\sum_b\res{\bmu+r}{b}\pt{\bmu+r}{b}\). So
\eqref{eq:displacement}, applied to \(\pt{\bmu}{r}\), says that
\((\bwt{r}-q_3\bwt{b})s_b=\bwt{r}\res{\bmu+r}{b}\) for every \(b\). Here \(q_3\bwt{b}=\bwt{r}\) would
mean \(\ceig{b}=\ceig{\bplus{r}}\), hence \(b=\bplus{r}\), which is not an addable corner of
\(\bmu+r\). So \(s_b=\res{\bmu+r}{b}\,\ceig{b}/(\ceig{b}-\ceig{\bplus{r}})\).
\end{proof}

\begin{lemma}\label{lem:shiftip}
For \(A,B\in\Wsp{m}\),
\[
 \ip{\Shift A}{\Shift B}=\Mq\,\ip{A}{B}.
\]
The classes \(\taut^{-j}\pulm\fix{\bl}\), \(0\leq j<m\), are orthogonal
to the image of \(\Shift\).
\end{lemma}

Thus \(\Mq^{-1/2}\Shift\) is an isometric embedding over \(\K(\Mq^{1/2})\).
The proof is given in Section~\ref{sec:form}, after the fixed-point norms and
corner moments have been computed.

\subsection{The fibre over a fixed point}\label{sec:fibre}

We prove the \(\K[\Shift]\)-module decomposition of Lemma~\ref{lem:generation}(i)
one source fixed point at a time. For the point of \(\Gies{m}{n}\) indexed by \(\bl\), put
\[
 \Wfib{\bl}:=\operatorname{span}_{\K}\{\pt{\bl}{a}:a\in\Add{\bl}\},
 \qquad \Wsp{m}=\bigoplus_\bl\Wfib{\bl}.
\]
Functions on the addable-corner contents give coordinates on \(\Wfib{\bl}\)
for the pullback classes and the shift images from the preceding degree.

By Proposition~\ref{prop:fixed}, the spectrum of \(\GLax{m}\) on \(\Wfib{\bl}\)
consists of the contents \(\ceig{a}\). We use \(x\) for the spectral variable.
The Nekrasov \(Y\)-function records these spectral points as the zeros of the
rational function
\[
 \Yfn{\bl}(1/x)=\frac{\prod_{a\in\Add{\bl}}(1-x/\ceig{a})}
                   {\prod_{b\in\Remp{\bl}}(1-x/\ceig{b})}.
\]
The poles of \(\Yfn{\bl}(1/x)\) are the outer-corner contents \(\ceig{b}\), and
\(\Yfn{\bl}(1/x)\) has value \(1\) at \(x=0\).
The residue formulas below give the coefficients in the function model.
Each colour has one more addable than removable corner, and the contents
telescope in each colour \cite{BCS}:
\begin{equation}\label{eq:cornerproducts}
 |\Add{\bl}|-|\Rem{\bl}|=m,\qquad
 \frac{\prod_{a\in\Add{\bl}}\ceig{a}}{\prod_{b\in\Remp{\bl}}\ceig{b}}
 =\prod_\alpha v_\alpha^{-1}.
\end{equation}
Hence \(\Yfn{\bl}(1/x)=-\betaz^{-1}x^m+O(x^{m-1})\) at infinity.

\begin{lemma}\label{lem:residues}
\[
 \Yfn{\bl}(1/x)^{-1}=\sum_{a\in\Add{\bl}}\frac{\res{\bl}{a}\,\ceig{a}}{\ceig{a}-x},
 \qquad
 \cores{\bl}{r}=-\Res_{x=\ceig{\bplus{r}}}\frac{\Yfn{\bl}(1/x)}{x}.
\]
\end{lemma}

\begin{proof}
In the variable \(x=1/w\), the conventions of \cite{BCS} give
\(\res{\bl}{a}=-\Res_{x=\ceig{a}}\Yfn{\bl}(1/x)^{-1}/x\) and the stated
residue formula for \(\cores{\bl}{r}\).
The shift by \(q_3\) in the definition of \(\cores{\bl}{r}\) becomes the outer
corner \(\bplus{r}\). The first formula is the partial fraction expansion of
\(\Yfn{\bl}(1/x)^{-1}\), which vanishes at infinity and has simple poles at the
\(\ceig{a}\).
\end{proof}

\begin{lemma}\label{lem:fibre}
\begin{enumerate}
\item[(i)] For every function \(f\) on \(\{\ceig{a}:a\in\Add{\bl}\}\),
\[
 f(\taut^{-1})\,\pulm\fix{\bl}
 =\sum_{a\in\Add{\bl}}\res{\bl}{a}\,f(\ceig{a})\,\pt{\bl}{a},
\]
and \(f\mapsto f(\taut^{-1})\pulm\fix{\bl}\) is an isomorphism from these
functions onto \(\Wfib{\bl}\).
\item[(ii)] \(\taut^{-j}\pulm\fix{\bl}\) corresponds to the function \(x^j\),
and \(\Shift\,\pt{\bl-r}{r}\), for \(r\in\Rem{\bl}\), to the function
\(x/(x-\ceig{\bplus{r}})\).
\end{enumerate}
\end{lemma}

\begin{proof}
(i) This is \eqref{eq:pullfixed}, since the line \(\taut^{-1}\) has weight \(\ceig{a}\) at
\(\pt{\bl}{a}\). The residues are nonzero and the \(\ceig{a}\) are distinct.
(ii) The first statement is (i) with \(f=x^j\). The second is
Lemma~\ref{lem:shiftfixed} with \(\bmu=\bl-r\).
\end{proof}

\subsection{The polynomial decomposition}\label{sec:freeness}

We prove Lemma~\ref{lem:generation}(i) using an interpolation basis of \(\Wfib{\bl}\).
The powers \(1,x,\dots,x^{m-1}\) give the pullback
generators, and the rational terms with poles at outer corners give the images of
the shift.

\begin{proposition}\label{prop:free}
The classes \(\taut^{-j}\pulm\fix{\bl}\), \(0\leq j<m\), and
\(\Shift\,\pt{\bl-r}{r}\), \(r\in\Rem{\bl}\), form a basis of \(\Wfib{\bl}\).
\end{proposition}

\begin{proof}
By Lemma~\ref{lem:fibre}, a linear relation among these classes is a
rational function
\[
 g(x)=P(x)+\sum_{r\in\Rem{\bl}}t_r\,\frac{x}{x-\ceig{\bplus{r}}},\qquad\deg P<m,
\]
that vanishes at every \(\ceig{a}\). The polynomial
\(g(x)\prod_{b\in\Remp{\bl}}(x-\ceig{b})\) has degree at most
\(m-1+|\Rem{\bl}|=|\Add{\bl}|-1\) and \(|\Add{\bl}|\) distinct zeros, so \(g=0\).
The residues of \(g\) at the distinct nonzero \(\ceig{\bplus{r}}\) give \(t_r=0\),
and then \(P=0\). There are \(m+|\Rem{\bl}|=|\Add{\bl}|=\dim\Wfib{\bl}\) classes.
\end{proof}

\begin{proof}[Proof of Lemma~\ref{lem:generation}]
(i) Every fixed point of \(\Zc_{N-1,N}\) is \(\pt{\bl-r}{r}\) for exactly one pair
\((\bl,r)\) with \(|\bl|=N\) and \(r\in\Rem{\bl}\). So \(\Shift\) maps the
fixed-point basis of \(\Wdeg{m}{N-1}\) to the classes \(\Shift\,\pt{\bl-r}{r}\),
and Proposition~\ref{prop:free} gives
\[
 \Wdeg{m}{N}=\bigoplus_{j<m}\taut^{-j}\pulm\Fdeg{m}{N}\;\oplus\;\Shift\,\Wdeg{m}{N-1},
\]
with \(\Shift\) injective. By induction
\(\Wdeg{m}{N}=\bigoplus_{k=0}^N\Shift^k\bigoplus_{j<m}\taut^{-j}\pulm\Fdeg{m}{N-k}\).
This proves freeness over \(\K[\Shift]\) with the stated generator space.

(ii) Multiply \eqref{eq:displacement} by \(\taut^{-1}\) on the left and on the right.
\end{proof}

This decomposition identifies \(\Wsp{m}\) with polynomials in \(z^{-1}\) whose
coefficients lie in the generator space
\(\bigoplus_{j=0}^{m-1}\taut^{-j}\pulm\Fock{m}\): the monomial \(z^{-k}u\)
corresponds to \(\Shift^k u\). Under this identification, \(\Shift\) acts by
multiplication by \(z^{-1}\).

Lemma~\ref{lem:generation}(ii) gives, by induction on \(k\),
\begin{equation}\label{eq:iterated}
 \dtaut\,\Shift^k=q_3^k\,\Shift^k\,\dtaut
 +\sum_{i=0}^{k-1}q_3^{\,i}\,\Shift^i\,\dtaut\,\pulm\pushp\,\Shift^{k-1-i}.
\end{equation}

\subsection{Newton polynomials}\label{sec:newton}

The elementary Newton polynomials \(\Newton{r}\) of \eqref{eq:newtoncorners}, with
\(\Newton{r}=0\) for \(r<0\), are the elementary symmetric functions of the difference of
two alphabets: the contents of the addable corners and the contents of the outer corners.
The corresponding diagonal operators act on \(\fix{\bl}\), \(\kket{\bl}\) and
\(\pt{\bl}{a}\) by \(\Newton{r}(\bl)\). Factoring
\(-x/c\) out of each factor \(1-x/c\) of \(\Yfn{\bl}(1/x)\) and using \eqref{eq:cornerproducts}
gives, at \(x=\infty\),
\begin{equation}\label{eq:Yexpansion}
 -\betaz\,\Yfn{\bl}(1/x)=\sum_{r\geq0}(-1)^r\,\Newton{r}(\bl)\,x^{m-r}.
\end{equation}
Lemma~\ref{lem:residues} gives the corresponding partial fraction decomposition
\[
 \Yfn{\bl}(1/x)=-\betaz^{-1}p_\bl(x)
 -\sum_{r\in\Rem{\bl}}\cores{\bl}{r}\,
 \frac{x}{x-\ceig{\bplus{r}}},
\]
with \(p_\bl\) as in \eqref{eq:companionpolynomial}.
The rational terms represent images of the shift by Lemma~\ref{lem:fibre}(ii).
The rational terms vanish at \(x=0\), so \(\Yfn{\bl}(1/x)|_{x=0}=1\) forces
\(p_\bl(0)=-\betaz\). This explains the constant term in the companion polynomial.

The Cartan operators act diagonally on the same bases, with eigenvalues
given by the power sums of this difference of alphabets,
\[
 \Cartan{k}(\bl):=\sum_{a\in\Add{\bl}}\ceig{a}^k-\sum_{b\in\Remp{\bl}}\ceig{b}^k.
\]
These definitions give
\(\sum_r\Newton{r}(\bl)u^r=\exp\bigl(\sum_k(-1)^{k-1}\Cartan{k}(\bl)u^k/k\bigr)\) and
\(\Newton{1}=\Cartan{1}\). Put \(\epsilon:=(1-q_1^{-1})(1-q_2^{-1})\). By \eqref{eq:Wcharacter},
\(\Cartan{1}(\bl)\) is the value at \(\bl\) of the class \([\Wc^\vee]\), and
\begin{equation}\label{eq:Hone}
 \Cartan{1}(\bl)=\sum_\alpha v_\alpha^{-1}-\epsilon\sum_{b\in\bl}\ceig{b}.
\end{equation}

\section{The horizontal realization}\label{sec:horizontal}

The Fock isomorphism \(\Hiso\) identifies fixed-point classes with generalized
Macdonald functions. The normalization of \(\Hiso\) determines the constants in the
horizontal currents and in the transported Euler pairing.

\subsection{The Fock module}\label{sec:bcs}

\paragraph{Generalized Macdonald functions.}
The horizontal module \(\La^{\otimes m}\) of level \((m,0)\) is the ring of
symmetric functions in \(m\) alphabets \(x^{(1)},\dots,x^{(m)}\), with power sums
\(\pk{k}{\alpha}\) normalized as in \cite[Section~2.2.2]{BCS} and weights
\(u_\alpha\). The colour \(\alpha\) now labels the alphabet \(x^{(\alpha)}\).
Let \(\GMP{\lambda}\) be the Macdonald function
with parameters \((\Mq,\Mt)\). The generalized Macdonald function
\(\GMP{\bl}(x^\bullet|\bu)\) is the eigenvector of the zero mode
\(\rho^{(m,0)}_{\bu}(\xplus{0})\) of \cite{BCS} of the form
\[
 \GMP{\bl}=\GMP{\lambda^{(1)}}(x^{(1)})\cdots\GMP{\lambda^{(m)}}(x^{(m)})
 +\text{lower terms}.
\]
The lower terms are products \(\prod_\alpha\GMP{\mu^{(\alpha)}}(x^{(\alpha)})\)
in which boxes have moved toward alphabets with smaller indices,
\(|\mu^{(1)}|+\dots+|\mu^{(k)}|\geq|\lambda^{(1)}|+\dots+|\lambda^{(k)}|\) for every
\(k\), with \(\bmu\neq\bl\). For a partition \(\lambda\), write
\(n(\lambda'):=\sum_i\binom{\lambda_i}{2}\), where \(\lambda'\) is the transposed
partition. With \(\pvac\) the vacuum specialization of \cite{BCS}, the
normalization is\footnote{We correct two misprints in
arXiv:2508.19704v2. In Section~2.2.2, the factor indexed by \(i'\) in the
product for \(\GMP{\lambda}(\pvac)\) must be inverted. In Section~3.1, the
ordering must transfer boxes from the second partition to the first, in agreement
with Appendix~A.2 and the example for \(\GMP{(\varnothing,(1))}\).}
\[
 \GMPn{\bl}:=\frac{\GMP{\bl}}{\prod_\alpha\GMP{\lambda^{(\alpha)}}(\pvac)},
 \qquad
 \GMP{\lambda}(\pvac)=\frac{(-1)^{|\lambda|}\Mq^{n(\lambda')}}
   {\prod_{s\in\lambda}(1-\Mq^{a(s)}\Mt^{l(s)+1})},
\]
where \(a(s)\) and \(l(s)\) are the arm and leg of the box \(s\). Thus
\(\GMPn{\lambda}=(-1)^{|\lambda|}\Mq^{-n(\lambda')}J_\lambda\) in rank one.

The Nekrasov factors of \cite{BCS}, in the zero-based box indexing of
Section~\ref{sec:mainresults}, are
\begin{align}
 \Nekbare{\lambda}{\mu}(s)
 &:=\prod_{(i,j)\in\lambda}\bigl(1-s\,q_1^{\lambda'_{j+1}-i}q_2^{j+1-\mu_{i+1}}\bigr)
   \prod_{(i,j)\in\mu}\bigl(1-s\,q_1^{i+1-\mu'_{j+1}}q_2^{\lambda_{i+1}-j}\bigr),\notag\\
 \Nek{\lambda}{\mu}(s)
 &:=\Nekbare{\lambda}{\mu}(s)\prod_{(i,j)\in\lambda}
       \bigl(-q_3^{-1}s\,q_1^{i}q_2^{j}\bigr)^{-1}.\label{eq:nekrasovfactors}
\end{align}
Here \(\lambda_1\geq\lambda_2\geq\dots\) are the parts, with missing parts
equal to zero. The modified factor \(\Nek{\lambda}{\mu}\) specifies the
normalization of the horizontal fixed-point basis.

\begin{proposition}[{Bourgine--Cassia--Stoyan \cite{BCS}}]\label{prop:horizontal}
The map \(\Hiso\) is the graded linear isomorphism of \cite[Section~3.3]{BCS} from the
vertical to the horizontal Fock space, rescaled by \(\kappam{m}^{-N}\) in
degree \(N\), where \(\kappam{m}:=q_1\gamma^{1-m/2}\).
\begin{enumerate}
\item[(i)] The horizontal vector of the fixed point \(\fix{\bl}\) is
\begin{equation}\label{eq:kket}
 \kket{\bl}=(-1)^{|\bl|}
   \prod_{\alpha<\beta}\Nek{\lambda^{(\alpha)}}{\lambda^{(\beta)}}(v_\alpha/v_\beta)\,
   \GMPn{\bl}(x^\bullet|\bu).
\end{equation}
\item[(ii)] As maps from \(\Fock{m}\) to \(\La^{\otimes m}\),
\[
 \Hiso\,\pushp\pulm
 =(1-\Mt)\Bigl(\sum_{\alpha=1}^m\pk{1}{\alpha}\Bigr)\Hiso.
\]
\end{enumerate}
\end{proposition}

\begin{proof}
The isomorphism of Bourgine--Cassia--Stoyan \cite[Section~3.3]{BCS} sends the
vertical basis vector \(\fix{\bl}\), identified in Proposition~\ref{prop:vertical}, to
\begin{gather*}
 (-q_1\gamma^{1/2})^{|\bl|}\,
 \gamma^{-(m-1)|\bl|/2+\sum_\alpha(m-\alpha)|\lambda^{(\alpha)}|}\,
 G^*_{\bl}(\bu)^{-1}\,\GMPn{\bl},\\
 G^*_{\bl}(\bu)=\gamma^{\sum_\alpha(m-\alpha)|\lambda^{(\alpha)}|}
 \prod_{\alpha<\beta}\Nek{\lambda^{(\alpha)}}{\lambda^{(\beta)}}(v_\alpha/v_\beta)^{-1},
\end{gather*}
which is \(\kappam{m}^{|\bl|}(-1)^{|\bl|}\prod_{\alpha<\beta}
\Nek{\lambda^{(\alpha)}}{\lambda^{(\beta)}}(v_\alpha/v_\beta)\,\GMPn{\bl}\).
The BCS isomorphism intertwines the vertical zero mode \(\xplus{0}\), which is \(\pushp\pulm\) by
Proposition~\ref{prop:vertical}, with the image of \(\xplus{0}\) under the Miki
automorphism. That image is \(-(\gamma-\gamma^{-1})^{-1}a_{-1}\), which acts by
\(-\gamma^{1-m/2}(1-q_1)(\pk{1}{1}+\dots+\pk{1}{m})
=\kappam{m}(1-\Mt)(\pk{1}{1}+\dots+\pk{1}{m})\) \cite[Section~2]{BCS}.
Multiplying by \(\kappam{m}^{-N}\) in degree \(N\) removes \(\kappam{m}\) from both
statements.
\end{proof}

In rank one, the horizontal vector is
\(\Hiso\fix{\lambda}=\Mq^{-n(\lambda')}J_\lambda\), where \(J_\lambda\) is the
integral Macdonald function \cite{Macdonald}.

\begin{proof}[Proof of Theorem~\ref{thm:fock}(i)]
By Lemma~\ref{lem:generation}(i), the classes \(\Shift^k\taut^{-j}\pulm\fix{\bl}\) form a basis of
\(\Wsp{m}\), and by Proposition~\ref{prop:horizontal} the \(\kket{\bl}\) form a basis of
\(\La^{\otimes m}\). Thus \(\Omap\) is an isomorphism. By construction
\(\Omap(z^{-1}v)=\Shift\,\Omap(v)\), and \(\Omap\) maps the constant sector
onto \(\bigoplus_{j<m}\taut^{-j}\pulm\Fock{m}\).
\end{proof}

\paragraph{The horizontal current.}
The vertical and horizontal representations of the quantum toroidal algebra \(\mathcal E\) are
\(\rhov\) on \(\Fock{m}\) and \(\rhou\) on \(\La^{\otimes m}\)
\cite[Section~2.2]{BCS}. They are related by the Miki automorphism \(\Miki\)
\cite[(2.9)]{BCS}. By Proposition~\ref{prop:vertical},
\(\rhov(\xplus{d})=\pushp(\taut^d\pulm(-))\), and by \cite[(2.22) and Section~3.3]{BCS}, in the
normalization of Proposition~\ref{prop:horizontal},
\(\Hiso\rhov(e)\Hiso^{-1}=\kappam{m}^{-d}\rhou(\Miki(e))\) for \(e\) of degree \(d\). The horizontal
current \(\Xcur(w)=\rhou(x^-(w))=\sum_{n\in\mathbb Z}\Xcur_n\,w^{-n}\) is
\begin{equation}\label{eq:Xcur}
 \Xcur(w)=\sum_{\alpha=1}^mu_\alpha^{-1}
 \exp\Bigl(-\sum_{k\geq1}\frac{w^k}{k}\,\gamma^{mk}(1-q_1^k)q_3^{-(\alpha-1)k}\pk{k}{\alpha}\Bigr)
 \exp\Bigl(\sum_{k\geq1}w^{-k}q_3^k\gamma^{-mk}(1-q_2^k)\,\partial^{(\alpha)}_k\Bigr)
\end{equation}
\cite[(2.52)]{BCS}, where \(\partial_{\alpha,k}=\partial/\partial\pk{k}{\alpha}\) and
\(\partial^{(\alpha)}_k=q_3^{(\alpha-1)k}\partial_{\alpha,k}
+(1-q_3^{-k})\sum_{\beta>\alpha}q_3^{(\beta-1)k}\partial_{\beta,k}\).
The first Cartan operator is represented by the rescaled zero mode
\(\Hiso\Cartan{1}\Hiso^{-1}=-\gamma\Xcur_0\), since
\(\Cartan{1}=-\rhov(a_{-1})/(1-q_3^{-1})\) by \cite[(2.24)]{BCS} and
\(\Miki(a_{-1})=(\gamma-\gamma^{-1})x^-_0\).

For a series \(C(w)=\sum_nC_nw^{-n}\) we
write \(C(w)_{>0}=\sum_{n\geq1}C_{-n}w^n\), \(C(w)_{<0}=\sum_{n\geq1}C_nw^{-n}\) and
\(C(w)_{\leq0}=\sum_{n\geq0}C_nw^{-n}\). For a series in \(z\), \(\Dq C(z)=C(z)-C(q_3z)\).

\subsection{The pairing}\label{sec:form}

The rescaled vertical norms \(\srcnorm{\bl}\) of Section~\ref{sec:mainresults}
have the product formula
\begin{equation}\label{eq:sourcenormproduct}
 \srcnorm{\bl}=\Mq^{|\bl|}\prod_{\alpha,\beta}
 \Nek{\lambda^{(\alpha)}}{\lambda^{(\beta)}}(v_\alpha/v_\beta)
\end{equation}
by \cite[(2.25)]{BCS}. In particular, \(\srcnorm{\varnothing}=1\).

\begin{lemma}\label{lem:normratio}
\(\srcnorm{\bmu+a}/\srcnorm{\bmu}=-\Mq\,\cores{\bmu+a}{a}/\res{\bmu}{a}\) for
\(a\in\Add{\bmu}\).
\end{lemma}

\begin{proof}
The pairing of \cite[(2.25)]{BCS} is contravariant for \((\xplus{k})^\dagger=-x^-_k\).
The vertical current acts by
\[
 \rhov(x^-(z))\fix{\bmu+a}
 =\gamma^{-m}\sum_{r\in\Rem{\bmu+a}}
   \delta(\bwt{r}/z)\cores{\bmu+a}{r}\fix{\bmu+a-r}
\]
\cite[Section~2.2.1]{BCS}, and \(\xplus{k}\) is given by
Proposition~\ref{prop:vertical}. Pairing \(\xplus{k}\fix{\bmu}\) with
\(\fix{\bmu+a}\) gives the ratio of the unrescaled norms as
\(-\gamma^{-m}\cores{\bmu+a}{a}/\res{\bmu}{a}\).
The degree rescaling multiplies this ratio by \(\Mq\gamma^m\), giving
the stated formula.
\end{proof}

Write \(\Tc\) for the tautological bundle, with character
\(\sum_{b\in\bl}\bwt{b}\) at \(\bl\). The tangent space at \(\bl\) has the
character
\[
 T_\bl\Gies{m}{n}=\bigl(\chi(\framing,\framing)-\chi(\Wc|_\bl,\Wc|_\bl)\bigr)^\vee,
 \qquad
 \chi(A,B)=\frac{A^\vee B}{(1-q_1^{-1})(1-q_2^{-1})},
\]
as in \cite[Section~2]{N1}.

\begin{lemma}\label{lem:pplus}
At the fixed point \(\pt{\bl}{a}\),
\(T_{\pt{\bl}{a}}\Hecke{n}-T_{\bl+a}\Gies{m}{n+1}
=-q_3\,\taut^{-1}\otimes[\Wc|_{\bl+a}]-1\), and
\[
 \pplus{}^*\fix{\bmu}
 =\sum_{r\in\Rem{\bmu}}(-1)^{m-1}\prod_\beta\frac{\bwt{r}}{q_3v_\beta}\;
   \cores{\bmu}{r}\,\pt{\bmu-r}{r}.
\]
\end{lemma}

\begin{proof}
Since \([\Wc|_{\bl+a}]=[\Wc|_\bl]-(1-q_1)(1-q_2)\bwt{a}\), the character formula
gives
\(T_{\bl+a}\Gies{m}{n+1}-T_\bl\Gies{m}{n}
=\taut\otimes[\Wc|_\bl]^\vee+q_3\,\taut^{-1}\otimes[\Wc|_{\bl+a}]\), and the
first statement follows from Proposition~\ref{prop:pullback}. In
\(-q_3\taut^{-1}[\Wc|_{\bl+a}]\) the outer corner \(\bplus{a}\) of \(\bl+a\)
contributes the weight \(1\), which cancels the \(-1\). The remaining weights,
and \eqref{eq:localisation}, give the coefficient
\(\prod_{b\in\Add{\bl+a}}(1-\bwt{a}/(q_3\bwt{b}))
/\prod_{b\in\Remp{\bl+a},\,b\neq\bplus{a}}(1-\bwt{a}/(q_3\bwt{b}))\). By the
formula for \(\cores{\bl+a}{a}\) as a residue and \eqref{eq:cornerproducts},
this is the displayed multiple of \(\cores{\bl+a}{a}\).
\end{proof}

\paragraph{Normalization of the pairing.}
For this comparison with \Negut{} \cite[(3.35)]{N3}, identify the parameters
\((q_1,q_2,\mathbf u)\) of that paper with our
\((\Mt^{-1},\Mq,\bv)\), so that the tautological bundle of \cite{N3} is \(\Tc\).
On \(\Kloc{\Gies{m}{n}}\), our pairing is the resulting pairing multiplied
in degree \(n\) by \((\Mq q_3^m)^n\),
\[
 \ip{\alpha}{\beta}:=(-1)^{mn}\Mq^{n}\,\chi_{\Gies{m}{n}}\bigl(\alpha\otimes\beta\otimes
 \det(q_3\,\framing\otimes\Tc^\vee)\bigr).
\]
The pairing \eqref{eq:eulerpairing} on \(\Wdeg{m}{n}\) uses the same twist,
pulled back along \(\pminus\).

\begin{lemma}[{\Negut{} \cite{N3}}]\label{lem:negut}
Let \(x^-_k\) be the box-removing modes of the vertical current in \cite{BCS}.
In the parameter convention above, the one-box adjunction
\cite[(3.39)]{N3} gives, for \(k\in\mathbb Z\),
\[
 \bigl(\rhov(\xplus{k})\bigr)^\dagger=
 (-1)^m\Mq\,q_3^m\prod_\alpha v_\alpha\cdot\pminus_*\bigl(\taut^{k-m}\cdot\pplus{}^*(-)\bigr)
 =-\Mq\gamma^m\rhov(x^-_k).
\]
\end{lemma}

\begin{proof}
On \(\Hecke{n}\), \(\pplus{}^*\Tc=\pulm\Tc+\taut\) in \(K\)-theory, so
\[
 \pplus{}^*\det(q_3\framing\otimes\Tc^\vee)
 =q_3^m\prod_\alpha v_\alpha\cdot\taut^{-m}\,\pulm\det(q_3\framing\otimes\Tc^\vee).
\]
The projection formulas for \(\pplus\) and \(\pminus\), with the ratio \((-1)^m\Mq\) of the constants
\((-1)^{mn}\Mq^n\) in degrees \(n+1\) and \(n\), give the first expression for the adjoint. By
Lemma~\ref{lem:pplus} and \eqref{eq:localisation},
\[
 \pminus_*\bigl(\taut^{k-m}\cdot\pplus{}^*\fix{\bmu}\bigr)
 =\frac{(-1)^{m-1}}{q_3^m\prod_\alpha v_\alpha}\sum_{r\in\Rem{\bmu}}\bwt{r}^k\,\cores{\bmu}{r}\,\fix{\bmu-r},
\]
The vertical action of \cite[Section~2.2.1]{BCS} is
\[
 \rhov(x^-_k)\fix{\bmu}
 =\gamma^{-m}\sum_{r\in\Rem{\bmu}}\bwt{r}^k\cores{\bmu}{r}\fix{\bmu-r},
\]
which gives the second equality.
\end{proof}

\begin{proposition}\label{prop:euler}
The fixed-point classes are orthogonal, with
\begin{equation}\label{eq:fixednorm}
 \ip{\pt{\bl}{a}}{\pt{\bl}{a}}
 =\srcnorm{\bl}\,
 \frac{\prod_{b\in\Add{\bl},\,b\neq a}(1-\bwt{b}/\bwt{a})}
      {\prod_{b\in\Remp{\bl}}(1-\bwt{b}/\bwt{a})}
 =\frac{\srcnorm{\bl}}{\res{\bl}{a}}.
\end{equation}
The rescaled vertical norm \(\srcnorm{\bl}\) is the norm of the source
fixed-point class for the twisted Euler pairing.
\end{proposition}

\begin{proof}
For \(|\bl|=n\), the geometric source norm is
\begin{equation}\label{eq:sourcenormEuler}
 N_\bl:=\Mq^{n}\prod_{b\in\bl}\prod_\beta\Bigl(-\frac{q_3v_\beta}{\bwt{b}}\Bigr)\,
 \Lambda_{-1}\bigl(T^\vee_\bl\Gies{m}{n}\bigr).
\end{equation}
Both \(\srcnorm{\bl}\) and \(N_\bl\) equal \(1\) at \(\bl=\varnothing\).
Proposition~\ref{prop:pullback}, Lemma~\ref{lem:pplus} and \eqref{eq:localisation} give
\[
 \frac{\Lambda_{-1}(T^\vee_{\bl}\Gies{m}{n})}
      {\Lambda_{-1}(T^\vee_{\pt{\bl}{a}}\Hecke{n})}=\res{\bl}{a},
 \qquad
 \frac{\Lambda_{-1}(T^\vee_{\bl+a}\Gies{m}{n+1})}
      {\Lambda_{-1}(T^\vee_{\pt{\bl}{a}}\Hecke{n})}
 =(-1)^{m-1}\prod_\beta\frac{\bwt{a}}{q_3v_\beta}\,\cores{\bl+a}{a}.
\]
Consequently,
\[
 \frac{N_{\bl+a}}{N_\bl}
 =\Mq\prod_\beta\Bigl(-\frac{q_3v_\beta}{\bwt{a}}\Bigr)\cdot
  (-1)^{m-1}\prod_\beta\frac{\bwt{a}}{q_3v_\beta}\cdot
  \frac{\cores{\bl+a}{a}}{\res{\bl}{a}}
 =-\Mq\,\frac{\cores{\bl+a}{a}}{\res{\bl}{a}},
\]
which is the ratio of Lemma~\ref{lem:normratio}. Induction gives \(N_\bl=\srcnorm{\bl}\).

The classes of distinct fixed points have disjoint supports, so their product vanishes and they are
orthogonal. For a fixed point, \(\srcnorm{\bl}/\res{\bl}{a}\) equals
\((-1)^{mn}\Mq^n\det(q_3\framing\otimes\Tc^\vee)|_\bl\,
\Lambda_{-1}(T^\vee_{\pt{\bl}{a}}\Hecke{n})\), by Proposition~\ref{prop:pullback}, and for a line
bundle \(\mathcal N\) the Euler characteristic \(\chi(\pt{\bl}{a}\otimes\pt{\bl}{a}\otimes\mathcal N)\) is the
weight of \(\mathcal N\) at \(\pt{\bl}{a}\) times \(\Lambda_{-1}(T^\vee_{\pt{\bl}{a}}\Hecke{n})\).
The product in \eqref{eq:fixednorm} follows by evaluating the residue of
\(1/(w\Yfn{\bl}(w))\) at \(w=\bwt{a}\).
\end{proof}

The corner measure of \(\bl\) is
\(\cornermeas{\bl}:=\sum_{a\in\Add{\bl}}\res{\bl}{a}\,\delta_{\ceig{a}}\).
The moments of \(\cornermeas{\bl}\) are
\begin{equation}\label{eq:cornermoments}
 \mom{s}{\bl}:=\sum_{a\in\Add{\bl}}\res{\bl}{a}\ceig{a}^{s}
 =\int x^s\,d\cornermeas{\bl},\qquad s\geq0.
\end{equation}
Let \(f\) and \(g\) be functions on the finite set
\(\{\ceig{a}:a\in\Add{\bl}\}\). By Lemma~\ref{lem:fibre}(i), they represent the
classes \(f(\taut^{-1})\pulm\fix{\bl}\) and \(g(\taut^{-1})\pulm\fix{\bl}\)
in \(\Wfib{\bl}\), with fixed-point coefficients
\(\res{\bl}{a}f(\ceig{a})\) and \(\res{\bl}{a}g(\ceig{a})\).
Proposition~\ref{prop:euler} therefore gives
\begin{equation}\label{eq:fibrepairing}
 \ip{f(\taut^{-1})\pulm\fix{\bl}}{g(\taut^{-1})\pulm\fix{\bl}}
 =\srcnorm{\bl}\sum_{a\in\Add{\bl}}\res{\bl}{a}f(\ceig{a})g(\ceig{a})
 =\srcnorm{\bl}\int fg\,d\cornermeas{\bl}.
\end{equation}
The integral denotes this finite weighted sum, and different summands
\(\Wfib{\bl}\) are orthogonal.

\begin{lemma}\label{lem:moments}
For a formal variable \(u\), the moments satisfy
\begin{equation}\label{eq:momentseries}
 \sum_{s\geq0}\mom{s}{\bl}\,u^s=1-\Yfn{\bl}(u)^{-1},
\end{equation}
with the right side expanded at \(u=0\).
The initial moments are \(\mom{0}{\bl}=1\), \(\mom{s}{\bl}=0\) for \(1\leq s<m\),
and \(\mom{m}{\bl}=\betaz\), where \(\betaz\) is the framing scalar of
Section~\ref{sec:mainresults}.
For \(r,s\in\Rem{\bl}\) and \(0\leq j<m\),
\[
 \int x^j\,\frac{x}{x-\ceig{\bplus{r}}}\,d\cornermeas{\bl}=0,\qquad
 \int\frac{x}{x-\ceig{\bplus{r}}}\,\frac{x}{x-\ceig{\bplus{s}}}\,d\cornermeas{\bl}
 =-\frac{\delta_{rs}}{\cores{\bl}{r}}.
\]
\end{lemma}

\begin{proof}
The partial fraction formula of Lemma~\ref{lem:residues}, evaluated at \(x=0\), gives
\(\sum_a\res{\bl}{a}=1\). Substituting \(x=u^{-1}\) gives
\[
 1-\Yfn{\bl}(u)^{-1}
 =\sum_{a\in\Add{\bl}}\frac{\res{\bl}{a}}{1-u\,\ceig{a}}.
\]
Expanding at \(u=0\) proves \eqref{eq:momentseries}, and
\eqref{eq:Yexpansion} gives the initial moments.

For the remaining integrals, the same partial fraction formula gives
\[
 \int f\,d\cornermeas{\bl}
 =-\sum_a\Res_{x=\ceig{a}}f(x)\,\Yfn{\bl}(1/x)^{-1}dx/x.
\]
The function
\(\Yfn{\bl}(1/x)^{-1}\) has value \(1\) at \(0\), simple zeros at the \(\ceig{\bplus{r}}\), and the
leading term \(-\betaz x^{-m}\) at infinity, by \eqref{eq:Yexpansion}. The
residue theorem on \(\mathbb P^1\) gives \(0\) in these two integrals,
except at \(x=\ceig{\bplus{r}}\) when \(r=s\). There the residue is the derivative of
\(x\,\Yfn{\bl}(1/x)^{-1}\), which is \(-1/\cores{\bl}{r}\) by Lemma~\ref{lem:residues}.
\end{proof}

For the corner measure, the monomials \(1,x,\dots,x^{m-1}\) have the Hankel Gram matrix
\begin{equation}\label{eq:Gram}
 \Gram{\bl}:=(\mom{i+j}{\bl})_{0\leq i,j<m}.
\end{equation}

\begin{proof}[Proof of Lemma~\ref{lem:shiftip}]
Fix \(\bl\) and put \(w_r:=\pt{\bl-r}{r}\) for \(r\in\Rem{\bl}\).
By Lemma~\ref{lem:fibre}(ii), \(\Shift w_r\) lies in \(\Wfib{\bl}\) and
corresponds to the function \(x/(x-\ceig{\bplus{r}})\).
Equation~\eqref{eq:fibrepairing} and Lemma~\ref{lem:moments} give
\[
 \ip{\Shift w_r}{\Shift w_s}
 =\srcnorm{\bl}\int
   \frac{x^2}{(x-\ceig{\bplus{r}})(x-\ceig{\bplus{s}})}\,d\cornermeas{\bl}
 =-\delta_{rs}\frac{\srcnorm{\bl}}{\cores{\bl}{r}}.
\]
For the diagonal entries, Lemma~\ref{lem:normratio} and
Proposition~\ref{prop:euler} give
\[
 -\frac{\srcnorm{\bl}}{\cores{\bl}{r}}
 =\Mq\,\frac{\srcnorm{\bl-r}}{\res{\bl-r}{r}}
 =\Mq\,\ip{w_r}{w_r}.
\]
Distinct \(w_r\) are orthogonal, so
\(\ip{\Shift w_r}{\Shift w_s}=\Mq\,\ip{w_r}{w_s}\).
Every fixed-point class has the form \(w_r\) for a unique pair \((\bl,r)\).
For two such classes with different \(\bl\), both their pairing and the pairing
of their images vanish, since different fixed-point classes and different
summands \(\Wfib{\bl}\) are orthogonal. Bilinearity now proves the first assertion.

For \(0\leq j<m\), the other identity in Lemma~\ref{lem:moments} gives
\[
 \ip{\taut^{-j}\pulm\fix{\bl}}{\Shift w_r}
 =\srcnorm{\bl}\int x^j\,\frac{x}{x-\ceig{\bplus{r}}}\,d\cornermeas{\bl}=0.
\]
Pairings with shift images in other summands also vanish. This proves the
orthogonality assertion.
\end{proof}

\begin{proof}[Proof of Theorem~\ref{thm:fock}(ii)]
By Lemmas~\ref{lem:shiftip} and~\ref{lem:fibre}(ii), classes with different powers of \(z^{-1}\) are
orthogonal, and
\[
 \ip{\Shift^k\taut^{-i}\pulm\fix{\bl}}{\Shift^k\taut^{-j}\pulm\fix{\bmu}}
 =\delta_{\bl\bmu}\,\Mq^k\,\srcnorm{\bl}\,\mom{i+j}{\bl}.
\]
The pairing is nondegenerate, being the pullback under the isomorphism \(\Omap\) of a diagonal
pairing with nonzero entries.
\end{proof}

On \(\La^{\otimes m}\) we use the pairing \(\ip{\kket{\bl}}{\kket{\bmu}}=\delta_{\bl\bmu}\srcnorm{\bl}\),
which is that of Theorem~\ref{thm:fock}(ii) on \(\La^{\otimes m}\otimes\auxe{0}\).

\begin{lemma}\label{lem:sourceadjoint}
The adjoint of \(\rhou(e)\), with \(e\in\mathcal E\) homogeneous and \(\rhou(e)\) raising the
degree by \(d\), is \(\Mt^{-d}\rhou(\sigH(e))\), where \(\sigH\) is the horizontal
anti-involution of \cite[(2.21)]{BCS}.
\end{lemma}

\begin{proof}
The pairing \cite[(2.25)]{BCS} on \(\Fock{m}\) satisfies \(\rhov(e)^\dagger=\rhov(\sigV(e))\), with
\(\sigV\) the vertical anti-involution of \cite[(2.21)]{BCS}, and \(\srcnorm{\bl}\) is
its norm rescaled by \((\Mq\gamma^m)^{|\bl|}\), which multiplies the adjoint of an
operator of degree \(d\) by \((\Mq\gamma^m)^d\). The two anti-involutions satisfy
\(\Miki\sigV\Miki^{-1}=\sigH\), since both sides agree on \(a_{\pm1}\), \(x^\pm_0\) and the
central elements \cite[(2.9) and (2.21)]{BCS}. Since \(\Hiso\) is an isometry and
intertwines \(\rhov(e)\) with \(\kappam{m}^{-d}\rhou(\Miki(e))\), the adjoint of
\(\rhou(e)\) is \((\kappam{m}^2\Mq\gamma^m)^d\rhou(\sigH(e))\), and
\(\kappam{m}^2\Mq\gamma^m=\Mt^{-1}\).
\end{proof}

For operators on \(\La^{\otimes m}\), write \(C^*\) for the adjoint under the
pairing of Lemma~\ref{lem:sourceadjoint}. For operators on \(\Hsp{m}\), use
\(\dagger\) for the adjoint under the auxiliary pairing of Theorem~\ref{thm:fock}(ii).
If \(C\) has degree \(n\), put \(\vartheta(C):=\Mq^{-n}C^*\).
Then \(\vartheta\) reverses products, fixes the Cartans,
and \(\vartheta(\rhou(e))=(\Mq\Mt)^{-n}\rhou(\sigH(e))\). For a series
\(C(z)=\sum_{n\geq0}C_nz^n\) with \(C_n\) raising the degree by \(n\), put
\(C^\vartheta(z^{-1})=\sum_{n\geq0}\vartheta(C_n)z^{-n}\).

\subsection{Raising and lowering series}\label{sec:raising}

For \(j\geq0\), the horizontal series of iterated Hecke operators is
\[
 \Pser{j}(z):=\sum_{n\geq1}P_{n,j}\,z^n,\qquad
 P_{n,j}:=\Hiso\,\pushp\,\Shift^{n-1}\,\taut^{-j}\,\pulm\,\Hiso^{-1}.
\]
The coefficient \(P_{n,j}\) raises the Fock degree by \(n\).
For a chain \(\bl=\bl_0\subset\bl_1\subset\dots\subset\bl_d\) adding the boxes
\(a_0,\dots,a_{d-1}\) in this order, put
\[
 \kappa_k(\bl_0,\dots,\bl_d):=\ceig{a_0}^{\,k}\prod_{l=1}^{d-1}
 \Bigl(1-q_3\frac{\ceig{a_{l-1}}}{\ceig{a_l}}\Bigr)^{-1},\qquad
 r(\bl_0,\dots,\bl_d):=\prod_{l<d}\res{\bl_l}{a_l}.
\]
By Propositions~\ref{prop:fixed} and~\ref{prop:pullback} and Lemma~\ref{lem:shiftfixed}, the coefficient of
\(\kket{\bmu}\) in \(P_{d,k}\kket{\bl}\) is
\begin{equation}\label{eq:Pchains}
 \sum\kappa_k(\bl_0,\dots,\bl_d)\,r(\bl_0,\dots,\bl_d),
\end{equation}
the sum over the chains from \(\bl\) to \(\bmu\). By Lemma~\ref{lem:normratio} and the definition of
\(\vartheta\), the coefficient of \(\kket{\bl}\) in \(\vartheta(P_{d,k})\kket{\bmu}\) is
\begin{equation}\label{eq:Pvartheta}
 (-1)^d\sum\kappa_k(\bl_0,\dots,\bl_d)\prod_{l<d}\cores{\bl_{l+1}}{a_l},
\end{equation}
the sum again over the chains from \(\bl\) to \(\bmu\).

\begin{samepage}
\begin{theorem}[The raising series]\label{thm:raising}
For every \(m\):
\begin{enumerate}
\item[(i)] \(\displaystyle 1+\Pser{0}(z)=\exp\Bigl(\sum_{k\geq1}\frac{1-\Mt^k}{k}\,z^k\,
\bigl(\pk{k}{1}+\dots+\pk{k}{m}\bigr)\Bigr)\).
\item[(ii)] \(\Pser{1}(z)=-\gamma\,\Xcur(q_2\gamma^mz)_{>0}\).
\item[(iii)] for \(j\geq0\),
\(\displaystyle\Dq\Pser{j+1}=\frac{1-q_1q_2}{(1-q_1)(1-q_2)}\,[\Cartan{1},\Pser{j}]
-(\Dq\Pser{1})\,\Pser{j}\).
\end{enumerate}
\end{theorem}
\end{samepage}

\begin{proof}[Proof of Theorem~\ref{thm:raising}(iii)]
On \(\Fock{m}\) let \(H\) be multiplication by the class \([\Wc^\vee]\), and on
\(\Wsp{m}\) multiplication by its pullback by \(\pminus\). By \eqref{eq:Wcharacter},
\(H\) acts by \(\Cartan{1}(\bl)\) on \(\fix{\bl}\) and on \(\pt{\bl}{a}\).
The pullbacks of the tautological bundles satisfy
\(\pplus{}^*[\Tc]=\pulm[\Tc]+[\taut]\) in \(K\)-theory, so
\(\pplus{}^*[\Wc]=\pulm[\Wc]-(1-q_1)(1-q_2)\,\taut\) on \(\Hecke{n}\), and similarly on
\(\HeckeTwo{n}\) with \(\tautone\). The projection formula gives
\[
 H\pushp=\pushp(H-\epsilon\,\dtaut),\qquad H\Shift=\Shift(H-\epsilon\,\dtaut),\qquad
 H\pulm=\pulm H,
\]
and \(H\) commutes with \(\taut\). Let \(R(z)=(1-z\Shift)^{-1}\), so that
\(\Pser{j}(z)=z\,\Hiso\pushp R(z)\taut^{-j}\pulm\Hiso^{-1}\). Then
\([H,R(z)]=-\epsilon zR(z)\Shift\dtaut R(z)\), and, since \(1+zR(z)\Shift=R(z)\),
\[
 [\Cartan{1},\Pser{j}(z)]
 =-\epsilon\,z\,\Hiso\pushp R(z)\dtaut R(z)\taut^{-j}\pulm\Hiso^{-1}.
\]
Lemma~\ref{lem:generation}(ii) gives
\(\dtaut R(z)=R(q_3z)\dtaut+zR(q_3z)\dtaut\pulm\pushp R(z)\). Substituting this and
\(R(z)R(q_3z)=(R(z)-q_3R(q_3z))/(1-q_3)\), the right side becomes
\(-\epsilon\bigl(\Dq\Pser{j+1}+(\Dq\Pser{1})\Pser{j}\bigr)/(1-q_3)\). Finally
\(-(1-q_3)/\epsilon=(1-q_1q_2)/((1-q_1)(1-q_2))\).
\end{proof}

\paragraph{The shuffle algebra.}
Let \(\Shuf\) be the shuffle algebra of \cite[Section~4.2]{N1}. Its elements of degree
\(k\) are symmetric rational functions of \(z_1,\dots,z_k\) of the form
\cite[(4.4)]{N1}, and its product \cite[(4.2)]{N1} is built from
\[
 \omega(x)=\frac{(x-1)(x-q_1q_2)}{(x-q_1)(x-q_2)}.
\]
It is generated by the functions \(z^d\) of one variable, \(d\in\mathbb Z\), and spanned
by the products
\(z^{n_1}*\dots*z^{n_k}=\Sym\bigl[z_1^{n_1}\cdots z_k^{n_k}\prod_{i<j}\omega(z_i/z_j)\bigr]\),
where \(\Sym\) is the sum over the permutations of the variables \cite[Section~4.2
and (4.13)]{N1}, \cite{Negutshuffle}. By \cite[Theorem~4.7]{N1}, composed with the
automorphism \(P\mapsto P\cdot(z_1\cdots z_k)^m\) of \(\Shuf\), the algebra \(\Shuf\) acts
on \(\Fock{m}\) with \(z^d\) acting by \(\xplus{d}\) (by \cite[Theorem~4.10]{N1} with
\(k=1\)). We write \(P|_{\Fock{m}}\) for the action of \(P\in\Shuf\).

\begin{lemma}\label{lem:shuffle}
\begin{enumerate}
\item[(i)] Let \(P=\Sym\bigl[K(z_1,\dots,z_k)\prod_{i<j}\omega(z_i/z_j)\bigr]\in\Shuf\),
with \(K=K_0/\prod_{i<k}(1-q_1q_2z_{i+1}/z_i)\) and \(K_0\) a Laurent polynomial. The
coefficient of \(\fix{\bl_k}\) in
\(P|_{\Fock{m}}\fix{\bl}\) is the sum over the chains from \(\bl\) to \(\bl_k\) of
\(K(\bwt{a_{k-1}},\dots,\bwt{a_0})\,r(\bl_0,\dots,\bl_k)\).
\item[(ii)] For \(n\geq1\), \(\pushp\Shift^{n-1}\taut^{-1}\pulm=\mathsf P_n|_{\Fock{m}}\), where
\[
 \mathsf P_n:=(-q_1q_2)^{n-1}\,
 \Sym\Bigl[\frac{z_1^{-1}}{\prod_{i=1}^{n-1}(1-q_1q_2z_{i+1}/z_i)}
 \prod_{i<j}\omega(z_i/z_j)\Bigr]\in\Shuf.
\]
\item[(iii)] For \(n\geq1\), \(z^0*\mathsf P_n-\mathsf P_n*z^0=\epsilon\,\mathsf P_{n+1}\). Hence
\([P_{1,0},P_{n,1}]=\epsilon P_{n+1,1}\).
\end{enumerate}
\end{lemma}

By \cite[(4.18)]{N1}, \(\mathsf P_n\) is a multiple of the generator of bidegree \((n,-1)\) of the
elliptic Hall algebra \cite{Negutshuffle}, and (iii) is a relation of that algebra.
Appendix~\ref{sec:shuffleproof} verifies its normalization.

\begin{proof}[Proof of (i) and (ii)]
(i) \cite[Proposition~5.5]{N1}, after the twist by \((z_1\cdots z_k)^m\), writes the coefficient as
the sum over the chains of \(K(\bwt{a_{k-1}},\dots,\bwt{a_0})\) times
\(\frac{g_{\bl}}{g_{\bl_k}}\prod_{i<j}\omega(\chi_i/\chi_j)\prod_i\Lambda(\Wc|_{\bl_k},\chi_iq_1q_2)\),
where \(\chi_i=\bwt{a_{k-i}}\). In the notation of \cite{N1}, \(g_\bl\) is the
cotangent Euler class at \(\bl\), and \(\Lambda(\mathcal V,u)\) is the
alternating exterior-power polynomial of \(\mathcal V^\vee\).
The vanishing factors are removed as in \cite[Remark~5.6]{N1}. Removing a box of weight
\(\chi\) adds \((1-q_1)(1-q_2)\chi\) to \([\Wc]\), which multiplies
\(\Lambda(\Wc,uq_1q_2)\) by \(\omega(\chi/u)\). Thus the factor multiplying \(K\) is the product along the
chain of the one-box coefficients \(\frac{g_{\bl_l}}{g_{\bl_{l+1}}}
\Lambda(\Wc|_{\bl_{l+1}},\bwt{a_l}q_1q_2)\), which are the
\(\res{\bl_l}{a_l}\) by the case of one box (Section~\ref{sec:fixedpoints}).

(ii) By Proposition~\ref{prop:pullback} and Lemma~\ref{lem:shiftfixed}, the coefficient
of \(\fix{\bl_n}\) in \(\pushp\Shift^{n-1}\taut^{-1}\pulm\fix{\bl}\) is the sum over chains of
\(\ceig{a_0}\prod_{l=1}^{n-1}(1-q_3\bwt{a_l}/\bwt{a_{l-1}})^{-1}\,r(\bl_0,\dots,\bl_n)\).
With \(z_i=\bwt{a_{n-i}}\), the identity
\((1-q_3x)^{-1}=-q_1q_2x^{-1}/(1-q_1q_2x^{-1})\) turns the factor before \(r\) into the kernel of
\(\mathsf P_n\), so (i) gives the claim. The element \(\mathsf P_n\) lies in
\(\Shuf\) by \cite[Section~4.9]{N1}.
\end{proof}

\begin{proof}[Proof of Theorem~\ref{thm:raising}(ii) and~(i)]
Put \(\mathsf p_k=\pk{k}{1}+\dots+\pk{k}{m}\). For every \(\alpha\),
\(\partial^{(\alpha)}_k\) applied to \(\mathsf p_k\) gives
\(q_3^{(\alpha-1)k}+(1-q_3^{-k})\sum_{\beta>\alpha}q_3^{(\beta-1)k}=q_3^{(m-1)k}\).

(ii) The case \(k=1\), with \eqref{eq:Xcur} and \(\gamma^2=q_3\), gives
\([\mathsf p_1,\Xcur(w)]=-w^{-1}\gamma^m(1-q_2)\Xcur(w)\). By
Proposition~\ref{prop:horizontal}(ii), \(P_{1,0}=(1-\Mt)\mathsf p_1\). Since
\(\epsilon=(1-\Mt)(1-\Mq^{-1})\), for every \(n\)
\[
 [P_{1,0},\Xcur_{-n}]=\epsilon\,q_2\gamma^m\,\Xcur_{-n-1}.
\]
The relations \(H\pushp=\pushp(H-\epsilon\,\dtaut)\) and \(H\pulm=\pulm H\) established in the proof of
Theorem~\ref{thm:raising}(iii) give \([P_{1,0},\Hiso\Cartan{1}\Hiso^{-1}]=\epsilon P_{1,1}\), and
\(\Hiso\Cartan{1}\Hiso^{-1}=-\gamma\Xcur_0\). With Lemma~\ref{lem:shuffle}(iii), induction on \(n\)
gives \(P_{n,1}=-\gamma(q_2\gamma^m)^n\Xcur_{-n}\) for \(n\geq1\).

(i) Write \(\mathcal V(z)\) for the right side of (i). On \(\mathcal V(z)\) the
annihilation part of the \(\alpha\)-th summand of \(\Xcur(w)\) acts by the factor
\(\exp\bigl(\sum_k(1-\Mq^k)(1-\Mt^k)s^k/k\bigr)\), \(s=\gamma^mz/w\), by the first paragraph. So
\(\Xcur(w)\,\mathcal V(z)=f(s)\,\mathcal V(z)\,\Xcur(w)\) with
\(f(s)=(1-\Mq s)(1-\Mt s)/((1-s)(1-\Mq\Mt s))\). Let
\(c=(1-q_1q_2)/((1-q_1)(1-q_2))\). Then
\[
 c\,[\Cartan{1},\mathcal V]\mathcal V^{-1}
 =c\,\gamma\,[w^0]\bigl((f(s)^{-1}-1)\Xcur(w)\bigr)
 =-\gamma\sum_{n\geq1}(\Mq^n-\Mt^n)(\gamma^mz)^n\,\Xcur_{-n}=\Dq\Pser{1}.
\]
The first equality uses \(\Hiso\Cartan{1}\Hiso^{-1}=-\gamma\Xcur_0\).
For the second, expand
\(1-f(s)^{-1}=(1-\Mq)(1-\Mt)\sum_{n\geq1}\frac{\Mq^n-\Mt^n}{\Mq-\Mt}s^n\)
and use \(c(1-\Mq)(1-\Mt)=\Mq-\Mt\).
Coefficient extraction gives \([w^0](s^n\Xcur(w))=(\gamma^mz)^n\Xcur_{-n}\).
The last equality follows from part~(ii) and \(\Mq^n-\Mt^n=\Mq^n(1-q_3^n)\).

Thus \(Y=\mathcal V\) satisfies \(c[\Cartan{1},Y]=(\Dq\Pser{1})Y\).
The series \(Y=1+\Pser{0}\) satisfies the same equation by
part~(iii) with \(j=0\). Both series have constant term \(1\). If the series differed, the
lowest coefficient of their difference would raise the degree and commute with \(\Cartan{1}\). But the
values \eqref{eq:Hone} of \(\Cartan{1}\) on the \(\kket{\bl}\) are distinct, since the contents of
distinct boxes are distinct monomials in the parameters.
An operator of positive degree commuting with \(\Cartan{1}\) must therefore
vanish, so \(\mathcal V=1+\Pser{0}\).
\end{proof}

Put \(\Qcr(z):=1+\Pser{0}(z)\) and \(\Qan:=1+\Pser{0}^\vartheta(z^{-1})\). For \(m=1\) and \(v_1=1\) these are
the operators of Section~\ref{sec:rankone}.

\begin{corollary}\label{cor:Qdagger}
\[
 \Qan=\exp\Bigl(\sum_{k\geq1}(1-\Mq^k)\,\Mq^{-k}z^{-k}\sum_{\alpha=1}^mq_3^{(\alpha-1)k}\,\partial_{\alpha,k}\Bigr),
 \qquad
 \Pser{1}^\vartheta(z^{-1})=-\gamma\,\Xcur(q_2\gamma^mz/q_3^{m-1})_{<0}.
\]
\end{corollary}

\begin{proof}
This is Theorem~\ref{thm:raising}(i) and (ii) with Lemma~\ref{lem:sourceadjoint}, by
\(\sigH(a_{-k})=-a_k\), \(\sigH(x^-_{-n})=x^-_n\) \cite[(2.21)]{BCS} and \cite[(2.50)]{BCS}.
\end{proof}

\section{Explicit forms of the Lax operator}\label{sec:laxforms}

We prove the two Lax forms in the reverse order to their statements in the introduction.
The q-Miura conjugacy of Theorem~\ref{thm:companion} comes first.
The last-column identity \eqref{eq:Um} from that proof, together with the
auxiliary pairing, gives the anti-normal factorization
of Theorem~\ref{thm:factorization}.

\subsection{The q-Miura form}\label{sec:qmiuraform}

\subsubsection{Fused currents}\label{sec:qmiura}

Write the horizontal current \eqref{eq:Xcur} as \(\Xcur(w)=\sum_\alpha u_\alpha^{-1}\eta^-_\alpha(w)\).
So \(\eta^-_\alpha(w)\) is the exponential in the power sums of colour \(\alpha\), which raises the
degree, times the exponential in the derivatives \(\partial^{(\alpha)}_k\), which lowers the degree. In rank
one \(\eta^-_1\) is the vertex operator \(\eta^-\) of \cite[(2.34)]{BCS}. A normal-ordered product
\({:}\eta^-_{\alpha_1}(w_1)\cdots\eta^-_{\alpha_k}(w_k){:}\) puts all the exponentials in the power sums
to the left of all the exponentials in the derivatives.
The normal-ordered product is unchanged by permuting its factors. Put
\[
 \phi(x):=\frac{(1-x/q_1)(1-x/q_2)}{(1-x)(1-q_3x)}
 =1+\frac{\epsilon}{1-q_3}\Bigl(\frac{1}{1-x}-\frac{1}{1-q_3x}\Bigr).
\]

\begin{lemma}\label{lem:normalorder}
With \(x=w/z\), and the right sides expanded in powers of \(x\),
\[
 \eta^-_\alpha(z)\,\eta^-_\beta(w)=
 \begin{cases}
 \phi(x)^{-1}\,{:}\eta^-_\alpha(z)\eta^-_\beta(w){:},&\alpha=\beta,\\
 \phi(1/x)\phi(x)^{-1}\,{:}\eta^-_\alpha(z)\eta^-_\beta(w){:},&\alpha<\beta,\\
 {:}\eta^-_\alpha(z)\eta^-_\beta(w){:},&\alpha>\beta.
 \end{cases}
\]
\end{lemma}

\begin{proof}
Move the exponential in the derivatives of \(\eta^-_\alpha(z)\) to the right of the exponential
in the power sums of \(\eta^-_\beta(w)\). The logarithm of the resulting scalar factor is
\(-\sum_kq_3^k(1-q_1^k)(1-q_2^k)x^k/k\) for \(\alpha=\beta\).
For \(\alpha<\beta\), the coefficient of \(x^k\) has the additional factor \(1-q_3^{-k}\).
For \(\alpha>\beta\), the derivatives of colours \(\geq\alpha\) annihilate
polynomials in the \(\beta\)-th alphabet.
\end{proof}

\begin{definition}[Fused currents]\label{def:fused}
For \(k\geq1\), the fused current \(\Xk{k}(w)=\sum_{n\in\mathbb Z}\Xk{k}_n\,w^{-n}\) is
\[
 \Xk{k}(w):=\sum_{\alpha_1<\dots<\alpha_k}(u_{\alpha_1}\cdots u_{\alpha_k})^{-1}\,
 {:}\eta^-_{\alpha_1}(w)\,\eta^-_{\alpha_2}(q_3w)\cdots\eta^-_{\alpha_k}(q_3^{k-1}w){:}.
\]
\end{definition}

So \(\Xk{1}=\Xcur\), \(\Xk{k}=0\) for \(k>m\), and the mode \(\Xk{k}_n\) lowers the degree by \(n\). The
arguments \(w,q_3w,\dots,q_3^{k-1}w\) are assigned in increasing colour order. The same products appear in
\cite[Section~3.1]{AFHKSY} and \cite{FHSSY}.

Write \(\DD=\Twist{q_3}\), so that \(\DD C(z)=C(z/q_3)\DD\) for a series \(C\), and \(\pi=\pile\). In the
projected q-Miura operator of Section~\ref{sec:mainresults},
\[
 \Xhat{k}=-\gamma^k\,\DD^k\,\pi\,\Xk{k}(q_2\gamma^mz)=-\gamma^k\,\pi\,\Xk{k}(q_2\gamma^mz/q_3^k)\,\DD^k,
\]
since \(\DD\) commutes with \(\pi\). The projection makes \(\Xhat{k}\) an operator on
the auxiliary polynomial space \(V\) of Section~\ref{sec:mainresults}. Multiplying out
the product below, moving the powers of \(\DD\) to the left, normal ordering, and using
\(u_\alpha=-\gamma v_\alpha\), gives
\[
 \Miuraop(y)={:}\bigl(1-y\,v_m^{-1}\DD\,\eta^-_m(q_2\gamma^mz)\bigr)\cdots
 \bigl(1-y\,v_1^{-1}\DD\,\eta^-_1(q_2\gamma^mz)\bigr){:}.
\]
This is the q-Miura transformation that defines the currents of the deformed \(W\)-algebra of
type \(\mathfrak{gl}_m\) \cite{FR,FF,AKOS}, \cite[(5.3)]{N3}. Section~\ref{sec:negut} identifies the
fused currents with the currents of \cite{N3}.

\begin{lemma}\label{lem:corner}
\(\betaz\,\Qcr(z)\,\Qan=-\gamma^m\Xk{m}(q_2\gamma^mz/q_3^{m-1})\), with \(\Qan\) acting first.
\end{lemma}

\begin{proof}
Both sides are an exponential in the power sums times an exponential in the
derivatives, and \(-\gamma^m(u_1\cdots u_m)^{-1}=\betaz\) since \(u_\alpha=-\gamma v_\alpha\). Put
\(w=q_2\gamma^mz/q_3^{m-1}\). Then \(\Xk{m}(w)\) is \((u_1\cdots u_m)^{-1}\) times
the normal-ordered product of the \(\eta^-_\alpha(q_3^{\alpha-1}w)\). In the power sums of colour
\(\alpha\) the factors \(q_3^{\pm(\alpha-1)k}\) of \cite[(2.52)]{BCS} cancel. In the derivatives,
the coefficient of \(\partial_{\beta,k}\) in \(\sum_\alpha q_3^{-(\alpha-1)k}\partial^{(\alpha)}_k\) is
\(1+(1-q_3^{-k})\sum_{s=1}^{\beta-1}q_3^{sk}=q_3^{(\beta-1)k}\). Since \(\gamma^2=q_3=(q_1q_2)^{-1}\),
\(q_1=\Mt^{-1}\) and \(q_2=\Mq\),
\[
 -w^k\gamma^{mk}(1-q_1^k)=(1-\Mt^k)z^k,\qquad
 w^{-k}q_3^k\gamma^{-mk}(1-q_2^k)=(1-\Mq^k)\Mq^{-k}z^{-k},
\]
so the two exponents are those of \(\Qcr(z)\) in Theorem~\ref{thm:raising}(i) and of
\(\Qan\) in Corollary~\ref{cor:Qdagger}.
\end{proof}

\subsubsection{\Negut's factorization}\label{sec:negut}

The raising and lowering factors of \Negut{} \cite[Theorem~1.2]{N3} are expressed
in corner contents. For the current comparison in this subsection, we therefore
use the inverted geometric parameters, in place of the direct identification
used for the pairing in Section~\ref{sec:form}.
The parameters \((q_1,q_2,q)\) of \cite{N3} are now our \((\Mq^{-1},\Mt,q_3)\).
The constant \((1-q_1)(1-q_2)/(1-q)\) of \cite{N3} becomes
\(\epsilon/(1-q_3)\), and the framing weights become
\(v_1^{-1},\dots,v_m^{-1}\). The weight \(\chi_b\) of a box \(b\) \cite[(3.5)]{N3} is its content
\(\ceig{b}\), and the function \(\zeta\) of \cite[(2.1)]{N3} is the function
\(\phi(x)=(1-x/\Mq)(1-\Mt x)/((1-x)(1-q_3x))\) of Section~\ref{sec:qmiura}, which is
\(\omega(x)^{-1}\). The function \(\tau(z)\) of \cite[(3.8)]{N3} is \(\prod_\alpha(1-v_\alpha z)\), and
for every \(m\)-tuple \(\bl\)
\begin{equation}\label{eq:negutY}
 \tau(z)\prod_{b\in\bl}\phi(\ceig{b}/z)=\Yfn{\bl}(1/z),
\end{equation}
since \(\Yfn{\varnothing}(1/z)=\prod_\alpha(1-v_\alpha z)\) and
\(\Yfn{\bl+b}(w)=\Yfn{\bl}(w)\,\omega(w/\bwt{b})^{-1}\) by \eqref{eq:Wcharacter}.

Let \(K\) be the localized \(K\)-theory of the spaces \(\Gies{m}{n}\), with the basis \(|\bl\rangle\)
of fixed points of \cite[(3.6)]{N3}, and write \(\langle\bmu|R|\bl\rangle\) for the coefficient of
\(|\bmu\rangle\) in \(R|\bl\rangle\). The currents \(W_k(x)=\sum_{d\in\mathbb Z}W_{d,k}x^{-d}\) of
\cite[(2.48)]{N3} act on \(K\), with \(W_0=1\) and \(W_k=0\) for \(k>m\) \cite[(3.25)]{N3}, and
\(W_{d,k}\) changes the number of boxes by \(-d\).
Let \(D_x\) be the dilation \(f(x)\mapsto f(xq)\).

\begin{theorem}[\Negut{} \cite{N3}]\label{thm:negutfactorization}
The currents on \(K\) satisfy the identity of formal difference operators
\cite[Theorem~1.2]{N3}
\begin{equation}\label{eq:LDU}
 \sum_{k\geq0}W_k(x)\,(-yD_x)^{-k}
 =T(x^{-1},yD_x)^{\leftarrow}\,E(yD_x)\,T(xq,yD_x)^{\rightarrow},
\end{equation}
where \(E(y)=\sum_{k\geq0}E_{0,k}(-y)^{-k}\) and
\[
 T(x^{-1},yD_x)^{\leftarrow}=1+\sum_{d,k\geq1}x^d(-yD_x)^{-k}\,T^{\leftarrow}_{d,k},\qquad
 T(xq,yD_x)^{\rightarrow}=1+\sum_{d,k\geq1}(-yD_x)^{-k}(xq)^{-d}\,T^{\rightarrow}_{d,k}.
\]
The powers of \(D_x\) are placed after the powers of \(x\) on the left side and in the first
factor, and before them in the last factor \cite[Remark~2.17]{N3}. The operator
\(T^{\leftarrow}_{d,k}\) adds \(d\) boxes, \(E_{0,k}\) is diagonal, and \(T^{\rightarrow}_{d,k}\)
removes \(d\) boxes \cite[(2.50)--(2.55)]{N3}.
\end{theorem}

\begin{proposition}\label{prop:negutdictionary}
Let \(\Xi\colon K\to\La^{\otimes m}\) be the linear isomorphism with
\(\Xi|\bl\rangle=(-\Mq)^{|\bl|}\srcnorm{\bl}^{-1}\kket{\bl}\). For \(d,k\geq1\),
\begin{enumerate}
\item[(i)] \(\Xi\,T^{\leftarrow}_{d,k}\,\Xi^{-1}=(-1)^{k-1}P_{d,k}\).
\item[(ii)] \(\Xi\,T^{\rightarrow}_{d,k}\,\Xi^{-1}=(-1)^{k-1}q_3^{(1-m)d}\,\vartheta(P_{d,k})\).
\item[(iii)] \(\Xi\,E_{0,k}\,\Xi^{-1}=\Newton{k}\).
\item[(iv)] \(\Xi\,W_k(x)\,\Xi^{-1}=(-\gamma)^k\,\Xk{k}(q_2\gamma^mx/q_3^{k-1})\) for every \(k\geq1\).
\end{enumerate}
\end{proposition}

\begin{proof}
(i), (ii) The differential \(\omega(u)^{-1}du/u\) has residues \(-\epsilon/(1-q_3)\) at \(u=1\) and
\(\epsilon/(1-q_3)\) at \(u=q_1q_2\). Lemma~\ref{lem:residues} and
\eqref{eq:Wcharacter} therefore give the one-box conversions
\[
 \frac{\epsilon}{1-q_3}\,\Yfn{\bl}(\bwt{\bplus{a}})=\cores{\bl+a}{a},\qquad
 \frac{\epsilon}{1-q_3}\,\Yfn{\bl+a}(\bwt{a})^{-1}=-\res{\bl}{a}.
\]
By \cite[(3.21), (3.22)]{N3}, the matrix coefficients of
\(T^{\leftarrow}_{d,k}\) and \(T^{\rightarrow}_{d,k}\) are sums over chains
\(\bl=\bl_0\nearrow\cdots\nearrow\bl_d=\bmu\), with
\(\bl_{l+1}=\bl_l+a_l\). The box with label \(i\) in \cite{N3} is \(a_{d-i}\).
Each chain contribution to either matrix coefficient contains
\[
 (-1)^{k-1}\kappa_k\prod_{l<l'}\phi(\ceig{a_{l'}}/\ceig{a_l}),
\]
where \(\kappa_k\) is the chain factor in \eqref{eq:Pchains}.
For each \(l\), the remaining factor is
\[
 \frac{\epsilon}{1-q_3}\tau(q_3\ceig{a_l})
   \prod_{b\in\bl}\phi(\ceig{a_l}/\ceig{b})
 \quad\text{or}\quad
 \frac{\epsilon q_3^{1-m}}{1-q_3}\Yfn{\bmu}(\bwt{a_l})^{-1},
\]
respectively. In the first expression,
\(\phi(x)=\phi(1/(q_3x))\) \cite[(2.2)]{N3} and \eqref{eq:negutY}
identify the product with
\(\epsilon(1-q_3)^{-1}\Yfn{\bl}(\bwt{\bplus{a_l}})\).
Multiplying by the factors with \(l'<l\) replaces \(\bl\) by \(\bl_l\).
In the second expression, the identity
\[
 \Yfn{\bmu}(\bwt{a_l})
 =\Yfn{\bl_{l+1}}(\bwt{a_l})
   \prod_{l'>l}\phi(\ceig{a_{l'}}/\ceig{a_l})
\]
cancels the factors with \(l'>l\). The one-box conversions now give
\begin{align*}
 \langle\bmu|T^{\leftarrow}_{d,k}|\bl\rangle&=(-1)^{k-1}\sum\kappa_k\prod_{l<d}\cores{\bl_{l+1}}{a_l},\\
 \langle\bl|T^{\rightarrow}_{d,k}|\bmu\rangle&=(-1)^{k-1}\bigl(-q_3^{1-m}\bigr)^d\sum\kappa_k\prod_{l<d}\res{\bl_l}{a_l}.
\end{align*}
For an operator \(R\) on \(K\), the coefficient of \(\kket{\bmu}\) in \(\Xi R\,\Xi^{-1}\kket{\bl}\) is
\((-\Mq)^{|\bmu|-|\bl|}(\srcnorm{\bl}/\srcnorm{\bmu})\langle\bmu|R|\bl\rangle\). Along a chain,
\(\prod_l\res{\bl_l}{a_l}=(-\Mq)^d(\srcnorm{\bl}/\srcnorm{\bmu})\prod_l\cores{\bl_{l+1}}{a_l}\) by
Lemma~\ref{lem:normratio}. Now compare \eqref{eq:Pchains} and \eqref{eq:Pvartheta}.

(iii) By \cite[(2.55) and (3.23)]{N3} and \eqref{eq:negutY}, on \(|\bl\rangle\) the operator \(E(y)\) acts by
\(-\betaz y^{-m}\Yfn{\bl}(1/y)\), which is \(\sum_k\Newton{k}(\bl)(-y)^{-k}\) by \eqref{eq:Yexpansion}.

(iv) For \(k=1\), \cite[(2.50)]{N3} gives
\(W_1(x)=\sum_{d\geq1}T^{\leftarrow}_{d,1}x^d+E_{0,1}+\sum_{d\geq1}T^{\rightarrow}_{d,1}x^{-d}\). By (i)
and Theorem~\ref{thm:raising}(ii), \(\Xi T^{\leftarrow}_{d,1}\Xi^{-1}=P_{d,1}=-\gamma(q_2\gamma^m)^d\Xcur_{-d}\).
By (iii), \(\Xi E_{0,1}\Xi^{-1}=\Newton{1}=\Cartan{1}=-\gamma\Xcur_0\). By (ii) and
Corollary~\ref{cor:Qdagger},
\(\Xi T^{\rightarrow}_{d,1}\Xi^{-1}=q_3^{(1-m)d}\vartheta(P_{d,1})
=-\gamma(q_2\gamma^m)^{-d}\Xcur_d\). So
\(\Xi W_1(x)\Xi^{-1}=-\gamma\Xcur(q_2\gamma^mx)\).

By Lemma~\ref{lem:normalorder}, with \(\phi=\zeta\), the vertex operators
\(v_\alpha^{-1}\eta^-_\alpha(q_2\gamma^mx)\) have the normal-ordering relations
of the vertex operators of colour \(m+1-\alpha\) in \cite[proof of Proposition~5.4]{N3}, so their
modes satisfy \cite[(5.1)]{N3} with the colours reversed. Since \(-\gamma u_\alpha^{-1}=v_\alpha^{-1}\),
the current \((-\gamma)^k\Xk{k}(q_2\gamma^mx/q_3^{k-1})\) is the current \(W_k\) of \cite[(5.4)]{N3}
built from them. These currents satisfy the relation \cite[(2.62)]{N3} by
\cite[Proposition~5.4]{N3} and \cite{Odake}, and the currents \(W_k\) on \(K\) satisfy the same relation by
\cite[Proposition~2.20 and Theorem~3.12]{N3}. For either family, the case
\(k'=1\), written in our parameters, is
\begin{align*}
 &W_k(x')W_1(x)\phi(q_3^{k-1}x/x')
   -W_1(x)W_k(x')\phi(x'/x)\\
 &\qquad=\frac{\epsilon}{1-q_3}
 \left[
 \delta\!\left(\frac{x}{q_3x'}\right)W_{k+1}(x)
 -\delta\!\left(\frac{x'}{q_3^kx}\right)W_{k+1}(x')
 \right].
\end{align*}
The first scalar factor is expanded in nonnegative powers of \(x/x'\),
and the second in nonnegative powers of \(x'/x\), as in
\cite[Remark~2.21]{N3}. Each coefficient then acts by a finite sum.
Multiplication by \(1-x'/(xq_3^k)\) annihilates the second delta-function.
The coefficient of \(x'^0\) on the right is
\[
 \frac{\epsilon(1-q_3^{-k-1})}{1-q_3}\,W_{k+1}(x),
\]
while the left side depends only on \(W_k\) and \(W_1\).
The scalar coefficient of \(W_{k+1}(x)\) is nonzero over \(\K\). The equality for \(k=1\) therefore
determines the equality for every \(k\) by induction, as in
\cite[(4.32)--(4.33)]{N3}.
\end{proof}

For \(k\geq1\), define the lowering series
\begin{equation}\label{eq:Rk}
 \Rlast{k}(z):=(-1)^{k+1}\Newton{k}+\sum_{j=1}^k(-1)^{k+j}\,\Newton{k-j}\,
 \Pser{j}^\vartheta\bigl((q_3^{m-1-k}z)^{-1}\bigr),
\end{equation}
whose coefficients act on \(\La^{\otimes m}\), with the Newton polynomials on the left.

\begin{proposition}\label{prop:WH}
\(\Miuraop(y)=\bigl(1-\sum_ky^k\Pser{k}(z/q_3)\DD^k\bigr)\bigl(1-\sum_ky^k\Rlast{k}(z)\DD^k\bigr)\). That is,
for every \(k\geq1\), with \(\Xk{k}=0\) for \(k>m\),
\begin{equation}\label{eq:WH}
 -\gamma^k\Xk{k}(q_2\gamma^mz/q_3^k)=\Pser{k}(z/q_3)+\Rlast{k}(z)
 -\sum_{j=1}^{k-1}\Pser{j}(z/q_3)\,\Rlast{k-j}(z/q_3^j).
\end{equation}
\end{proposition}

\begin{proof}
Conjugate the factorization of Theorem~\ref{thm:negutfactorization} by \(\Xi\),
put \(x=z/q_3\) and replace \(y\) by \(y^{-1}\). Then
\(D_x=\DD^{-1}\), \(xq=z\), and \((-yD_x)^{-k}\) becomes \((-y)^k\DD^k\). By
Proposition~\ref{prop:negutdictionary}(iv) the left side becomes \(\Miuraop(y)\). By (i), the first
factor becomes \(1-\sum_ky^k\Pser{k}(z/q_3)\DD^k\). By (ii) and (iii), and since
\(\DD C(z)=C(z/q_3)\DD\), the product of the other two factors becomes
\[
 \Bigl(\sum_{j\geq0}\Newton{j}(-y)^j\DD^j\Bigr)
 \Bigl(1-\sum_{k\geq1}y^k\DD^k\,\Pser{k}^\vartheta\bigl((q_3^{m-1}z)^{-1}\bigr)\Bigr)
 =1-\sum_{k\geq1}y^k\Rlast{k}(z)\DD^k.
\]
Comparing the coefficients of \(y^k\) and cancelling the common right factor
\(\DD^k\) gives \eqref{eq:WH}.
\end{proof}

Proposition~\ref{prop:negutdictionary}(iv) identifies the geometric currents of
\cite{N3} with the free-field currents. The intertwiner uses the BCS normalization
of the generalized Macdonald basis. The intertwiner in \cite[Claim~5.9]{N3} is defined
through the stable basis. That claim is stated without proof, and the
identification of the two intertwiners remains separate from our calculation.

\subsubsection{Conjugacy to the companion matrix}\label{sec:proof}\label{sec:companion}

The columns of the map \(\Omap\Vconj{m}\) satisfy a companion recurrence
under multiplication by \(\dtaut\). The displacement relation determines the
successive columns, while interpolation and the current factorization close
the recurrence after \(m\) steps.

For \(j\geq0\) let \(U_j\colon V\to\Wsp{m}\) be the map
\(U_j(z^{-k}f)=\Shift^k\taut^{-j}\pulm\Hiso^{-1}f\). So \(U_jz^{-1}=\Shift U_j\), and
\(\Omap(\xi\otimes\auxe{j})=U_j\xi\) for \(j<m\).

\begin{lemma}\label{lem:move}
For \(j\geq0\), \(\dtaut\,U_j=U_{j+1}\DD+U_1\DD\pi\,\Pser{j}(z)\) on \(V\). For \(j=0\) this reads
\(\dtaut\,U_0=U_1\DD\pi\,\Qcr(z)\).
\end{lemma}

\begin{proof}
Apply \eqref{eq:iterated} to \(\taut^{-j}\pulm\Hiso^{-1}f\):
\[
 \dtaut\,U_j(z^{-k}f)=q_3^k\,U_{j+1}(z^{-k}f)
 +\sum_{i=0}^{k-1}q_3^{\,i}\,U_1\bigl(z^{-i}P_{k-i,j}f\bigr).\qedhere
\]
\end{proof}

Interpolation on each fixed-point fibre expresses \(U_{m+1}\) in terms of
\(U_1,\dots,U_m\). The coefficients are obtained by expanding
\(\taut^{-m-1}\pulm\fix{\bl}\) and the fixed-point classes in the powers
\(\taut^{-j}\pulm\fix{\bl}\), \(1\leq j\leq m\), and the image of \(\Shift\).
For an \(m\)-tuple \(\bmu\) and a number \(c\), put
\[
 h_j(\bmu,c):=\betaz^{-1}\sum_{r+s=m-j}(-1)^r\Newton{r}(\bmu)\,c^s.
\]

\begin{lemma}\label{lem:expansion}
For every \(\bl\), and for every \(\bmu\) and \(a\in\Add{\bmu}\),
\begin{gather*}
 \taut^{-m-1}\pulm\fix{\bl}=\sum_{k=1}^m(-1)^{k+1}\Newton{k}(\bl)\,\taut^{-(m+1-k)}\pulm\fix{\bl}
 -\betaz\,\Shift\sum_{r\in\Rem{\bl}}\cores{\bl}{r}\,\ceig{\bplus{r}}\,\pt{\bl-r}{r},\\
 \pt{\bmu}{a}=\sum_{j=1}^mh_j(\bmu,\ceig{a})\,\taut^{-j}\pulm\fix{\bmu}
 +\Shift\sum_{r\in\Rem{\bmu}}\cores{\bmu}{r}\,\frac{\ceig{\bplus{r}}}{\ceig{\bplus{r}}-\ceig{a}}\,\pt{\bmu-r}{r}.
\end{gather*}
\end{lemma}

\begin{proof}
By Lemma~\ref{lem:fibre}, on the contents of the addable corners of an \(m\)-tuple \(\boldsymbol\nu\), the class
\(\taut^{-j}\pulm\fix{\boldsymbol\nu}\) is the function \(x^j\), \(\Shift\,\pt{\boldsymbol\nu-r}{r}\) is
\(x/(x-\ceig{\bplus{r}})\), and \(\pt{\boldsymbol\nu}{a}\) has the value \(\delta_{ab}/\res{\boldsymbol\nu}{b}\) at
\(\ceig{b}\). So the two left sides are the restrictions of \(g(x)=x^{m+1}+\betaz x\Yfn{\bl}(1/x)\), for
\(\boldsymbol\nu=\bl\), and of \(F(x)=x\Yfn{\bmu}(1/x)/(\ceig{a}-x)\), for \(\boldsymbol\nu=\bmu\). Both vanish at
\(x=0\), have simple poles at the \(\ceig{\bplus{r}}\), and have polynomial parts
of degree at most \(m\). Thus each function is
\(\sum_{j=1}^mp_jx^j+\sum_rt_r\,x/(x-\ceig{\bplus{r}})\). The \(p_j\) come from \eqref{eq:Yexpansion}. By
Lemma~\ref{lem:residues}, \(\Res_{x=\ceig{\bplus{r}}}\Yfn{\boldsymbol\nu}(1/x)=-\ceig{\bplus{r}}\cores{\boldsymbol\nu}{r}\),
and the coefficient \(t_r\) for either \(g\) or \(F\) is its residue divided by
\(\ceig{\bplus{r}}\).
\end{proof}

Substitute the fixed-point expansion of Lemma~\ref{lem:expansion} repeatedly
into the expansion of \(\taut^{-m-1}\pulm\fix{\bl}\).
The coefficient of \(\Shift^d\), \(d\geq1\), is the
sum over the chains \(\bl_0\subset\dots\subset\bl_d=\bl\) of
\[
 -\betaz\,\ceig{\bplus{a_{d-1}}}\prod_{l=1}^{d-1}\frac{\ceig{\bplus{a_{l-1}}}}{\ceig{\bplus{a_{l-1}}}-\ceig{a_l}}
\prod_{l<d}\cores{\bl_{l+1}}{a_l}\ \sum_{j=1}^mh_j(\bl_0,\ceig{a_0})\,\taut^{-j}\pulm\fix{\bl_0}.
\]

Since \(\ceig{\bplus{a}}=q_3\ceig{a}\), the factors before \(\prod_{l<d}\) multiply to
\((-q_3)^d\betaz\,\kappa_1\), and \(\betaz\,c\,h_{m+1-k}(\bl_0,c)\) is
\(\sum_{j=1}^k(-1)^{k-j}\Newton{k-j}(\bl_0)\,c^j\). So by \eqref{eq:Pvartheta} and \eqref{eq:Rk}, the coefficient
of \(\Shift^d\taut^{-(m+1-k)}\pulm\fix{\bl_0}\) is the coefficient of \(z^{-d}\kket{\bl_0}\) in
\(\Rlast{k}(q_3^{k-m}z)\kket{\bl}\), where \(\Pser{j}^\vartheta\) is evaluated at \((z/q_3)^{-1}\). The terms
with \(d=0\) are the \(\Newton{k}\) in the expansion of \(\taut^{-m-1}\pulm\fix{\bl}\).

Consequently, for \(f\in\La^{\otimes m}\),
\[
 \taut^{-m-1}\pulm\Hiso^{-1}f=\sum_{k=1}^mU_{m+1-k}\bigl(\Rlast{k}(q_3^{k-m}z)f\bigr).
\]
Since each \(\Rlast{k}\) is a series in \(z^{-1}\), the equality extends to the identity of maps on \(V\)
\begin{equation}\label{eq:Eprime}
 U_{m+1}=\sum_{k=1}^mU_{m+1-k}\,\Rlast{k}(q_3^{k-m}z).
\end{equation}

For \(j\geq1\) put
\[
 \Theta_j=U_j\DD^{j-1}-\sum_{k=1}^{j-1}U_{j-k}\,\DD^{j-k-1}\,\Rlast{k}(z)\,\DD^k,
\]
so that \(\Theta_1=U_1\). Multiplying \eqref{eq:Eprime} by \(\DD^m\) on the right gives
\begin{equation}\label{eq:Thetam}
 \Theta_{m+1}=0.
\end{equation}

\begin{proposition}\label{prop:Theta}
For \(1\leq j\leq m\), \(\dtaut\,\Theta_j=\Theta_{j+1}+U_1\Xhat{j}\).
\end{proposition}

\begin{proof}
Lemma~\ref{lem:move}, applied to each term of \(\Theta_j\), gives \(\Theta_{j+1}\) plus
\[
 U_1\,\pi\Bigl(\Pser{j}(z/q_3)+\Rlast{j}(z)
 -\sum_{k=1}^{j-1}\Pser{j-k}(z/q_3)\,\Rlast{k}(z/q_3^{j-k})\Bigr)\DD^j,
\]
using \(\DD\pi=\pi\DD\) and \(\Rlast{j}=\pi\Rlast{j}\). By \eqref{eq:WH} this is \(U_1\Xhat{j}\).
\end{proof}

Let \(\Vconj{m}\colon V^m\to\Hsp{m}\) be
\begin{equation}\label{eq:Vm}
 \Vconj{m}(\xi_1,\dots,\xi_m)=\betaz\Qan\DD^{m-1}\xi_m\otimes\auxe{0}
 +\sum_{j=1}^{m-1}\Bigl(\DD^{j-1}\xi_j-\sum_{l=j+1}^{m-1}
 \Rlast{l-j}(z/q_3^{j-1})\,\DD^{l-1}\xi_l\Bigr)\otimes\auxe{j}.
\end{equation}
In the source order \(\xi_1,\dots,\xi_{m-1},\xi_m\) and target order
\(\auxe{1},\dots,\auxe{m-1},\auxe{0}\), the map \(\Vconj{m}\) is triangular.
The diagonal entries \(\DD^{j-1}\) and \(\betaz\Qan\DD^{m-1}\) are invertible,
so \(\Vconj{m}\) is invertible. The map preserves the total degree, as do
\(\Miura{m}\) and \(\Lax{m}\).

\begin{proof}[Proof of Theorem~\ref{thm:companion}]
Write \(\xi^{[l]}\in V^m\) for the vector with \(\xi\) in place \(l\) and \(0\) elsewhere. We show that
\begin{equation}\label{eq:OmegaV}
 \Omap\Vconj{m}\,\xi^{[l]}=\Theta_l\,\xi,\qquad1\leq l\leq m.
\end{equation}
For \(l<m\) this is \eqref{eq:Vm} term by term, since
\(\Rlast{k}(z/q_3^{j-1})\DD^{l-1}=\DD^{j-1}\Rlast{k}(z)\DD^{l-j}\). For \(l=m\), Lemma~\ref{lem:corner} and
Lemma~\ref{lem:move} with \(j=0\) give
\[
 \dtaut\,U_0\,\betaz\Qan\DD^{m-1}=U_1\DD\pi\bigl(\betaz\Qcr(z)\Qan\bigr)\DD^{m-1}=U_1\Xhat{m}
 =\dtaut\,\Theta_m,
\]
the last step by Proposition~\ref{prop:Theta} and \eqref{eq:Thetam}. Since \(\dtaut\) is invertible,
\(\Theta_m=U_0\,\betaz\Qan\DD^{m-1}\), which is \eqref{eq:OmegaV} for \(l=m\). Now
\(\Miura{m}\xi^{[l]}=(\Xhat{l}\xi)^{[1]}+\xi^{[l+1]}\), with \(\xi^{[m+1]}=0\), so
Proposition~\ref{prop:Theta} and \eqref{eq:Thetam} are the identity
\(\dtaut\,\Omap\Vconj{m}=\Omap\Vconj{m}\Miura{m}\), and \(\dtaut\,\Omap=\Omap\Lax{m}\). For \(m=1\) and
\(v_1=1\), \(\betaz=1\), and Lemma~\ref{lem:corner} is the statement
\(\eta^-(q_2\gamma z)=\Qcr(z)\Qan\).
\end{proof}

The case \(l=m\) can be written
\begin{equation}\label{eq:Um}
 U_m=U_0\,\betaz\Qan+\sum_{i=1}^{m-1}U_i\,\Rlast{m-i}(z/q_3^{i-1}).
\end{equation}

\subsection{The anti-normal form}\label{sec:construction}

We combine the companion conjugacy with the auxiliary pairing to obtain the
anti-normal factorization. The key calculation is the adjoint of the raising
matrix with respect to the full Hankel pairing.

The raising matrix \(\RaiseFull{m}\) is the \(m\times m\) identity matrix with
\(\Pser{0}(z),\dots,\Pser{m-1}(z)\) added to its first row:
\((\RaiseFull{m})_{ij}=\delta_{ij}+\delta_{i0}\Pser{j}(z)\).
Multiplication by a matrix series here means that the coefficient operators
act on the Fock variables, while the powers of \(z\) act on the auxiliary variable.
Thus \(\RaiseFull{m}\) sends \(\Hsp{m}\) to Laurent series in \(z\) with finite
coefficient sums, and \(\Raise{m}=\pile\RaiseFull{m}\) is an operator on \(\Hsp{m}\).
The operator \(\Comp{m}\) sends
\(z^{-k}\kket{\bl}\otimes v\) to \(z^{-k}\kket{\bl}\otimes\Comp{m}(\bl)v\), where
\begin{equation}\label{eq:companionentries}
 \Comp{m}(\bl)_{ij}=\delta_{i,j+1}\quad(j<m-1),\qquad
 \Comp{m}(\bl)_{i,m-1}=\begin{cases}\betaz,&i=0,\\
 (-1)^{m-i+1}\Newton{m-i}(\bl),&1\leq i<m.\end{cases}
\end{equation}
Let \(\Akernel{m}\) be the matrix of operators with \((\Akernel{m})_{ij}=\delta_{i,j+1}\) for \(j<m-1\),
\((\Akernel{m})_{0,m-1}=\betaz\Qan\), and \((\Akernel{m})_{i,m-1}=\Rlast{m-i}(z/q_3^{i-1})\) for \(1\leq i<m\).
The series \(\Akernel{m}\) acts on \(\Hsp{m}\) by multiplication and commutes with \(z^{-1}\).

\begin{proposition}\label{prop:currentkernel}
On \(\Hsp{m}\), \(\Lax{m}=\Akernel{m}\,\DD\,\Raise{m}\).
\end{proposition}

\begin{proof}
Let \(G=\sum_{j<m}G_j\otimes\auxe{j}\) and \(H=\Raise{m}G\), so that \(H_j=G_j\) for \(j\geq1\). Lemma~\ref{lem:move}
gives \(\dtaut\,\Omap G=\sum_{j<m}\dtaut\,U_jG_j=\sum_{j=1}^mU_j\DD H_{j-1}\). Replacing \(U_m\) by
\eqref{eq:Um}, this is \(\Omap(\Akernel{m}\DD H)\).
\end{proof}

\begin{theorem}[Adjoint of the raising factor]\label{thm:adjointexplicit}
\(\Raise{m}^\dagger\) is multiplication by the matrix series
\[
 G^{-1}\QQ_m^{\vartheta}(z^{-1})^{\mathsf T}G.
\]
Here \(G\) acts by \(\Gram{\bl}\) on
\(\kket{\bl}\otimes\K^m\), and \(\vartheta\) is applied to the entries of \(\RaiseFull{m}\).
The entries of \(\Raise{m}^\dagger\) are \((\Raise{m}^\dagger)_{ij}=\delta_{ij}\) for \(j\geq1\), \((\Raise{m}^\dagger)_{00}=\Qan\), and
\[
 (\Raise{m}^\dagger)_{i0}=\betaz^{-1}\sum_{j=1}^{m-i}(-1)^{i+j+m}\,\Newton{m-i-j}\,
 \Pser{j}^\vartheta(z^{-1})\qquad(1\leq i<m),
\]
with the Newton polynomials acting on the left.
\end{theorem}

\begin{proof}
Let \(M(z)=\sum_{n\geq0}z^nM_n\) be a matrix series whose entries in \(M_n\) raise the degree
by \(n\). By Theorem~\ref{thm:fock}(ii), the pairing on \(z^{-k}\Csec{m}\) is \(\Mq^k\) times
\(\ip{f\otimes\auxe{i}}{g\otimes\auxe{j}}=\ip{f}{G_{ij}g}\), with the pairing of
Lemma~\ref{lem:sourceadjoint} and \(G\) self-adjoint. Comparing
\(\ip{\pile M(z^{-k}f)}{z^{-l}g}\) with \(\ip{z^{-k}f}{Y(z^{-l}g)}\) shows that the adjoint
of \(\pile M\) is multiplication by
\[
 Y=\sum_{n\geq0}z^{-n}\,G^{-1}\vartheta(M_n)^{\mathsf T}G,
\]
with no projection. For \(M=\RaiseFull{m}\), the transpose
\(\QQ_m^\vartheta(z^{-1})^{\mathsf T}\) is the identity with the column
\((\Pser{0}^\vartheta,\dots,\Pser{m-1}^\vartheta)^{\mathsf T}\) added to its first column.
By Lemma~\ref{lem:moments}, the first row and column of \(G\) are those of the identity.
Conjugating by \(G\) therefore multiplies this column on the
left by \(G^{-1}\) and leaves the other columns alone.

The inverse Gram matrices have entries \(\Gram{\bl}^{-1}_{00}=1\),
\(\Gram{\bl}^{-1}_{0j}=\Gram{\bl}^{-1}_{j0}=0\) for \(j\geq1\), and
\[
 \Gram{\bl}^{-1}_{ij}=(-1)^{i+j+m}\betaz^{-1}\Newton{m-i-j}(\bl)
 \qquad\text{for }1\leq i,j<m.
\]
Indeed, \eqref{eq:momentseries} and \eqref{eq:Yexpansion} give
\[
 \left(\sum_{s\geq1}\mom{s}{\bl}\,u^s\right)
 \left(\sum_{a\geq0}(-1)^a\Newton{a}(\bl)\,u^a\right)=\betaz u^m.
\]
For \(1\leq i,k<m\), comparing the coefficients of \(u^{m+i-k}\) gives
\[
 \sum_{j=1}^{m-1}\mom{i+j}{\bl}\,
 (-1)^{j+k+m}\betaz^{-1}\Newton{m-j-k}(\bl)=\delta_{ik},
\]
using the vanishing initial moments and the vanishing of Newton polynomials with negative indices.
Substituting these inverse entries gives the stated adjoint entries.
\end{proof}

\begin{proof}[Proof of Theorem~\ref{thm:factorization}]
By Theorem~\ref{thm:adjointexplicit}, \eqref{eq:companionentries} and \eqref{eq:Rk},
\(\Raise{m}^\dagger\Comp{m}=\Akernel{m}\). In each column \(j<m-1\), both matrices
have a single \(1\) in row \(j+1\). For the last column,
\((\Raise{m}^\dagger\Comp{m})_{0,m-1}=\betaz\Qan\) and, for \(1\leq i<m\),
\[
 (\Raise{m}^\dagger\Comp{m})_{i,m-1}=\betaz(\Raise{m}^\dagger)_{i0}+(\Comp{m})_{i,m-1}
 =\Rlast{m-i}(z/q_3^{i-1}).
\]
So Proposition~\ref{prop:currentkernel} is the formula of the theorem. By Theorem~\ref{thm:fock}(ii), the
pairing is the orthogonal sum of its restrictions to the subspaces \(z^{-k}\Csec{m}\), and on
\(z^{-k}\Csec{m}\) it is \(\Mq^k\) times its restriction to the constant sector. So \(\DD\) is self-adjoint,
and \(\piz\) is the orthogonal projection onto the constant sector. On the constant sector \(\Raise{m}\) and
\(\DD\) are the identity, and the constant term of \(\Akernel{m}\) in \(z^{-1}\) is \(\Comp{m}\), so
\(\piz\Lax{m}=\Comp{m}\) there. Since \(\Lax{m}\) is self-adjoint, being diagonal in an orthogonal basis, so is
\(\Comp{m}\). The operators \(\Comp{m}\) and \(\DD\) commute. For \(m=1\) and \(v_1=1\), \(\Comp{1}=1\), \(\Raise{1}=\pile\Qcr\) and
\(\Raise{1}^\dagger=\Qan\).
\end{proof}

\section{Module structure and spectral formulas}\label{sec:consequences}

The conjugacy gives the module presentations and the corner eigenvectors.
We then prove the transfer determinant directly in the fibre model and work out
the resulting operators and eigenvectors in rank two.

\subsection{The module presentation}\label{sec:structure}

The module \(\Gmod{m}\) of \eqref{eq:modulequotient} is the quotient by a monic
operator polynomial of degree \(m\), with coefficients acting on the auxiliary
polynomial space \(V\) of Section~\ref{sec:mainresults}.
For \(\xi\in V\), define polynomials in \(V[\Llet]\) recursively by
\[
 \Trunc{1}(\xi):=\xi,\qquad
 \Trunc{j+1}(\xi):=\Llet\,\Trunc{j}(\xi)-\Xhat{j}\xi,\quad1\leq j\leq m.
\]
Equivalently,
\[
 \Trunc{j}(\xi)=\Llet^{j-1}\xi-\sum_{k=1}^{j-1}\Llet^{j-1-k}\Xhat{k}\xi.
\]
These polynomials define the coordinate map
\[
 V^m\longrightarrow\Gmod{m},\qquad
 (\xi_1,\dots,\xi_m)\longmapsto\sum_{j=1}^m\Trunc{j}(\xi_j).
\]

\begin{corollary}\label{cor:structure}
\begin{enumerate}
\item[(i)] The map
\[
 V[\Llet]\to\Wsp{m},\qquad
 \Llet^iz^{-k}\kket{\bl}\longmapsto\taut^{-i}\Shift^k\taut^{-1}\pulm\fix{\bl},
\]
is \(\K[\Llet]\)-linear and surjective, with kernel
\(\Llet^m\Miuraproj(\Llet^{-1})\,V[\Llet]\).
The induced map \(\OmapG\) is an isomorphism from \(\Gmod{m}\) to
\(\Wsp{m}\). The map \(\OmapG\) identifies \(\GLax{m}\) with multiplication by \(\Llet\).
\item[(ii)] In companion coordinates, \(\Llet\) acts by \(\Miura{m}\) and
\(\Shift\) by \(z^{-1}\Delta\). Here \(\Delta\) is the diagonal matrix with
entries \(1,q_3^{-1},\dots,q_3^{1-m}\).
The maps
\[
 \bigl(V^m;\,\Miura{m},\,z^{-1}\Delta\bigr)\xrightarrow{\;\Vconj{m}\;}
 \bigl(\Hsp{m};\,\Lax{m},\,z^{-1}\bigr)\xrightarrow{\;\Omap\;}
 \bigl(\Wsp{m};\,\dtaut,\,\Shift\bigr)
\]
are isomorphisms of \(\K\langle\Llet,\Shift\rangle\)-modules, each triple
being a space with the actions of \(\Llet\) and \(\Shift\).
Their composite is \(\OmapG\) in companion coordinates.
\end{enumerate}
\end{corollary}

\begin{proof}
Write \(\OmapG\) also for the map \(V[\Llet]\to\Wsp{m}\) of part~(i), so that
\(\OmapG(\Llet^i\xi)=(\dtaut)^iU_1\xi\). Then \(\OmapG\,\Trunc{1}=U_1=\Theta_1\).
The recursion and Proposition~\ref{prop:Theta} give
\(\OmapG\,\Trunc{j}=\Theta_j\) for \(1\leq j\leq m+1\).
By \eqref{eq:OmegaV}, the restriction of \(\OmapG\) to the image of \(V^m\) in
\(V[\Llet]\) is \(\Omap\Vconj{m}\).

By \eqref{eq:Thetam}, the map \(\OmapG\) vanishes on
\(\Trunc{m+1}(V)=\Llet^m\Miuraproj(\Llet^{-1})V\).
The map intertwines multiplication by \(\Llet\) with \(\dtaut\), so
\(\OmapG\) vanishes on \(\Llet^m\Miuraproj(\Llet^{-1})V[\Llet]\).
Monic division gives the direct sum of this subspace and the image of \(V^m\).
Since \(\Omap\Vconj{m}\) is an isomorphism, the kernel and surjectivity assertions
of part~(i) follow, and \(\Omap^{-1}\OmapG=\Vconj{m}\) in companion coordinates.
Multiplication by \(\Llet\) is represented by \(\Miura{m}\), since
\(\Llet\,\Trunc{j}(\xi)=\Trunc{j+1}(\xi)+\Trunc{1}(\Xhat{j}\xi)\).

For the shift, each term of \eqref{eq:Vm} applies to \(\xi_l\) an operator \(A\,\DD^{l-1}\), where \(A\) is
\(\betaz\Qan\), \(1\) or \(-\Rlast{l-j}(z/q_3^{j-1})\). Each such \(A\) is multiplication by a series in
\(z^{-1}\), so \(A\) commutes with \(z^{-1}\). Since
\(z^{-1}\DD^{l-1}=q_3^{1-l}\DD^{l-1}z^{-1}\), we obtain
\(z^{-1}\Vconj{m}=\Vconj{m}\,z^{-1}\Delta\).
The map \(\Omap\) intertwines multiplication by \(z^{-1}\) with \(\Shift\),
and \(\Lax{m}\) with \(\dtaut\).
Thus \(\Vconj{m}\) and \(\Omap\) are isomorphisms of \(\K\langle\Llet,\Shift\rangle\)-modules, and
\(\Shift\) acts by \(z^{-1}\Delta\) in the coordinates from \(V^m\).
\end{proof}

On \(\Wfib{\bl}\), multiplication by \(\dtaut\) has
\(m+|\Rem{\bl}|\) distinct eigenvalues. If \(\Rem{\bl}\neq\varnothing\), the operator
therefore satisfies no monic relation of degree \(m\) with scalar coefficients
on \(\Wfib{\bl}\). The degree-\(m\) relation in \eqref{eq:modulequotient} has the
endomorphisms \(\Xhat{k}\) of \(V\) as coefficients. For \(m=1\), this
relation is linear in \(\Llet\), and \(\Gmod{1}=V\).

For the full presentation, take the free module
\(\mathcal F:=\K\langle\Llet,\Shift\rangle\otimes\La^{\otimes m}\).
Embed \(V\) into \(\mathcal F\) by \(z^{-k}f\mapsto\Shift^k\otimes f\).
The formulas for \(\Trunc{j}(\xi)\) also define elements of \(\mathcal F\), with each
\(\Xhat{k}\) acting on \(\xi\in V\) before this embedding.

\begin{samepage}
\begin{corollary}[Presentation of \(\Wsp{m}\)]\label{cor:presentation}
The module \(\Wsp{m}\) is the quotient of \(\mathcal F\) by the submodule generated by
\[
 \Llet^m\xi-\sum_{k=1}^m\Llet^{m-k}\Xhat{k}\xi,\qquad
 \Shift\,\Trunc{j}(\xi)-q_3^{1-j}\Trunc{j}(z^{-1}\xi),\quad2\leq j\leq m,
\]
for \(\xi\in V\). The quotient map is determined by
\(1\otimes\kket{\bl}\mapsto\taut^{-1}\pulm\fix{\bl}\).
\end{corollary}
\end{samepage}

The shift relations encode the action in Corollary~\ref{cor:structure}(ii).

\begin{proof}
Let \(\Phi\colon\mathcal F\to\Wsp{m}\) be the map specified on the generators, and
let \(\mathcal R\subset\mathcal F\) be the submodule generated by the given elements. On \(V\), \(\Phi\) is
\(\OmapG\), so in the coordinates from \(V^m\), \(\Phi\,\Trunc{j}(\xi)\) is \(\xi\) in place \(j\) for
\(j\leq m\), and \(\Phi\,\Trunc{m+1}(\xi)=0\). By Corollary~\ref{cor:structure}, \(\mathcal R\subset\ker\Phi\).
By the recursion and the displayed shift relations, the image of
\(\Trunc{1}(V)+\dots+\Trunc{m}(V)\) in \(\mathcal F/\mathcal R\) is stable under \(\Llet\) and \(\Shift\).
This subspace contains the image of the generating space \(\La^{\otimes m}\),
so the subspace equals \(\mathcal F/\mathcal R\). Thus the map
\(V^m\to\mathcal F/\mathcal R\), \((\xi_j)\mapsto\sum_j\Trunc{j}(\xi_j)\), is onto, and its composite with \(\Phi\)
is the isomorphism \(\OmapG\) in the coordinates from \(V^m\). Hence \(\Phi\) induces an isomorphism
\(\mathcal F/\mathcal R\to\Wsp{m}\).
\end{proof}

\subsection{Eigenvectors and the kernel}\label{sec:eigen}

The projected q-Miura operator \(\Miuraproj(y)\) acts on the auxiliary
polynomial space \(V\), with the total grading of
Definition~\ref{def:auxiliary}. For \(v\in V\), the companion identity is
\begin{equation}\label{eq:block}
 (1-y\Miura{m})(v,yv,\dots,y^{m-1}v)=(\Miuraproj(y)v,0,\dots,0).
\end{equation}
For \(c\neq0\), the eigenvectors of \(\Miura{m}\) with eigenvalue \(c\) are
therefore the vectors \((c^{m-1}\varphi,\dots,c\varphi,\varphi)\) with
\(\varphi\neq0\) and \(\Miuraproj(1/c)\varphi=0\).
For an \(m\)-tuple \(\bl\) and \(a\in\Add{\bl}\), define
\(\Meig{\bl}{a}\in V\) by the recursion
\begin{equation}\label{eq:phicorners}
 \Meig{\bl}{a}:=\betaz^{-1}\kket{\bl}
 +q_3^{1-m}z^{-1}\sum_{r\in\Rem{\bl}}\cores{\bl}{r}\,
 \frac{\ceig{\bplus{r}}}{\ceig{\bplus{r}}-\ceig{a}}\,\Meig{\bl-r}{r}.
\end{equation}
The empty multipartition gives the constant polynomial \(\betaz^{-1}\).
For general \(\bl\), only the states \(\kket{\bmu}\) with \(\bmu\subseteq\bl\) occur.

\begin{theorem}[Eigenvectors and the kernel]\label{thm:kernel}
Let \(\bl\) be an \(m\)-tuple and \(a\in\Add{\bl}\).
\begin{enumerate}
\item[(i)] The vector \((\ceig{a}^{m-1}\Meig{\bl}{a},\dots,\ceig{a}\Meig{\bl}{a},\Meig{\bl}{a})\) is an eigenvector
of \(\Miura{m}\) with eigenvalue \(\ceig{a}\), and \(\Vconj{m}\) sends it to the eigenvector
\(\Omap^{-1}\pt{\bl}{a}\) of \(\Lax{m}\). For \(|\bl|=N\) these vectors form an eigenbasis of \(\Miura{m}\) on total
degree \(N\), with determinant
\[
 \det(1-y\Miura{m})=\det\Miuraproj(y)=\prod_{|\bmu|=N}\ \prod_{b\in\Add{\bmu}}(1-y\,\ceig{b}).
\]
\item[(ii)] Let \(c\neq0\) and \(N\geq0\). On total degree \(N\), the kernel of \(\Miuraproj(1/c)\) has the basis
\(\Meig{\bmu}{b}\), one for each \(m\)-tuple \(\bmu\) of size \(N\) with an addable corner \(b\) of content
\(\ceig{b}=c\). The operator \(\Miuraproj(1/c)\) is invertible there unless \(c\) is such a content.
For each such pair \((\bmu,b)\), \(\Meig{\bmu}{b}\) is the unique kernel vector
with constant term \(\betaz^{-1}\kket{\bmu}\) in \(z^{-1}\).
\item[(iii)] Without the projection,
\[
 \Miuraop(1/\ceig{a})\,\Meig{\bl}{a}
 =-\ceig{a}^{-m}\sum_{n\geq1}(z/q_3)^n\,\Hiso\pushp\Shift^{n-1}\pt{\bl}{a},
\]
and the coefficient of \(z\) here is \(-\ceig{a}^{-m}q_3^{-1}\kket{\bl+a}\).
\item[(iv)] As series in \(y\),
\[
 \Miuraproj(y)^{-1}\kket{\bl}
 =\sum_{\tilde a\in\Add{\bl}}\frac{\res{\bl}{\tilde a}\,\ceig{\tilde a}^m}{1-y\,\ceig{\tilde a}}\,
 \Meig{\bl}{\tilde a}.
\]
\end{enumerate}
\end{theorem}

Part~(iv) expresses the first diagonal block of
\((1-y\Miura{m})^{-1}\), applied to \(\kket{\bl}\), as a finite sum of
corner vectors with simple poles at \(y=\ceig{a}^{-1}\).
The kernel in part~(ii) can have dimension greater than one.
In rank one with \(v_1=1\), the partitions \((3,1)\) and \((2,2)\)
both have the addable corner \(\bx{2}{0}{}\), of content \(\Mt^2\).

\begin{proof}[Proof of Theorem~\ref{thm:kernel}]
(i) Let \(\nu\colon\Wsp{m}\to V\) send a class \(w\) to the last coordinate of
\((\Omap\Vconj{m})^{-1}w\).
By \eqref{eq:Vm}, the \(\auxe{0}\)-component of \(\Vconj{m}\xi\) is \(\betaz\Qan\DD^{m-1}\xi_m\), and
\(\betaz\Qan\DD^{m-1}\) is invertible. So the last coordinate of \(\Vconj{m}^{-1}(\xi\otimes\auxe{j})\) vanishes
for \(1\leq j<m\). Since \(U_j\xi=\Omap(\xi\otimes\auxe{j})\), we have \(\nu\,U_j=0\) for \(1\leq j<m\).
By \eqref{eq:OmegaV} for \(l=m\), \(\nu\,\Theta_m=1\). The definition of \(\Theta_m\) therefore gives
\(\nu\,U_m\DD^{m-1}=1\), hence \(\nu\,U_m=\DD^{1-m}\).
Corollary~\ref{cor:structure}(ii) gives \(\nu\,\Shift=q_3^{1-m}z^{-1}\nu\).

Apply \(\nu\) to the second formula of Lemma~\ref{lem:expansion}, with \(\bmu=\bl\).
Since \(\taut^{-j}\pulm\fix{\bl}=U_j\kket{\bl}\) and \(\nu\,U_j=0\) for \(1\leq j<m\),
the polynomial sum contributes only the term \(j=m\), in which
\(h_m(\bl,\ceig{a})=\betaz^{-1}\) and \(\DD^{1-m}\kket{\bl}=\kket{\bl}\). So the vectors \(\nu\,\pt{\bl}{a}\)
satisfy \eqref{eq:phicorners}, and \(\nu\,\pt{\bl}{a}=\Meig{\bl}{a}\) by induction on \(|\bl|\). By
Proposition~\ref{prop:fixed} and Theorem~\ref{thm:companion}, \(\xi=(\Omap\Vconj{m})^{-1}\pt{\bl}{a}\) is an
eigenvector of \(\Miura{m}\) with eigenvalue \(\ceig{a}\). The rows \(2,\dots,m\) of
\(\Miura{m}\xi=\ceig{a}\xi\) say that \(\xi_j=\ceig{a}^{m-j}\xi_m\), and \(\xi_m=\Meig{\bl}{a}\).
The fixed-point classes form a basis of \(\Wdeg{m}{N}\), which gives the eigenbasis
and the product formula for \(\det(1-y\Miura{m})\).
Multiplying \(1-y\Miura{m}\) on the right by the identity matrix with first column \((1,y,\dots,y^{m-1})\), of
determinant \(1\), gives by \eqref{eq:block} a matrix with first column \((\Miuraproj(y),0,\dots,0)\) and a
unitriangular lower right block of size \(m-1\). So \(\det(1-y\Miura{m})=\det\Miuraproj(y)\).

(ii) By \eqref{eq:block}, \(\varphi\mapsto(c^{m-1}\varphi,\dots,c\varphi,\varphi)\) is an isomorphism from the
kernel of \(\Miuraproj(1/c)\) onto the eigenspace of \(\Miura{m}\) with eigenvalue \(c\), and (i) gives a basis of
this eigenspace. Distinct boxes have distinct contents, so each \(\bmu\) has at most one addable corner of
content \(c\). By \eqref{eq:phicorners} the constant term of \(\Meig{\bmu}{b}\) in \(z^{-1}\) is
\(\betaz^{-1}\kket{\bmu}\), and the \(\kket{\bmu}\) are linearly independent.

(iii) Let \(\pi_{>0}=1-\pi\), and let
\(\mathbf P=\sum_{n\geq1}z^n\,\Hiso\pushp\Shift^{n-1}\), a map from \(\Wsp{m}\) to series in \(z\). By the
definitions of \(U_j\) and \(P_{n,j}\), \(\mathbf P\,U_j=\pi_{>0}\Pser{j}(z)\) on \(V\). So, by the definition of
\(\Theta_j\) and \(\DD C(z)=C(z/q_3)\DD\), for \(1\leq j\leq m\)
\[
 \DD\,\mathbf P\,\Theta_j
 =\pi_{>0}\Bigl(\Pser{j}(z/q_3)-\sum_{k=1}^{j-1}\Pser{j-k}(z/q_3)\,\Rlast{k}(z/q_3^{j-k})\Bigr)\DD^j
 =-\pi_{>0}\,\gamma^j\,\DD^j\,\Xk{j}(q_2\gamma^mz),
\]
the last step by \eqref{eq:WH}, since \(\Rlast{j}(z)\) has no positive powers of \(z\). By (i) and
\eqref{eq:OmegaV}, \(\pt{\bl}{a}=\sum_l\ceig{a}^{m-l}\,\Theta_l\,\Meig{\bl}{a}\), so
\[
 \DD\,\mathbf P\,\pt{\bl}{a}
 =-\ceig{a}^m\,\pi_{>0}\sum_{l=1}^m\ceig{a}^{-l}\gamma^l\,\DD^l\,\Xk{l}(q_2\gamma^mz)\,\Meig{\bl}{a}
 =-\ceig{a}^m\,\pi_{>0}\,\Miuraop(1/\ceig{a})\,\Meig{\bl}{a}.
\]
Since \(\pi\,\Miuraop(1/\ceig{a})\,\Meig{\bl}{a}=\Miuraproj(1/\ceig{a})\,\Meig{\bl}{a}=0\), this is the formula of
(iii), and \(\pushp\pt{\bl}{a}=\fix{\bl+a}\) gives the coefficient of \(z\).

(iv) By \eqref{eq:OmegaV} for \(l=1\) and Lemma~\ref{lem:fibre}(i),
\(\Omap\Vconj{m}(\kket{\bl},0,\dots,0)=U_1\kket{\bl}=\sum_{\tilde a}\res{\bl}{\tilde a}\ceig{\tilde a}\pt{\bl}{\tilde a}\).
With Theorem~\ref{thm:companion} and (i) this gives
\[
 (1-y\Miura{m})^{-1}(\kket{\bl},0,\dots,0)
 =\sum_{\tilde a\in\Add{\bl}}\frac{\res{\bl}{\tilde a}\,\ceig{\tilde a}}{1-y\,\ceig{\tilde a}}\,
 \bigl(\ceig{\tilde a}^{m-1}\Meig{\bl}{\tilde a},\dots,\ceig{\tilde a}\Meig{\bl}{\tilde a},\Meig{\bl}{\tilde a}\bigr).
\]
By \eqref{eq:block} with \(v=\Miuraproj(y)^{-1}\kket{\bl}\), the first coordinate of the left side is
\(\Miuraproj(y)^{-1}\kket{\bl}\), and the first coordinate of the right side is the right side of (iv).
\end{proof}

\paragraph{Corner operators.}
The three-factor decomposition of \(\Miuraop(y)\) in the proof of
Proposition~\ref{prop:WH} ends with the lowering factor
\[
 T_-(y):=1-\sum_{k\geq1}y^k\DD^k\,
 \Pser{k}^\vartheta\bigl((q_3^{m-1}z)^{-1}\bigr).
\]
Every nonconstant term lowers the Fock degree and increases the degree in \(z^{-1}\)
by the same amount. On total degree \(N\), the operator \(T_-(y)-1\) is nilpotent
of exponent at most \(N+1\). The inverse of \(T_-(y)\) is the finite geometric sum
\(\sum_{j=0}^N(1-T_-(y))^j\), whose coefficients are formal series in \(y\).
Regard \(\La^{\otimes m}\) as the constant polynomials in \(V\).
For an independent spectral parameter \(c\), define
\begin{equation}\label{eq:corneroperator}
 X(c):=\left.\betaz^{-1}T_-(c^{-1})^{-1}\right|_{\La^{\otimes m}}.
\end{equation}
Let \(\CY(w)\) be the diagonal operator determined by
\(\CY(w)\kket{\bl}:=\Yfn{\bl}(w)\kket{\bl}\).

\begin{corollary}[Corner operators]\label{cor:corneroperator}
On each homogeneous component, \(X(c)\) is a rational family of maps
\(\La^{\otimes m}\to V\otimes_\K\K(c)\), with constant term
\(\betaz^{-1}\) times the input, satisfying
\begin{equation}\label{eq:corneroperatoridentity}
 \Miuraproj(c^{-1})\,X(c)=-c^{-m}\CY(c^{-1}).
\end{equation}
For a fixed coloured box \(a\), the family \(X(c)\) is regular at \(c=\ceig{a}\)
on the span of the \(\kket{\bl}\) for which \(a\) is addable.
The specialization \(X_a:=X(\ceig{a})\) on this subspace satisfies
\[
 X_a\kket{\bl}=\Meig{\bl}{a},
 \qquad
 \Miuraproj(\ceig{a}^{-1})X_a=0.
\]
\end{corollary}

\begin{proof}
Multiplying \(\Miuraop(y)\) by \(T_-(y)^{-1}\) cancels the lowering factor.
On constant input, the middle factor is
\(\sum_{j\geq0}(-y)^j\Newton{j}=-\betaz y^m\CY(y)\) by \eqref{eq:Yexpansion}.
Projection removes all the nonconstant terms of the remaining raising factor,
giving \eqref{eq:corneroperatoridentity}. Since \(\Miuraproj(0)=1\),
\[
 X(c)=-c^{-m}\Miuraproj(c^{-1})^{-1}\CY(c^{-1}),
\]
which proves rationality on each homogeneous component. The lowering factor has
constant term \(1\) in \(z^{-1}\), giving the stated constant term of \(X(c)\).

Since \(a\in\Add{\bl}\), the function \(\Yfn{\bl}(1/c)\) is regular and vanishes
at \(c=\ceig{a}\). If \(X(c)\kket{\bl}\) had a pole there, equation~\eqref{eq:corneroperatoridentity}
would place the leading pole coefficient of \(X(c)\kket{\bl}\) in
\(\ker\Miuraproj(\ceig{a}^{-1})\).
That coefficient would have zero constant term in \(z^{-1}\), since the constant
term of \(X(c)\kket{\bl}\) is independent of \(c\).
The constant-term map is injective on this kernel by Theorem~\ref{thm:kernel}(ii),
so such a pole is impossible. Specializing \eqref{eq:corneroperatoridentity}
therefore gives a vector in the projected q-Miura kernel. The specialized vector has constant term
\(\betaz^{-1}\kket{\bl}\), which identifies \(X(\ceig{a})\kket{\bl}\) with \(\Meig{\bl}{a}\) by the
uniqueness statement in Theorem~\ref{thm:kernel}(ii).
\end{proof}

In the unprojected product, the raising factor has a pole at \(c=\ceig{a}\)
on \(\kket{\bl}\), by \eqref{eq:Pchains}. This pole cancels the diagonal zero,
and Theorem~\ref{thm:kernel}(iii) describes the resulting positive powers of \(z\).

\subsection{The transfer determinant}\label{sec:transfer}

The fixed-point fibre model expresses the transfer matrix
\(\TransferMat{\bl}{u}\) in the interpolation basis. The resulting column
relations determine \(\det(I_m-\TransferMat{\bl}{u})\).

\begin{proof}[Proof of Theorem~\ref{thm:transfer}]
By Lemma~\ref{lem:fibre}, \(\Omap\) identifies the span of the \(\kket{\bl}\otimes\auxe{j}\) and the
\(z^{-1}\Omap^{-1}\pt{\bl-r}{r}\) with the functions on \(\{\ceig{a}:a\in\Add{\bl}\}\): \(\kket{\bl}\otimes\auxe{j}\) is the
function \(x^j\), \(z^{-1}\Omap^{-1}\pt{\bl-r}{r}\) the function \(x/(x-\ceig{\bplus{r}})\), and \(\Lax{m}\)
multiplication by \(x\). So \(\Transfer{m}{u}\) preserves \(\kket{\bl}\otimes\K^m\).
The entries of \(W=\TransferMat{\bl}{u}\auxe{0}\) are the coefficients of the polynomial
\(P\) of degree less than \(m\) in the expansion of \(u/(u-x)\) in this basis.
Since \(\Yfn{\bl}(1/x)\) vanishes at the \(\ceig{a}\), there \(u/(u-x)\) equals
\[
 g(x)=\frac{u}{u-x}\Bigl(1-\frac{\Yfn{\bl}(1/x)}{\Yfn{\bl}(1/u)}\Bigr).
\]
The function \(g\) has no pole at \(x=u\), has simple poles at the \(\ceig{\bplus{r}}\),
and has polynomial part of degree at most \(m-1\).
Thus \(g=P+\sum_rt_r\,x/(x-\ceig{\bplus{r}})\), and the \(00\) entry is
\(W_0=P(0)=g(0)=1-\Yfn{\bl}(1/u)^{-1}\).

By definition of \(\Omap\), \(\kket{\bl}\otimes\auxe{j}=\Lax{m}^j(\kket{\bl}\otimes\auxe{0})\) for \(j<m\).
With \(y=1/u\), column \(j\) of \(\TransferMat{\bl}{u}\) is
\(y^{-j}\bigl(W-\sum_{n<j}y^n\auxe{n}\bigr)\).

Consequently,
\(I_m-\TransferMat{\bl}{u}=A-W\,r^{\mathsf T}\), where \(A_{ij}=y^{i-j}\) for \(i\leq j\), \(A_{ij}=0\) for \(i>j\), and
\(r=(1,y^{-1},\dots,y^{1-m})^{\mathsf T}\). The rank-one determinant identity,
with \(\det A=1\) and \(r^{\mathsf T}A^{-1}=(1,0,\dots,0)\), gives
\(\det(I_m-\TransferMat{\bl}{u})=1-W_0=\Yfn{\bl}(1/u)^{-1}\).
\end{proof}

\subsection{An example in rank two}\label{sec:example}

Let \(m=2\), with arbitrary generic framing weights \(v_1,v_2\).
Both Lax forms act on two copies of the auxiliary polynomial space
\(V=\La^{\otimes2}[z^{-1}]\). For the NS form, these copies are the
\(\auxe{0}\)- and \(\auxe{1}\)-components of \(\Hsp{2}\).
The power sums \(\pk{k}{1},\pk{k}{2}\) and derivatives
\(\partial_{\alpha,k}=\partial/\partial\pk{k}{\alpha}\) give explicit
formulas for every entry. We use \(\betaz=-(v_1v_2)^{-1}\),
\(\DD=\Twist{q_3}\), and \(\pi=\pile\). Thus
\(\DD f(z)=f(z/q_3)\), where \(q_3=\Mt/\Mq\).
On \(V^2\), the dilation acts componentwise.
The projection \(\pi\) retains the nonpositive powers of \(z\) after
the operators to its right have acted on the input.

\paragraph{Power-sum operators.}
The two half-vertex series of Theorem~\ref{thm:raising} and
Corollary~\ref{cor:Qdagger} specialize to
\begin{equation}\label{eq:ranktwohalves}
 \begin{aligned}
 \Qcr(z)&=\exp\Bigl(\sum_{k\geq1}\frac{1-\Mt^k}{k}\,
                   z^k(\pk{k}{1}+\pk{k}{2})\Bigr),\\
 \Qan&=\exp\Bigl(\sum_{k\geq1}(1-\Mq^k)\Mq^{-k}z^{-k}
                   (\partial_{1,k}+q_3^k\partial_{2,k})\Bigr).
 \end{aligned}
\end{equation}
The first series acts by multiplication. The second acts by the translations
\(\pk{k}{1}\mapsto\pk{k}{1}+(1-\Mq^k)\Mq^{-k}z^{-k}\) and
\(\pk{k}{2}\mapsto\pk{k}{2}+(1-\Mq^k)\Mq^{-k}q_3^kz^{-k}\).
In particular, \(\Qan\) preserves \(V\), and its inverse is obtained by
negating the exponent.

Substituting \(m=2\), \(q_1=\Mt^{-1}\), and \(q_2=\Mq\) into
\eqref{eq:Xcur} gives the two vertex operators
\begin{equation}\label{eq:ranktwovertices}
 \begin{aligned}
 \eta^-_1(w)
 &=\exp\Bigl(\sum_{k\geq1}\frac{1-\Mt^k}{k}
                    (w/\Mq)^k\pk{k}{1}\Bigr)\\
 &\qquad{}\times
   \exp\Bigl(\sum_{k\geq1}(1-\Mq^k)w^{-k}
                    \bigl(\partial_{1,k}+(q_3^k-1)\partial_{2,k}\bigr)\Bigr),\\
 \eta^-_2(w)
 &=\exp\Bigl(\sum_{k\geq1}\frac{1-\Mt^k}{k}
                    (w/\Mt)^k\pk{k}{2}\Bigr)
   \exp\Bigl(\sum_{k\geq1}(1-\Mq^k)(q_3/w)^k\partial_{2,k}\Bigr).
 \end{aligned}
\end{equation}
In each product the differential exponential acts first.
Since \(u_\alpha=-\gamma v_\alpha\), the current is
\(-\gamma\Xcur(w)=v_1^{-1}\eta^-_1(w)+v_2^{-1}\eta^-_2(w)\).
The raising series, lowering series, and first Cartan are therefore
\begin{equation}\label{eq:ranktwoparts}
 \begin{aligned}
 \Pser{1}(z)
   &=\bigl(v_1^{-1}\eta^-_1(\Mt z)+v_2^{-1}\eta^-_2(\Mt z)\bigr)_{>0},\\
 \Pser{1}^\vartheta(z^{-1})
   &=\bigl(v_1^{-1}\eta^-_1(\Mq z)+v_2^{-1}\eta^-_2(\Mq z)\bigr)_{<0},\\
 \Cartan{1}
   &=[w^0]\bigl(v_1^{-1}\eta^-_1(w)+v_2^{-1}\eta^-_2(w)\bigr).
 \end{aligned}
\end{equation}
Here \(>0\) and \(<0\) select powers of \(z\) in the operator series.
The arguments \(\Mt z\) and \(\Mq z\) are different because the lowering
series includes the degree rescaling in \(\vartheta\).

\paragraph{The NS form.}
The auxiliary Gram matrix of \eqref{eq:Gram} is independent of the Fock state:
\[
 \Gram{\bl}=\begin{pmatrix}1&0\\0&\betaz\end{pmatrix}.
\]
For \(f,g\in\La^{\otimes2}\), the full pairing is
\[
 \ip{z^{-k}f\otimes\auxe{i}}{z^{-l}g\otimes\auxe{j}}
 =\delta_{kl}\Mq^k\Gram{\bl}_{ij}\ip{f}{g},
\]
where the pairing on the right is the Fock pairing of
Lemma~\ref{lem:sourceadjoint}.
Since \(\Newton{1}=\Cartan{1}\), the companion polynomial is
\(x^2-\Cartan{1}(\bl)x-\betaz\). The three matrix factors are
\begin{equation}\label{eq:ranktwofactors}
 \Raise{2}=
 \begin{pmatrix}
  \pi\Qcr(z)&\pi\Pser{1}(z)\\0&1
 \end{pmatrix},
 \quad
 \Raise{2}^\dagger=
 \begin{pmatrix}
  \Qan&0\\ \betaz^{-1}\Pser{1}^\vartheta(z^{-1})&1
 \end{pmatrix},
 \quad
 \Comp{2}=
 \begin{pmatrix}0&\betaz\\1&\Cartan{1}\end{pmatrix}.
\end{equation}
The adjoint matrix follows from the displayed Gram matrix and
Theorem~\ref{thm:adjointexplicit}. All entries have the power-sum expressions
\eqref{eq:ranktwohalves}--\eqref{eq:ranktwoparts}.
The operator \(\Cartan{1}\) is self-adjoint, and
\(\Gram{\bl}\Comp{2}\) is symmetric, so \(\Comp{2}\) is self-adjoint.
The dilation \(\DD\) is self-adjoint and commutes with \(\Comp{2}\).
Multiplication of the four factors gives
\begin{equation}\label{eq:ranktwolax}
 \Lax{2}
 =\Raise{2}^\dagger\Comp{2}\DD\Raise{2}
 =\begin{pmatrix}
   0&\betaz\Qan\DD\\
   \DD\pi\Qcr(z)&
   \DD\pi\Pser{1}(z)+\bigl(\Cartan{1}+\Pser{1}^\vartheta(z^{-1})\bigr)\DD
  \end{pmatrix}.
\end{equation}
Thus \(\Lax{2}\) acts in every total degree by explicit power-sum
multiplication, differential operators, projection, and dilation.

\paragraph{The q-Miura form.}
There are only two nonzero fused currents:
\[
 \Xk{1}(w)=-\gamma^{-1}
       \bigl(v_1^{-1}\eta^-_1(w)+v_2^{-1}\eta^-_2(w)\bigr),
 \qquad
 \Xk{2}(w)=\frac{{:}\eta^-_1(w)\eta^-_2(q_3w){:}}{\gamma^2v_1v_2}.
\]
The q-Miura operator is the normal-ordered product
\begin{equation}\label{eq:ranktwomiuraproduct}
 \Miuraop(y)
 ={:}\bigl(1-yv_2^{-1}\DD\eta^-_2(\Mt z)\bigr)
      \bigl(1-yv_1^{-1}\DD\eta^-_1(\Mt z)\bigr){:}.
\end{equation}
The order of \(\DD\) and the vertex operators matters.
Moving the dilations to the left gives the coefficient of \(y^2\) as
\[
 \frac{\DD^2}{v_1v_2}\,{:}\eta^-_1(\Mt z)\eta^-_2(\Mt q_3z){:}.
\]
Adding the creation and derivative exponents in \eqref{eq:ranktwovertices} gives
\begin{equation}\label{eq:ranktwofusion}
 {:}\eta^-_1(\Mq z)\eta^-_2(\Mt z){:}=\Qcr(z)\Qan.
\end{equation}
Using \(\DD C(z)=C(z/q_3)\DD\) and projecting therefore gives
\begin{equation}\label{eq:ranktwomiuracoefficients}
 \begin{aligned}
 \Miuraproj(y)&=1-y\Xhat{1}-y^2\Xhat{2},\\
 \Xhat{1}
   &=\pi\bigl(v_1^{-1}\eta^-_1(\Mq z)+v_2^{-1}\eta^-_2(\Mq z)\bigr)\DD\\
   &=\DD\pi\Pser{1}(z)+
       \bigl(\Cartan{1}+\Pser{1}^\vartheta(z^{-1})\bigr)\DD,\\
 \Xhat{2}&=\betaz\DD\pi\Qcr(z)\Qan\DD.
 \end{aligned}
\end{equation}
For the second expression for \(\Xhat{1}\), split the vertex series into
positive, zero, and negative powers, and use \eqref{eq:ranktwoparts}.
In the last line, the projection follows the product \(\Qcr(z)\Qan\).
The companion matrix on \(V^2\) is
\[
 \Miura{2}=\begin{pmatrix}\Xhat{1}&\Xhat{2}\\1&0\end{pmatrix}.
\]
In source coordinates \((\xi_1,\xi_2)\) and target coordinates
\((\auxe{0},\auxe{1})\), the conjugating map and its inverse are
\begin{equation}\label{eq:ranktwoconjugator}
 \Vconj{2}=
 \begin{pmatrix}0&\betaz\Qan\DD\\1&0\end{pmatrix},
 \qquad
 \Vconj{2}^{-1}=
 \begin{pmatrix}0&1\\\betaz^{-1}\DD^{-1}(\Qan)^{-1}&0\end{pmatrix}.
\end{equation}
The explicit NS matrix and \eqref{eq:ranktwomiuracoefficients} give
\[
 \Lax{2}\Vconj{2}
 =\begin{pmatrix}\betaz\Qan\DD&0\\\Xhat{1}&\Xhat{2}\end{pmatrix}
 =\Vconj{2}\Miura{2}.
\]
This verifies the conjugacy in every degree. The series \(\Rlast{k}\) of
\eqref{eq:Vm} first appear in the conjugating map at rank three.

\paragraph{Total degree one.}
The degree-one Fock component has the two fixed-point vectors
indexed by \(((1),\varnothing)\) and \((\varnothing,(1))\), and
\(\Lax{2}\) has six eigenvalues on total degree one.
In degree one \cite[Section~3.1]{BCS}, \eqref{eq:kket} gives
\[
 \mathsf a=\kket{(1),\varnothing}=(1-\Mt)\bigl(1-q_3v_2/v_1\bigr)\,\pk{1}{1},\qquad
 \mathsf b=\kket{\varnothing,(1)}=(1-\Mt)\Bigl(\bigl(1-v_1/v_2\bigr)\,\pk{1}{2}+(1-q_3)\,\pk{1}{1}\Bigr).
\]
So the part of \(V\) of total degree one has the basis \(\mathsf a,\mathsf b,z^{-1}\). The addable corners
\(\bx{1}{0}{1},\bx{0}{1}{1},\bx{0}{0}{2}\) of \(((1),\varnothing)\) have the contents \(\Mt/v_1\), \(1/(\Mq v_1)\),
\(1/v_2\), and those of \((\varnothing,(1))\) are \(1/v_1\), \(\Mt/v_2\), \(1/(\Mq v_2)\). Put \(d=v_1-v_2\) and
\[
 A=\frac{(1-\Mq)(1-\Mt)(1-q_3v_2/v_1)}{\Mq v_1},\qquad
 B=\frac{(1-\Mq)(1-\Mt)(1-q_3v_1/v_2)}{\Mq v_2}.
\]
Write \(H_{\mathsf a}=\Cartan{1}((1),\varnothing)\) and \(H_{\mathsf b}=\Cartan{1}(\varnothing,(1))\), so that
\(H_{\mathsf a}=v_1^{-1}(\Mt+\Mq^{-1}-\Mt\Mq^{-1})+v_2^{-1}\). On the part of \(\Hsp{2}\) of total degree one,
\begin{align*}
 \Lax{2}(\mathsf a\otimes\auxe{0})&=\mathsf a\otimes\auxe{1},\\
 \Lax{2}(\mathsf a\otimes\auxe{1})&=-(v_1v_2)^{-1}\mathsf a\otimes\auxe{0}+H_{\mathsf a}\,\mathsf a\otimes\auxe{1}
 -(A/v_2)\,z^{-1}\otimes(\auxe{0}-v_2\auxe{1}),\\
 \Lax{2}(z^{-1}\otimes\auxe{0})&=q_3\,z^{-1}\otimes\auxe{1}+d^{-1}(v_1\mathsf a-v_2\mathsf b)\otimes\auxe{1},\\
 \Lax{2}(z^{-1}\otimes\auxe{1})&=q_3\,z^{-1}\otimes\bigl(-(v_1v_2)^{-1}\auxe{0}+(v_1^{-1}+v_2^{-1})\auxe{1}\bigr)
 +d^{-1}(\mathsf a-\mathsf b)\otimes\auxe{1},
\end{align*}
and the action on \(\mathsf b\otimes\auxe{j}\) is obtained by exchanging
\(\mathsf a,A,v_1\) with \(\mathsf b,B,v_2\).

In the basis \(\mathsf a,\mathsf b,z^{-1}\), with columns the images,
\[
 \Xhat{1}=\begin{pmatrix}H_{\mathsf a}&0&d^{-1}\\[2pt]0&H_{\mathsf b}&-d^{-1}\\[2pt]A&B&q_3(v_1^{-1}+v_2^{-1})\end{pmatrix},\quad
 \Xhat{2}=\begin{pmatrix}
 -\frac{1}{v_1v_2}-\frac{Av_1}{v_2d}&-\frac{B}{d}&-\frac{q_3}{v_2d}\\[3pt]
 \frac{A}{d}&-\frac{1}{v_1v_2}+\frac{Bv_2}{v_1d}&\frac{q_3}{v_1d}\\[3pt]
 -\frac{q_3A}{v_2}&-\frac{q_3B}{v_1}&-\frac{q_3^2}{v_1v_2}\end{pmatrix},
\]
so \(\Miura{2}\) is a \(6\times6\) matrix.
To evaluate \(\Vconj{2}\) from \eqref{eq:ranktwoconjugator}, use
\(\Qan\mathsf a=\mathsf a+v_1Az^{-1}\), \(\Qan\mathsf b=\mathsf b+v_2Bz^{-1}\),
\(\Qan z^{-1}=z^{-1}\), and \(\DD z^{-1}=q_3z^{-1}\).

\paragraph{Two corner eigenvectors.}
For a q-Miura kernel vector \(\varphi\) at a corner of content \(c\), the
companion eigenvector is \((c\varphi,\varphi)\). The conjugating map gives
the NS eigenvector
\[
 \Vconj{2}(c\varphi,\varphi)
 =\betaz\Qan\DD\varphi\otimes\auxe{0}
   +c\varphi\otimes\auxe{1}.
\]
For the corner \(\bx{1}{0}{1}\) of \(((1),\varnothing)\), the content is
\(\Mt/v_1\). The recursion \eqref{eq:phicorners} gives the kernel vector
\begin{align*}
 \Meig{(1),\varnothing}{\bx{1}{0}{1}}
 &=-v_1v_2\kket{(1),\varnothing}\\
 &\quad+\frac{v_1v_2}{\Mt}(1-\Mt)(1-q_3v_2/v_1)\,
          z^{-1}\kket{\varnothing,\varnothing}.
\end{align*}
Using the displayed actions of \(\Qan\) and \(\DD\), the corresponding
NS eigenvector with eigenvalue \(\Mt/v_1\) is
\begin{align*}
 \Omap^{-1}\pt{((1),\varnothing)}{\bx{1}{0}{1}}
 &=\kket{(1),\varnothing}\otimes(\auxe{0}-\Mt v_2\auxe{1})\\
 &\quad-(1-\Mt)(1-q_3v_2/v_1)\,z^{-1}\kket{\varnothing,\varnothing}
             \otimes(\auxe{0}-v_2\auxe{1}).
\end{align*}
For the corner \(\bx{1}{0}{2}\) of \((\varnothing,(1))\), the content is
\(\Mt/v_2\). The kernel vector is
\begin{align*}
 \Meig{\varnothing,(1)}{\bx{1}{0}{2}}
 &=-v_1v_2\kket{\varnothing,(1)}\\
 &\quad+\frac{v_1v_2}{\Mt}(1-\Mt)(1-q_3v_1/v_2)\,
          z^{-1}\kket{\varnothing,\varnothing}.
\end{align*}
The same conjugating map gives the NS eigenvector with eigenvalue \(\Mt/v_2\):
\begin{align*}
 \Omap^{-1}\pt{(\varnothing,(1))}{\bx{1}{0}{2}}
 &=\kket{\varnothing,(1)}\otimes(\auxe{0}-\Mt v_1\auxe{1})\\
 &\quad-(1-\Mt)(1-q_3v_1/v_2)\,z^{-1}\kket{\varnothing,\varnothing}
             \otimes(\auxe{0}-v_1\auxe{1}).
\end{align*}
Substitution into the degree-one matrices gives
\[
 \Miuraproj(v_1/\Mt)\Meig{(1),\varnothing}{\bx{1}{0}{1}}=0,
 \qquad
 \Miuraproj(v_2/\Mt)\Meig{\varnothing,(1)}{\bx{1}{0}{2}}=0.
\]
Together with \(\Lax{2}\Vconj{2}=\Vconj{2}\Miura{2}\), these identities
verify both NS eigenvalue equations.

Theorem~\ref{thm:kernel}(iii), applied to the empty multipartition, also
determines the unprojected terms. The coefficient of \(z\) in
\(\Miuraop(v_1)\cdot1\) is \(v_1\mathsf a/(q_3v_2)\), and in
\(\Miuraop(v_2)\cdot1\) it is \(v_2\mathsf b/(q_3v_1)\).

The relation of Corollary~\ref{cor:structure} reads
\(\dtaut\,\dtaut\,U_1\xi=\dtaut\,U_1\Xhat{1}\xi+U_1\Xhat{2}\xi\).
On \(\Wfib{((1),\varnothing)}\), multiplication by \(\dtaut\) has three distinct eigenvalues.
Thus \(\dtaut\) satisfies no monic quadratic relation with
scalar coefficients on this subspace.

\section{Comparison with existing models}\label{sec:comparison}

\subsection{Rank one}\label{sec:rankone}

We specialize to \(m=1\) and \(v_1=1\), and check that the operator, pairing, and
eigenvectors agree with the rank-one construction of \cite{M26}.

The two half-vertex series in \eqref{eq:rankoneNS} are
\[
 \Qcr(z)=\exp\Bigl(\sum_{n\geq1}\frac{1-\Mt^n}{n}\,p_nz^n\Bigr),\qquad
 \Qan=\exp\Bigl(-\sum_{n\geq1}\frac{1-\Mq^{-n}}{n}\,p_n^\perp z^{-n}\Bigr),
\]
where \(p_n^\perp:=n\,\partial/\partial p_n\).

By Theorem~\ref{thm:factorization}, \(\Lax{1}\) is the operator
\(\LL=\Qan\Twist{\Mt/\Mq}\pile\Qcr\) of \cite[Definition~3.1]{M26}, and Sections~\ref{sec:operator}--\ref{sec:consequences}
do not use \cite{M26}.

The pairing \(\ipinf{z^{-k}f}{z^{-l}g}=\Mq^k\delta_{kl}\ip{f}{g}_{\Mq,\Mt}\)
of \cite[(7)]{M26} agrees with the pairing of Theorem~\ref{thm:fock}(ii).
On \(\La\), Lemma~\ref{lem:sourceadjoint}, together with
\(\sigH(a_{-k})=-a_k\) and \cite[(2.33)]{BCS}, gives the adjoint
\(k(1-\Mq^k)(1-\Mt^k)^{-1}\partial/\partial p_k\) of \(p_k\),
as for \(\ip{\cdot}{\cdot}_{\Mq,\Mt}\).

The eigenbasis of
\cite[Proposition~3.11]{M26} is \(\PsiNS{\lambda}{s}=\Mq^{n(\lambda')}\Omap^{-1}\pt{\lambda}{s}\). Both sides are eigenvectors
with eigenvalue \(\ceig{s}\) and constant term \(J_\lambda\), since \(\pt{\lambda}{s}\) is the function
\(\ceig{s}\Yfn{\lambda}(1/x)/(\ceig{s}-x)\) on the contents, of value \(1\) at \(x=0\).
For a fixed degree and eigenvalue, an eigenvector of \(\LL\) is determined by
its constant term, by \cite[Proposition~3.11]{M26} and the independence of the \(J_\mu\).

The weight of the corner \(s\) in Kerov's transition measure \cite{Kerov},
\cite[Definition~3.8]{M26} is \(\tau_\lambda^s=\res{\lambda}{s}\).
The creation formula \cite[Proposition~4.8]{M26} agrees with
Lemma~\ref{lem:shiftfixed} under the rank-one identification.
The transfer formula \cite[Proposition~3.7]{M26} agrees with the
\(00\) entry of Theorem~\ref{thm:transfer}.

In total degree one, the companion operator has the action
\begin{align*}
 \Xhat{1}p_1
 &=(\Mt+\Mq^{-1}-\Mt\Mq^{-1})p_1+\Mt\Mq^{-2}(1-\Mq)z^{-1},\\
 \Xhat{1}z^{-1}
 &=(1-\Mt)p_1+\Mt\Mq^{-1}z^{-1}.
\end{align*}
The eigenvalues of \(\Xhat{1}\) are \(\Mt\) and \(\Mq^{-1}\), the contents of the
addable corners of \((1)\).

\subsection{The Gonz\'alez--Gorsky--Simental model}\label{sec:ggs}

Gonz\'alez--Gorsky--Simental \cite{GGS} realize their operators on the \(K\)-theory
of parabolic Gieseker spaces. We compare their correspondences with the fixed-origin
Hecke correspondence and identify the inverse Lax operator with their generator
\(\zone\).

\begin{proposition}\label{prop:GGS}
Let \(\Parab{m}{n}{n+1}\) be the parabolic Gieseker space of \cite{GGS}, and
\(\incl\colon\Hecke{n}\hookrightarrow\Parab{m}{n}{n+1}\) the locus where the
quotient sits at the origin. Identify the torus of \cite{GGS} with that of
\cite{BCS} by \((q,t,a_\alpha)=(q_2,q_1,v_\alpha)\). The localized
pushforward \(\incl_*\) is an isomorphism intertwining
\[
 \taut\leftrightarrow\zone,\qquad
 \pulm\leftrightarrow-\dplus,\qquad
 \pushp\leftrightarrow\dminus,\qquad
 \Shift\leftrightarrow-(\yone+\dplus\dminus).
\]
In particular \(\GLax{m}\) is \(\zone^{-1}\), the operator \(\yone\) of
\cite{CGM,GGS} is the untwisted pull-push \(\pippush\pimpull\), and the relation
\(\zone(q\dplus\dminus-\dminus\dplus)=qt\,(\dplus\dminus-\dminus\dplus)\zone\)
of \cite{GGS} on \(\Wsp{m}\) is the excess formula for the square
(2.37) of \cite{N2}.
\end{proposition}

\begin{proof}
The two spaces have the same fixed points. By \cite[(3.8) and Theorem~6.21]{GGS}, \(\zone\) acts on
\(\incl_*\pt{\bl}{a}\) by the weight \(\bwt{a}\), as \(\taut\) acts on \(\pt{\bl}{a}\). By \eqref{eq:localisation},
\(\dminus\) sends \(\incl_*\pt{\bl}{a}\) to \(\fix{\bl+a}\), as \(\pushp\) sends \(\pt{\bl}{a}\). By \cite[Lemma~5.8]{GGS} with \(k=0\)
and \eqref{eq:Wcharacter}, the coefficient of \(\incl_*\pt{\bl}{a}\) in \(\dplus\fix{\bl}\) is \(-\res{\bl}{a}\), so
\(\incl_*\pulm=-\dplus\) by \eqref{eq:pullfixed}. Put \(S'=-(q\dplus\dminus-\dminus\dplus)/(q-1)=-(\yone+\dplus\dminus)\)
on \(V_1\). Since \(\dminus\) forgets the largest sheaf, \(S'\) sends \(\incl_*\pt{\bmu}{r}\) into the span of the
\(\incl_*\pt{\bmu+r}{b}\). The relation of \cite{GGS} on \(V_1\) \cite[(2.8) and Lemma~3.22]{GGS}, with
\(qt=q_3^{-1}\), is equivalent to \(S'\zone-q_3\zone S'=-\dplus\dminus\zone\), which corresponds to
\eqref{eq:displacement}. If \(s_b\) is the coefficient of
\(\incl_*\pt{\bmu+r}{b}\) in \(S'\incl_*\pt{\bmu}{r}\), this relation gives
\[
 (\bwt{r}-q_3\bwt{b})s_b=\bwt{r}\res{\bmu+r}{b}.
\]
The factor on the left is nonzero, as in the proof of Lemma~\ref{lem:shiftfixed},
so these are exactly the coefficients of \(\Shift\).
Finally, Lemma~\ref{lem:excess} identifies \(\pippush\pimpull\) with \(\yone\).
\end{proof}

The auxiliary module itself also occurs in the poset construction of
Gonz\'alez--Gorsky--Simental \cite[Definition~3.20 and Theorem~6.21]{GGS}.
Let \(E\) be the poset of \(m\)-tuples of partitions, with edge weight \(\bwt{a}\)
for \(\bl\subset\bl+a\). The component \(V_1(E)\) has basis
\([\bl;\bwt{a}]\), indexed by these one-step chains, and
\(\zone[\bl;\bwt{a}]=\bwt{a}[\bl;\bwt{a}]\).
Under the fixed-point identification
\(\pt{\bl}{a}\leftrightarrow[\bl;\bwt{a}]\), this is the action of \(\taut\).
Thus \(\Wsp{m}\) is their component \(V_1(E)\), with the Lax operator acting by
\(\zone^{-1}\) and the shift by the combination in Proposition~\ref{prop:GGS}.

\subsection{The Carlsson--Gorsky--Mellit model}\label{sec:cgm}

The polynomial model of Carlsson--Gorsky--Mellit \cite{CGM} gives a second
rank-one comparison, by a plethystic change of variables. Here the parameters
are \((q,t)=(\Mq^{-1},\Mt)\), the inverses of the parameters used in
Section~\ref{sec:ggs}. In this convention the NS operator corresponds to
\(\zone\), whereas the geometric comparison identifies the NS operator with \(\zone^{-1}\).

We use the action of \(\mathbb B_{q,t}\) on
\(\bigoplus_kV_k\), \(V_k=\La[y_1,\dots,y_k]\), constructed in
\cite[Theorem~3.4.2(2)]{CGM}. The maps \(\dplus\colon
V_k\to V_{k+1}\) and \(\dminus\colon V_k\to V_{k-1}\) are the raising and
lowering operators, \(\yy{i}\) acts on \(V_k\) by multiplication, and the
\(\zz{i}\) act by the operators of \cite{CGM}. Let \(\LL\) be the rank-one operator of
Section~\ref{sec:rankone}.

In the polynomial action we use the Hecke operators
\[
 T_iF=\frac{(q-1)y_iF+(y_{i+1}-qy_i)s_iF}{y_{i+1}-y_i},
\]
where \(s_i\) interchanges \(y_i\) and \(y_{i+1}\). These operators satisfy
\((T_i-1)(T_i+q)=0\).

The standard polynomial \(\mathbb A_{q,t}\) action of \cite[Proposition~3.4.1]{CGM} is obtained
from this \(\mathbb B_{q,t}\) action through the homomorphism
\(\beta\colon\mathbb A_{q,t}\to\mathbb B_{q,t}\) of \cite[Section~3.2]{CGM}. In particular,
\[
 \beta(\zone^{\mathbb A})=-qt\,\yone\zone,\qquad
 \zone^{\mathbb A}(1)=-qt\,y,\qquad \zone(1)=1
\]
on \(V_1=\La[y]\). In this subsection \(\zone\) denotes the generator in the
\(\mathbb B_{q,t}\) action.

For this action, use the extension to \((y_1\cdots y_k)^{-1}V_k\) with
\(d_+^*=-(qty_1)^{-1}d_{+,CM}^*\), where \(d_{+,CM}^*\) is the operator of
\cite[Proposition~3.4.1]{CGM}. The commutator preserves \(V_1=\La[y]\) and satisfies
\([d_+^*,\dminus]=(q^{-1}-1)\zone\) there
\cite[proof of Proposition~3.2.6]{CGM}.
\begin{theorem}[The CGM intertwiner]\label{thm:rankone}
Let \(\Vone=\La[y]\) be the summand with \(k=1\) of the \(\mathbb B_{q,t}\) representation of
\cite[Theorem~3.4.2(2)]{CGM}, with parameters \((q,t)=(\Mq^{-1},\Mt)\). The ring isomorphism
\[
 \Plethmap\colon\La[z^{-1}]\to\Vone,\qquad
 \Plethmap(G)(X;y)=G\Bigl[\frac{X+(\Mq^{-1}-1)y}{1-\Mt^{-1}};\ z^{-1}=\Mt y\Bigr],
\]
satisfies \(\Plethmap\,\LL=\zone\,\Plethmap\). The map \(\Plethmap\) sends multiplication by
\(z^{-1}\) to \(\Mt\yone\), and \(\La\subset\La[z^{-1}]\) to \(\dplus(\La)\)
after \(X\mapsto X/(1-\Mt^{-1})\).
\end{theorem}

\begin{proof}
Write \(w=z^{-1}\) and \(X'=(X+(q-1)ty)/(1-t^{-1})\), and let \(\Exp[A]=\exp\sum_{k\geq1}p_k[A]/k\). In the
parameters \((q,t)\), \(\Qcr(z)=\Exp[(1-t)X/w]\), the operator \(\Qan\) is the translation \(X\mapsto X+(q-1)w\), and
\(\Twist{\Mt/\Mq}\) is \(w\mapsto qtw\). So \(\LL G(X;w)=H(X+(q-1)w;\,qtw)\), where \(H\) is the part of
\(\Exp[(1-t)X/w]\,G\) with non-negative powers of \(w\). Under \(\Plethmap\), \(X+(q-1)w\) becomes \(X'\) and \(qtw\)
becomes \(qt^2y\). For a series \(K(w)\), the part with non-negative powers of \(w\), at \(w=c\), is
\([u^0]\,K(u)/(1-c/u)\). Since \((1-t)/(1-t^{-1})=-t\), this gives, after \(u\mapsto tu\),
\[
 \Plethmap\LL G=[u^0]\,\frac{\Exp[-(X+(q-1)ty)/u]}{1-qty/u}\,G[X';tu].
\]
With the corrected conventions above, the computation in the proof of
\cite[Theorem~3.4.2]{CGM} gives
\[
 \zone F=[u^0]\,F[X+(q-1)ty-(q-1)u;\,u]\,\Exp[(ty-X)/u].
\]
For \(F=\Plethmap G\), \(F[X+(q-1)ty-(q-1)u;u]=G[X';tu]\), and the
two kernels agree, since \(\Exp[(1-q)ty/u]=(1-qty/u)/(1-ty/u)\).

Finally, \(\Plethmap\) is a ring homomorphism with
inverse \(F\mapsto F[(1-t^{-1})X-(q-1)w/t;\,y=w/t]\), \(\Plethmap(z^{-1}G)=ty\,\Plethmap(G)\), and for \(G\in\La\),
\(\Plethmap(G)=\dplus\bigl(G[X/(1-t^{-1})]\bigr)\), since \(\dplus F=F[X+(q-1)y]\) on \(V_0\) \cite[Proposition~3.4.1]{CGM}.
\end{proof}

The first numerator term in \cite[Proposition~3.4.1]{CGM} is printed with
\(y_{i+1}\) in place of \(y_i\), as in \cite[(2.10)]{GGS}.
The printed formula gives \(T_i1=q\), hence
\((T_i-1)(T_i+q)1=2q(q-1)\), contrary to the defining Hecke relation.
The corrected formula gives \(T_i1=1\).
The proof of \cite[Theorem~3.4.2]{CGM} also prints the spectral shift as
\((1-q)ty_1\) instead of \((q-1)ty_1\), and the commutator coefficient as
\(1-q^{-1}\) instead of \(q^{-1}-1\). With the printed signs, the relation
\(\zone(q\dplus\dminus-\dminus\dplus)=qt(\dplus\dminus-\dminus\dplus)\zone\)
of \cite[Definition~3.2.4]{CGM} fails on \(V_1\).

\subsection{Macdonald polynomials in superspace}\label{sec:superspace}

In rank one, pairs consisting of a partition \(\lambda\) and an addable corner \(s\)
correspond bijectively to superpartitions of fermionic degree one.
Under the isometry of Proposition~\ref{prop:superspace}, the eigenvectors of \(\LL\)
correspond to Macdonald polynomials in superspace \cite{BDLM1,BDLM2,Comp}. Throughout
this subsection \(m=1\) and \(v_1=1\), as in Section~\ref{sec:rankone}.

Let \(\theta_1,\theta_2,\dots\) be anticommuting variables. The functions of degree one in the
\(\theta_i\) that are invariant under simultaneous permutations of the pairs \((x_i,\theta_i)\) form the
space \(\bigoplus_{k\geq0}\tilde p_k\La\), where \(\tilde p_k=\sum_i\theta_ix_i^k\). A superpartition of
fermionic degree one is a pair \((k;\nu)\) of an integer \(k\geq0\) and a partition \(\nu\). Put
\(\lambda=\nu\cup(k)\) and \(\lambda+s=\nu\cup(k+1)\). Then \(s\) is an addable corner of \(\lambda\) in
column \(k\), and \((k;\nu)\mapsto(\lambda,s)\) is a bijection onto the pairs of a partition and an
addable corner.

Let \(P_{(k;\nu)}\) be the Macdonald superpolynomials of fermionic degree one, orthogonal for the
pairing in which \(\tilde p_kf\) and \(\tilde p_lg\) pair to \(\Mq^k\delta_{kl}\ip{f}{g}_{\Mq,\Mt}\)
\cite[Section~4.1]{Comp}. Their integral form is \(J_{(k;\nu)}=c_{(k;\nu)}P_{(k;\nu)}\), with
\[
 c_{(k;\nu)}=\prod_{b\in\lambda}\bigl(1-\Mq^{a_{\lambda+s}(b)}\Mt^{l_\lambda(b)+1}\bigr),
\]
where \(a_{\lambda+s}(b)\) is the arm of \(b\) in \(\lambda+s\) and \(l_\lambda(b)\) is its leg in
\(\lambda\). The coefficient of \(\theta_1\) in \(J_{(k;\nu)}\) is the \(1\)-symmetric Macdonald polynomial
of \cite{Lapointe,Goodberry,CL} \cite[Remark~6]{Lapointe}, a polynomial in \(x_1\) with coefficients
symmetric in \(x_2,x_3,\dots\).

\begin{proposition}[The superspace isometry]\label{prop:superspace}
The map \(\Supmap\colon\La[z^{-1}]\to\bigoplus_k\tilde p_k\La\), \(z^{-k}f\mapsto\tilde p_kf\), is an
isometry for the pairing \(\ipinf{\cdot}{\cdot}\) of Section~\ref{sec:rankone} and the superspace pairing.
\begin{enumerate}
\item[(i)] Let \(\tau_i\) be the substitution \(x_i\mapsto\Mq x_i\), and let
\(A_j=\prod_{l\neq j}(\Mt^{-1}x_j-x_l)/(x_j-x_l)\). In \(N\) variables \(x_1,\dots,x_N\),
\(\Supmap\,\LL=E_{N,1}\,\Supmap\), where
\[
 E_{N,1}=(1-\Mt)\Mt^{N-1}\sum_{i,j}\tau_i^{-1}\,\frac{x_jA_j}{x_j-\Mt x_i}\,\theta_i\partial_{\theta_j}
\]
is the operator \(\frac{\Mq}{\Mq-1}(D_{N,1}-D_{N,2})\) of \cite[Section~4.5]{Comp}, from \cite{BDLM2}.
\item[(ii)] If \((k;\nu)\) corresponds to \((\lambda,s)\), then
\(\Supmap\bigl(\PsiNS{\lambda}{s}\bigr)=\Mq^{-k}J_{(k;\nu)}\).
\end{enumerate}
\end{proposition}

\begin{proof}
The isometry is the formula \(\ipinf{z^{-k}f}{z^{-l}g}=\Mq^k\delta_{kl}\ip{f}{g}_{\Mq,\Mt}\).

(i) For \(g\in\La[z^{-1}]\), the coefficient of \(\theta_i\) in \(\Supmap(g)\) is \(g\) at \(z^{-1}=x_i\). In
\(N\) variables \(\Qcr(z)=\prod_l(z^{-1}-\Mt x_l)/(z^{-1}-x_l)\) is a rational function of \(z^{-1}\), with
residue \((1-\Mt)\Mt^{N-1}x_jA_j\) at \(z^{-1}=x_j\), and \(\pile(\Qcr g)\) is the polynomial part of
\(\Qcr g\) in \(z^{-1}\). Since \(\Qcr\) vanishes at \(z^{-1}=\Mt x_i\), partial fractions give
\[
 \pile(\Qcr g)\big|_{z^{-1}=\Mt x_i}
 =(1-\Mt)\Mt^{N-1}\sum_j\frac{x_jA_j}{x_j-\Mt x_i}\;g\big|_{z^{-1}=x_j}.
\]
The operator \(\Twist{\Mt/\Mq}\) replaces \(z^{-1}\) by \((\Mt/\Mq)z^{-1}\), and \(\Qan\) replaces \(p_n\)
by \(p_n+(\Mq^{-n}-1)z^{-n}\), which at \(z^{-1}=x_i\) is \(\tau_i^{-1}p_n\). So
\((\LL g)|_{z^{-1}=x_i}=\tau_i^{-1}\bigl(\pile(\Qcr g)|_{z^{-1}=\Mt x_i}\bigr)\), which is the coefficient
of \(\theta_i\) in \(E_{N,1}\Supmap(g)\).

(ii) By \cite[Section~4.5]{Comp}, \(P_{(k;\nu)}\) is an eigenvector of \(E_{N,1}\) with eigenvalue
\(\ceig{s}\). Take \(N\) larger than \(|\lambda|\). By (i), \(\Supmap^{-1}J_{(k;\nu)}\) is an eigenvector of
\(\LL\) of degree \(|\lambda|\) with eigenvalue \(\ceig{s}\). Different boxes have different contents, so
by \cite[Proposition~3.11]{M26}, \(\Supmap^{-1}J_{(k;\nu)}\) is a combination
\(\sum_\mu c_\mu\PsiNS{\mu}{s}\) over the \(\mu\) with
\(|\mu|=|\lambda|\) and \(s\in\Add{\mu}\). Setting \(x_1=0\) in the coefficient of \(\theta_1\) in
\(\Supmap(g)\) gives \(\piz g\), and \(\piz\PsiNS{\mu}{s}=J_\mu\). So \(\sum_\mu c_\mu J_\mu=\Mq^kJ_\lambda\)
by \cite[Proposition~32]{Lapointe}, and \(c_\mu=\Mq^k\delta_{\mu\lambda}\).
\end{proof}

Under the identification and normalization of Proposition~\ref{prop:superspace}, the eigenbasis
and orthogonality of \cite[Propositions~3.11 and~3.13]{M26} correspond to the
results of \cite{BDLM2}. The constant-term identity \(\piz\PsiNS{\lambda}{s}=J_\lambda\)
agrees with \cite[Proposition~32]{Lapointe}, and the decomposition of \(J_\lambda\)
in \cite[Proposition~3.12]{M26} agrees with \cite[Theorem~6]{CL}.

With the same normalization, the norm formula \cite[Lemma~3.16]{M26} agrees with
\cite[Proposition~28]{GL} and \cite[Theorem~19]{CL}.
The principal-specialization formula \cite[Proposition~4.2]{M26}
agrees with \cite[Theorem~15]{GL} and \cite[Proposition~7]{CL}.
Compare \cite[Theorem~4.13]{Goodberry} for the creation formula
\cite[Proposition~4.8]{M26}.

We know of no earlier
statement of Proposition~\ref{prop:superspace}, or of the form \(\Qan\Twist{\Mt/\Mq}\pile\Qcr\) of the operators
\(E_{N,1}\).

\subsection*{Acknowledgements}
The author thanks Jean-Emile Bourgine for many helpful discussions
relating to this work. The author used Claude Opus 5.5 for much of the
investigation into the higher rank case and its relation to \Negut's work and
the \(q\)-Miura operator.

\appendix
\section{The shuffle commutator}\label{sec:shuffleproof}

Let \(\Shuf\) be the shuffle algebra of Section~\ref{sec:raising}, with the
elements \(\mathsf P_n\) of Lemma~\ref{lem:shuffle}. Its product is denoted by
\(*\), and \(z^0\) is the degree-one shuffle element with constant value \(1\).

\begin{proof}[Proof of Lemma~\ref{lem:shuffle}(iii)]
Put \(\mathsf g(u,v):=(1-q_3v/u)^{-1}\), so that \(\mathsf g(u,0)=1\).
In the variables \(u_l=z_{n-l}\), part~(ii) gives
\(\mathsf P_n=\Sym[K_n\,\Pi(u_0,\dots,u_{n-1})]\), where
\[
 K_n:=u_0^{-1}\prod_{l=1}^{n-1}\mathsf g(u_{l-1},u_l),
 \qquad
 \Pi(u_0,\dots,u_{n-1}):=\prod_{l<l'}\omega(u_{l'}/u_l).
\]
Let \(w\) be one more variable and put \(u_n:=0\). For \(0\leq p\leq n\), set
\[
 \Pi_p:=\Pi(u_0,\dots,u_{p-1},w,u_p,\dots,u_{n-1}),
 \qquad Z_p:=\Sym[K_n\Pi_p].
\]
The symmetrization includes the additional variable \(w\).
Then \(Z_n=z^0*\mathsf P_n\) and
\(Z_0=\mathsf P_n*z^0\). For \(0\leq p<n\), set
\[
 Y_p:=\Sym[K_n\,\mathsf g(u_p,w)\,\mathsf g(w,u_{p+1})\,\Pi_{p+1}],
 \qquad Y_{-1}:=0.
\]
Since \(Y_{n-1}=\mathsf P_{n+1}\), summing the identities
\(Z_{p+1}-Z_p=\epsilon(Y_p-Y_{p-1})\) proves the commutator.

Fix \(p\), write \(u:=u_p\) and \(u^+:=u_{p+1}\), and write
\(u^-:=u_{p-1}\) when \(p\geq1\). The identity
\(\omega(x)^{-1}=1+\epsilon x/((1-x)(1-q_3x))\), applied at \(w/u\) and \(u/w\), gives
\[
 \omega(w/u)-\omega(u/w)=\frac{\epsilon}{w-u}\bigl(w\,\omega(w/u)\,\mathsf g(u,w)
 +u\,\omega(u/w)\,\mathsf g(w,u)\bigr).
\]
Since \(\Pi_{p+1}/\omega(w/u)=\Pi_p/\omega(u/w)\), exchanging \(u\) and \(w\)
in the second term gives
\[
 Z_{p+1}-Z_p=\epsilon\,\Sym\Bigl[\frac{w}{w-u}\,(K_n-K_n^w)\,\mathsf g(u,w)\,\Pi_{p+1}\Bigr],
\]
where \(K_n^w\) is obtained from \(K_n\) by replacing \(u_p\) with \(w\).
For any \(a\),
\[
 \frac{w}{w-u}\bigl(\mathsf g(a,u)\mathsf g(u,u^+)-\mathsf g(a,w)\mathsf g(w,u^+)\bigr)
 =\mathsf g(a,u)\,\mathsf g(w,u^+)\bigl(\mathsf g(u,u^+)-\mathsf g(a,w)\bigr).
\]
For \(p\geq1\), set \(a=u^-\). For \(p=0\), multiply by \(-q_3/a\) and let
\(a\to0\), using \((-q_3/a)\mathsf g(a,v)\to v^{-1}\) and
\(\mathsf g(a,w)\to0\). After multiplication by the remaining factors
and symmetrization, the term with \(\mathsf g(u,u^+)\) is \(Y_p\).
For \(p\geq1\), the term with \(\mathsf g(u^-,w)\) becomes \(Y_{p-1}\) after
exchanging \(u\) and \(w\). This proves the identities.
The operator commutator follows from part~(ii), since \(z^0\) acts by
\(\xplus{0}=\pushp\pulm\).
\end{proof}


\end{document}